\documentclass[10pt,journal,twocolumn]{IEEEtran}

\usepackage[english]{babel}
\usepackage[latin1]{inputenc}
\usepackage[T1]{fontenc}
\usepackage{url}
\usepackage{graphicx}
\usepackage{subfigure}
\usepackage{color}
\usepackage[cmex10]{amsmath}
\usepackage{amsthm,amsfonts,amssymb}
\usepackage{colortbl}                 
\usepackage{url}
\usepackage{booktabs}
\usepackage{epsfig}
\usepackage{pstricks}
\usepackage{pst-node}
\usepackage{pst-grad}
\usepackage{pst-plot}
\usepackage{pstricks-add}
\usepackage{pst-sigsys}
\usepackage{ifthen}
\usepackage{algorithm}
\usepackage{algpseudocode}
\usepackage{float}
\usepackage{fp}
\usepackage{cite}  % For Computer Society papers

\newfloat{algorithm}{t}{lop}

\newboolean{draft}
\newcommand{\isdraft}[2]{\ifthenelse{\boolean{draft}}{#1}{#2}}

\isdraft{\usepackage{setspace}}{}                              % DRAFT
\isdraft{\usepackage[footnotesize]{caption}}{}                 % DRAFT
\isdraft{\usepackage{paralist}}{}                              % DRAFT

\def\today{\ifcase\month\or
  January\or February\or March\or April\or May\or June\or
  July\or August\or September\or October\or November\or December\fi \space \number\year}
\newcommand{\fref}[1]{Fig.~\ref{#1}}

\newcommand{\cref}[1]{Chapter~\ref{#1}}
\newcommand{\sref}[1]{Section~\ref{#1}}
\newcommand{\ssref}[1]{Subsection~\ref{#1}}
\newcommand{\aref}[1]{Appendix~\ref{#1}}

\newcommand{\thref}[1]{Theorem~\ref{#1}}

\newcommand{\lref}[1]{Lemma~\ref{#1}}
\newcommand{\coref}[1]{Corollary~\ref{#1}}
\newcommand{\Arxiv}[1]{preprint: \url{#1}}
\newcommand{\OptimizationOnline}[1]{preprint: \url{#1}}
\theoremstyle{plain}
\newtheorem{Theorem}{Theorem}

\newtheorem{Corollary}[Theorem]{Corollary}

\newtheorem{Lemma}[Theorem]{Lemma}

\definecolor{red}{RGB}{153,0,0}
\definecolor{green}{RGB}{0,153,0}			
\definecolor{blue}{RGB}{0,0,153}
\definecolor{darkred}{RGB}{90,0,0}
\definecolor{darkgreen}{RGB}{0,90,0}
\definecolor{darkblue}{RGB}{0,0,90}			

\newcommand{\mypar}[1]{{\bf #1.}}

\isdraft{}{}

\title{Adaptive-Rate Sparse Signal Reconstruction With Application in Compressive Background Subtraction}

\author{Jo\~ao F.\ C.\ Mota, Nikos Deligiannis, Aswin C.\ Sankaranarayanan, Volkan Cevher, Miguel R.\ D.\ Rodrigues
\IEEEcompsocitemizethanks{
	\IEEEcompsocthanksitem
		J. Mota, N. Deligiannis, and M. Rodrigues were supported by the EPSRC grant EP/K033166/1. A.\ C.\ Sankaranarayanan was supported in part by the NSF grant CCF-1117939. V.\ Cevher's work is supported in part by the European Commission under grants MIRG-268398 and ERC Future Proof and by the Swiss Science Foundation under grants SNF 200021-132548, SNF 200021-146750, and SNF CRSII2-147633.	Part of this work will be presented at the IEEE International Conference on Acoustics, Speech, and Signal Processing (ICASSP), Brisbane, 2015, \cite{Mota14-DynamicSparseStateEstimationUsingL1L1}.
	\IEEEcompsocthanksitem 
		J. Mota, N. Deligiannis, and M. Rodrigues are with the Department of Electronic and Electrical Engineering, University College London, UK.\protect\\
		E-mail: \{j.mota,n.deligiannis,m.rodrigues\}@ucl.ac.uk.
	\IEEEcompsocthanksitem 
		A. Sankaranarayanan is with the Department of Electrical and Computer Engineering, Carnegie Mellon University, PA 15213, USA.\protect\\
		E-mail: saswin@andrew.cmu.edu.
	\IEEEcompsocthanksitem 
		V. Cevher is with the Laboratory for Information and Inference Systems (LIONS), EPFL, Switzerland.\protect\\
		E-mail: volkan.cevher@epfl.ch.
}% <-this % stops an unwanted space
\thanks{}
}

\begin{document}

\maketitle
%\isdraft{\pagebreak}{}

\begin{abstract}
We propose and analyze an online algorithm for reconstructing a sequence of signals from a limited number of linear measurements. The signals are assumed sparse, with unknown support, and evolve over time according to a generic nonlinear dynamical model. Our algorithm, based on recent theoretical results for $\ell_1$-$\ell_1$ minimization, is recursive and computes the number of measurements to be taken at each time on-the-fly. As an example, we apply the algorithm to compressive video background subtraction, a problem that can be stated as follows: given a set of measurements of a sequence of images with a static background, simultaneously reconstruct each image while separating its foreground from the background. The performance of our method is illustrated on sequences of real images: we observe that it allows a dramatic reduction in the number of measurements with respect to state-of-the-art compressive background subtraction schemes. 
%We propose and analyze a novel method for compressive background subtraction that reconstructs the foreground of each frame from a small set of linear measurements. Given the background image, our method has two main components: a motion estimation module, which estimates the current frame from previous ones, and an optimization problem, namely $\ell_1$-$\ell_1$ minimization, that reconstructs sparse signals using side information. The use of side information in the reconstruction scheme enables reducing the number of required measurements dramatically. We model our problem as the estimation of the state of a dynamical system and tackle it using recent theoretical results for $\ell_1$-$\ell_1$ minimization. The performance of our method is illustrated on two sequences of real images, and we observe that it requires up to four times less measurements per frame than state-of-the-art compressive background subtraction schemes. 
\end{abstract}


\section{Introduction}
\label{Sec:Intro}

	%\IEEEPARstart{E}{stimating} a sequence of sparse, time-varying signals from a limited number of measurements arises in applications as diverse as object tracking in computer vision systems~\cite{Cevher08-CompressiveSensingForBackgroundSubtraction}, localization in wireless networks~\cite{Ribeiro10-KalmanFilterInWSN}, and statistical inference estimation in wireless sensor networks~\cite{Barbarossa}. In this paper, we consider 

	%\IEEEPARstart{E}{stimating} the state of a dynamical system from a limited number of observations is a central problem in control theory, finding numerous applications in object tracking in computer vision systems~\cite{Chellappa09-StatisticalMethodsAndModelsForVideoBasedTrackingModelingAndRecognition}, localization in wireless networks~\cite{Ribeiro10-KalmanFilterInWSN}, and statistical inference estimation in wireless sensor networks~\cite{Barbarossa}.  
	%~\cite{Sontag98-MathematicalControlTheory,Kailath00-LinearEstimation} and references therein. 	
	
	%Consider the estimation of the state of a dynamical system from a limited number of state measurements: 
	%[Paragraph with high-level motivation]
	
	\IEEEPARstart{C}{onsider} the problem of reconstructing a sequence of sparse signals from a limited number of measurements. Let~$x[k] \in \mathbb{R}^n$ be the signal at time~$k$ and $y[k] \in \mathbb{R}^{m_k}$ be the vector of signal measurements at time~$k$, where~$m_k \ll n$. Assume the signals evolve according to the dynamical model	
	\begin{subequations}\label{Eq:IntroStateSpace}
  \begin{align}
		x[k] &= f_k\big(\{x[i]\}_{i=1}^{k-1}\big) + \epsilon[k]
			\label{Eq:IntroStateSpaceModel}
    \\
    y[k] &= A_k\,x[k]\,,
    \label{Eq:IntroStateSpaceObservations}
  \end{align}
  \end{subequations}  
  where~$\epsilon[k] \in \mathbb{R}^n$ is modeling noise and~$A_k \in \mathbb{R}^{m_k \times n}$ is a sensing matrix. In~\eqref{Eq:IntroStateSpaceModel}, $f_k:(\mathbb{R}^n)^{k-1} \xrightarrow{} \mathbb{R}^n$ is a known, but otherwise arbitrary, map that describes~$x[k]$ as a function of past signals. We assume that each~$x[k]$ and~$\epsilon[k]$ is sparse, i.e., it has a small number of nonzero entries. Our goal is to reconstruct the signal sequence~$\{x[k]\}$ from the measurement sequence~$\{y[k]\}$. We require the reconstruction scheme to be recursive (or	 online), i.e., $x[k]$ is reconstructed before acquiring measurements of any future signal~$x[i]$, $i>k$, and also to use a minimal number of measurements. We formalize the problem as follows.
  
  \mypar{Problem statement}
  \textit{Given two unknown sparse sequences~$\{x[k]\}$ and~$\{\epsilon[k]\}$ satisfying~\eqref{Eq:IntroStateSpace}, design an online algorithm that \text{1)} uses a minimal number of measurements~$m_k$ at time~$k$, and \text{2)} perfectly reconstructs each~$x[k]$ from~$y[k]$ acquired as in~\eqref{Eq:IntroStateSpaceObservations}, and possibly $x[i]$, $i<k$.} 
  
  Note that our setting immediately generalizes from the case where each~$x[k]$ is sparse to the case where~$x[k]$ has a sparse representation in a linear, invertible transform.\footnote{If~$x[k]$ is not sparse but $z[k] := \Psi x[k]$ is, where~$\Psi$ is an invertible matrix, then redefine~$f_k$ as the composition $f_k^z = \Psi^{-1} \circ f_k \circ \Psi$ and~$A_k$ as $A_k^z:= A_k \Psi^{-1}$. The signal~$z[k]$ satisfies~\eqref{Eq:IntroStateSpace} with~$f_k^z$ and~$A_k^z$.} 
  %Our motivation for considering the above problem, however, stems from its application to compressive video background subtraction~\cite{Cevher08-CompressiveSensingForBackgroundSubtraction}, as explained next.
  
  %Although we also assume that the noise~$\epsilon[k]$ is Laplacian, and thus that each of its realizations is also a sparse vector, our algorithm and results can be easily adapted to the case where the noise is Gaussian.\footnote{In addition, they can also be adapted to the case where there is additive measurement noise (with bounded magnitude) in~\eqref{Eq:IntroStateSpaceObservations}. This case is, however, significantly more complicated.} Our motivation for considering~\eqref{Eq:IntroStateSpace} with Laplacian noise is its application to compressive video background subtraction~\cite{Cevher08-CompressiveSensingForBackgroundSubtraction}, as explained next.
  
	%We denote the unknown target signal at time~$k$ with~$x[k] \in \mathbb{R}^n$, from which we have access to~$m_k$ linear measurements $y[k] = A_k x[k]$, where~$A_k \in \mathbb{R}^{m_k \times n}$ is a measurement matrix. We assume~$x[k]$ is a function of past signals plus sparse noise, i.e., $x[k] = f_k\big(\{x[i]\}_{i=1}^{k-1}\big) + \epsilon[k]$, where~$f_k$ is a known function and~$\epsilon[k]$ is an unknown sparse vector. That is, the sequence~$\{x[k]\}$ evolves according to 
	
	\subsection{Applications}

	Many problems require estimating a sequence of signals from a sequence of measurements satisfying the model in~\eqref{Eq:IntroStateSpace}. These include classification and tracking in computer vision
	systems~\cite{Forsyth02-ComputerVision-AModernApproach,Chellappa09-StatisticalMethodsAndModelsForVideoBasedTrackingModelingAndRecognition},
	%~\cite{Lipton98-MovingTargetClassificationAndTrackingFromRealTimeVideo,Forsyth02-ComputerVision-AModernApproach,Chellappa09-StatisticalMethodsAndModelsForVideoBasedTrackingModelingAndRecognition},
	radar 
	%tracking~\cite{Pearson74-KalmanFilterApplicationsInAirborneRadarTracking,Herman09-HighResolutionRadarCS},
	tracking~\cite{Herman09-HighResolutionRadarCS},
	dynamic 	
	MRI~\cite{Eldar14-ApplicationCSLongitudinalMRI}
	%~\cite{Lustig07-SparseMRI,Chen08-PriorImageConstrainedCS,Eldar14-ApplicationCSLongitudinalMRI} 
	and several tasks in wireless sensor 
	networks~\cite{Ribeiro10-KalmanFilterInWSN}.
	%~\cite{Barbarossa07-DecentralizedMLEForSensorNetworksComposedOfNonlinearlyCoupledDynamicalSystems,Ribeiro10-KalmanFilterInWSN}. 
	
	Our application focus, however, is \textit{compressive background subtraction}~\cite{Cevher08-CompressiveSensingForBackgroundSubtraction}. Background subtraction is a key task for detecting and tracking objects in a video sequence and it has been applied, for example, in video 
	surveillance
	%~\cite{Gutchess01-BackgroundModelInitializationAlgorithmsVideoSurveillance,Maddalena08-SelfOrganizingApproachBackgroundSubtractionVisualSurveillanceApplications,Brutzer11-EvaluationBackgroundSubtractionTechniquesVideoSurveillance}, 
	\cite{Maddalena08-SelfOrganizingApproachBackgroundSubtractionVisualSurveillanceApplications,Brutzer11-EvaluationBackgroundSubtractionTechniquesVideoSurveillance}, 
	traffic monitoring~\cite{Tseng02-RealTimeVideoSurveillanceForTrafficMonitoringUsingVirtualLineAnalysis,Cheung03-RobustTechniquesForBackgroundSubtractionUrbanTrafficVideo}, and medical imaging~\cite{Profio86-DigitalBackgroundSubtractionForFluorescenceImaging,Otazo14-LowRankPlusSparseMatrixDecompositionAcceleratedDynamicMRI}. Although there are many background subtraction techniques, e.g., 
	%\cite{Cucchiara03-DetectingMovingObjectsGhostsAndShadowsInVideoStreams,Picardi04-BackgroundSubtractionTechniques-AReview,Lee05-EffectiveGaussianMixtureLearningForVideoBackgroundSubtraction,Chellappa09-StatisticalMethodsAndModelsForVideoBasedTrackingModelingAndRecognition,Candes11-RobustPrincipalComponentAnalysis}
	\cite{Picardi04-BackgroundSubtractionTechniques-AReview,Chellappa09-StatisticalMethodsAndModelsForVideoBasedTrackingModelingAndRecognition,Candes11-RobustPrincipalComponentAnalysis},
	most of them assume access to full frames and, thus, are inapplicable in compressive video %sensing~\cite{Wakin06-CompressiveImagingForVideoRepresentationAndCoding,Park09-MultiscaleFrameworkForCompressiveSensingOfVideo,Reddy11-P2C2,Sankaranarayanan12-CSMUVI,Sankaranarayanan13-CompressiveAcquisitionForLinearDynamicalSystems},
	sensing~\cite{Wakin06-CompressiveImagingForVideoRepresentationAndCoding,Sankaranarayanan12-CSMUVI,Sankaranarayanan13-CompressiveAcquisitionForLinearDynamicalSystems},
	a technology used in cameras where sensing is expensive (e.g., infrared, UV wavelengths).
	
	In compressive video sensing, one has access not to full frames as in conventional video, but only to a small set of linear measurements of each frame, as in~\eqref{Eq:IntroStateSpaceObservations}. Cevher et al.~\cite{Cevher08-CompressiveSensingForBackgroundSubtraction} noticed that background subtraction is possible in this context if the foreground pixels, i.e., those associated to a moving object, occupy a small area in each frame. Assuming the background image is known beforehand, compressed sensing techniques~\cite{Donoho06-CompressedSensing,Candes06-RobustUncertaintyPrinciplesExactSignalReconstructionHighlyIncomplete} such as $\ell_1$-norm minimization allow reconstructing each foreground.
	%using a small number of measurements. 
	This not only reconstructs the original frame (if we add the reconstructed foreground to the known background), but also performs background subtraction as a by-product~\cite{Cevher08-CompressiveSensingForBackgroundSubtraction}. 
	
	We mention that, with the exception of~\cite{Warnell12-AdaptiveRateCompressiveSensingBackgroundSubtraction,Warnell14-AdaptiveRateCompressiveSensingUsingSideInformation}, most approaches to compressive video sensing and to compressive background subtraction assume a fixed number of measurements for all frames~\cite{Wakin06-CompressiveImagingForVideoRepresentationAndCoding,Cevher08-CompressiveSensingForBackgroundSubtraction,Park09-MultiscaleFrameworkForCompressiveSensingOfVideo,Reddy11-P2C2,Sankaranarayanan12-CSMUVI,Sankaranarayanan13-CompressiveAcquisitionForLinearDynamicalSystems}. If this number is too small, reconstruction of the frames fails. If it is too large, reconstruction succeeds, but at the cost of spending unnecessary measurements in some or all frames. The work in~\cite{Warnell12-AdaptiveRateCompressiveSensingBackgroundSubtraction,Warnell14-AdaptiveRateCompressiveSensingUsingSideInformation} addresses this problem with an online scheme that uses cross validation to compute the number of required measurements. Given a reconstructed foreground, \cite{Warnell12-AdaptiveRateCompressiveSensingBackgroundSubtraction,Warnell14-AdaptiveRateCompressiveSensingUsingSideInformation} estimates the area of the true foreground using extra cross-validation measurements. Then, assuming that foreground areas of two consecutive frames are the same, the phase diagram of the sensing matrix, which was computed beforehand, gives the number of measurements for the next frame. This approach, however, fails to use information from past frames in the reconstruction process, information that, as we will see, can be used to significantly reduce the number of measurements.

  \subsection{Overview of our approach and contributions}
  \label{SubSubSec:OverviewAndContributions}
  
  \mypar{Overview}
  Our approach to adaptive-rate signal reconstruction is based on the recent theoretical results of~\cite{Mota14-CSSideInfo,Mota14-CSwithSideInfo-GlobalSIP}. These characterize the performance of sparse reconstructing schemes in the presence of side information. The scheme we are most interested in is the \textit{$\ell_1$-$\ell_1$ minimization}:
  \begin{equation}\label{Eq:L1L1}
		\begin{array}[t]{ll}
			\underset{x}{\text{minimize}} & \|x\|_1 + \beta \|x - w\|_1 \\
      \text{subject to} & Ax = y\,,
    \end{array}
  \end{equation}
  where~$x \in \mathbb{R}^n$ is the optimization variable and~$\|x\|_1:= \sum_{i=1}^n|x_i|$ is the $\ell_1$-norm. In~\eqref{Eq:L1L1}, $y \in \mathbb{R}^m$ is a vector of measurements and~$\beta$ is a positive parameter. The vector~$w \in \mathbb{R}^n$ is assumed known and is the so-called \textit{prior} or \textit{side information}: a vector similar to the vector that we want to reconstruct, say~$x^\star$. Note that if we set~$\beta = 0$ in \eqref{Eq:L1L1}, we obtain \textit{basis pursuit}~\cite{Donoho98-AtomicDecompositionBasisPursuit}, a well-known problem for reconstructing sparse signals and which is at the core of the theory of \textit{compressed sensing}~\cite{Candes06-RobustUncertaintyPrinciplesExactSignalReconstructionHighlyIncomplete,Donoho06-CompressedSensing}. Problem~\eqref{Eq:L1L1} generalizes basis pursuit by integrating the side information~$w$. The work in~\cite{Mota14-CSSideInfo,Mota14-CSwithSideInfo-GlobalSIP} shows that, if~$w$ has reasonable quality and the entries of~$A$ are drawn from an i.i.d.\ Gaussian distribution, the number of measurements required by~\eqref{Eq:L1L1} to reconstruct~$x^\star$ is much smaller than the number of measurements required by basis pursuit. Furthermore, the theory in~\cite{Mota14-CSSideInfo,Mota14-CSwithSideInfo-GlobalSIP} establishes that~$\beta = 1$ is an optimal choice, irrespective of any problem parameter. This makes the reconstruction problem~\eqref{Eq:L1L1} parameter-free. 
    
	We address the problem of recursively reconstructing a sequence of sparse signals satisfying~\eqref{Eq:IntroStateSpace} as follows.
  %\mypar{Contributions}
  Assuming the measurement matrix is Gaussian,\footnote{Although Gaussian matrices are hard to implement in practical systems, they have optimal performance. There are, however, other more practical matrices with a similar performance, e.g., \cite{Berinde08-SparseRecoveryUsingSparseMatrices,Liutkus14-ImagingWithNature-CompressiveImagingUsingMultiplyScatteringMedium}.} we propose an algorithm that uses~\eqref{Eq:L1L1} with~$w = f_k\big(\{x[i]\}_{i=1}^{k-1}\big)$ to reconstruct each signal~$x[k]$. And, building upon the results of~\cite{Mota14-CSSideInfo,Mota14-CSwithSideInfo-GlobalSIP}, we equip our algorithm with a mechanism to automatically compute an estimate on the number of required measurements. As application, we consider compressive background subtraction and show how to generate side information from past frames.
  
  \mypar{Contributions}
  We summarize our contributions as follows:
  \newcounter{Contributions}
	\begin{list}{\roman{Contributions})}{\usecounter{Contributions}}
  	\item 
			We propose an adaptive-rate algorithm for reconstructing sparse sequences satisfying the model in~\eqref{Eq:IntroStateSpace}.
		\item 
			We establish conditions under which our algorithm reconstructs a finite sparse sequence~$\{x[i]\}_{i=1}^k$ with large probability.
		\item
			We describe how to apply the algorithm to compressive background subtraction problems, using motion-compensated extrapolation to predict the next image to be acquired. In other words, we show how to generate side information.
		\item
			Given that images predicted by motion-compensated extrapolation are known to exhibit Laplacian noise, we then characterize the performance of~\eqref{Eq:L1L1} under this model.
		\item
			Finally, we show the impressive performance of our algorithm for performing 
			compressive background subtraction on a sequence of real images.
		%\item 
		%	Finally, we also describe a simple solver for the optimization problem in~\eqref{Eq:L1L1}.
  \end{list}
  
  Besides the incorporation of a scheme to compute a minimal number of measurements on-the-fly, there is another aspect that makes our algorithm fundamentally different from prior work. As overviewed in \sref{Sec:RelatedWork}, most prior algorithms for reconstructing dynamical sparse signals work well only when the sparsity pattern of~$x[k]$ varies slowly with time. Our algorithm, in contrast, operates well even when the sparsity pattern of~$x[k]$ varies arbitrarily between consecutive time instants, as shown by our theory and experiments. What is required to vary slowly is the ``quality'' of the prediction given by each~$f_k$ (i.e., the quality of the side information) and, to a lesser extent, not the sparsity pattern of~$x[k]$ but only its sparsity, i.e., the \textit{number of} nonzero entries.

  \subsection{Organization}
        
  \sref{Sec:RelatedWork} overviews related work. In \sref{Sec:Preliminaries}, we state the results from~\cite{Mota14-CSSideInfo,Mota14-CSwithSideInfo-GlobalSIP} that are used by our algorithm. \sref{Sec:DynamicalSignalReconstruction} describes the algorithm and establishes reconstruction guarantees. \sref{Sec:ApplicationToBackgroundSubtraction} concerns the application to compressive background subtraction. Experimental results illustrating the performance of our algorithm are shown in section~\ref{Sec:ExperimentalResults}; and section~\ref{Sec:Conclusions} concludes the paper. The appendix
  contains the proofs of our results.
  %contains the description of a simple solver for~\eqref{Eq:L1L1} and the proofs of some results. 

\section{Related work}
\label{Sec:RelatedWork}
	
	There is an extensive literature on reconstructing time-varying signals from limited measurements. Here, we provide an overview by referring a few landmark papers. 
	
	\mypar{The Kalman filter}
	The classical solution to estimate a sequence of signals satisfying~\eqref{Eq:IntroStateSpace} or, in the control terminology, the state of a dynamical system, is
  the Kalman filter~\cite{Kalman60-NewApproachLinearFilteringAndPredictionProblems}. The Kalman filter is an online algorithm that is least-squares optimal when the model is linear, i.e., $f_k\big(\{x[i]\}_{i=0}^{k-1}\big) = F_k  x[k]$, and the sequence~$\{\epsilon[k]\}$ is Gaussian and independent across time. Several extensions are available when these assumptions do not hold~\cite{Haykin01-KalmanFilteringAndNeuralNetworks,Geromel99-OptimalLinearFilteringUnderParameterUncertainty,Ghaoui01-RobustFilteringDiscreteTimeSystemsBoundedNoiseParametricUncertainty}. The Kalman filter and its extensions, however, are inapplicable to our scenario, as they do not easily integrate the additional knowledge that the state is sparse.
  %; in addition, the state in our scenario is not observable, because the number of measurements is very limited, i.e., $m_k \ll n$. 
  
  \mypar{Dynamical sparse signal reconstruction}
  Some prior work incorporates signal structure, such as sparsity, into online sparse reconstruction procedures. For example, \cite{Vaswani08-KalmanFilteredCS,Vaswani09-AnalyzingKalmanFilteredCS} adapts a Kalman filter to estimate a sequence of sparse signals. Roughly, we have an estimate of the signal's support at each time instant and use the Kalman filter to compute the (nonzero) signal values. When a change in the support is detected, the estimate of the support is updated using compressed sensing techniques. The work in~\cite{Vaswani08-KalmanFilteredCS,Vaswani09-AnalyzingKalmanFilteredCS}, however, assumes that the support varies very slowly and does not provide any strategy to update (or compute) the number of measurements; indeed, the number of measurements is assumed constant along time. Related work that also assumes a fixed number of measurements includes~\cite{Ziniel10-TrackingAndSmoothingTimeVaryingSparseSignalsBP}, which uses approximate belief propagation, and
  %\cite{Carmi10-MethodsForSparseSignalRecoveryUsingKalmanFiltering,Kanevsky10-KalmanFilteringForCompressedSensing}, 
	\cite{Carmi10-MethodsForSparseSignalRecoveryUsingKalmanFiltering}, 
	which integrates sparsity knowledge into a Kalman filter via a pseudo-measurement technique. The works in~\cite{Balzano10-OnlineIdentificationAndTrackingOfSubspacesFromHighlyIncompleteInformation,Balzano14-LocalConvergenceOfAnAlgorithmForSubspaceIdentificationFromPartialData-GROUSE} and~\cite{Chi13-PETRELS} propose online algorithms named GROUSE and PETRELS, respectively, for estimating signals that lie on a low-dimensional subspace. Their model can be seen as a particular case of~\eqref{Eq:IntroStateSpace}, where each map~$f_k$ is linear and depends only on the previous signal. Akin to most prior work, both GROUSE and PETRELS assume that the rank of the underlying subspace (i.e., the sparsity of~$x[k]$) varies slowly with time, and fail to provide a scheme to compute the number of measurements.
        
  We briefly overview the work in~\cite{Charles11-SparsityPenaltiesDynamicalSystemEstimation}, which is
  probably the closest to ours.
  %, although it does not provide any mechanism to compute the number of measurements. 
  Three dynamical reconstruction schemes are studied
  in~\cite{Charles11-SparsityPenaltiesDynamicalSystemEstimation}. The one with the best performance is 
  \begin{equation}\label{Eq:RombergProb}
  	\underset{x}{\text{minimize}}\,\,\|x\|_1 + \beta\|x - w\|_1 + \beta_2\|Ax - y\|_2^2\,,
  \end{equation} 
  where~$\beta_2 > 0$ and~$\|\cdot\|_2$ is the Euclidean $\ell_2$-norm. Problem~\eqref{Eq:RombergProb} is the Lagrangian version of the problem we obtain by replacing the constraints of~\eqref{Eq:L1L1} with~$\|Ax - y\|_2 \leq \sigma$, where~$\sigma$ is a bound on the measurement noise; see problem~\eqref{Eq:L1L1Noisy} below. For~$\beta_2$ in a given range, the solutions of~\eqref{Eq:L1L1Noisy} and~\eqref{Eq:RombergProb} coincide. This is why the approach in~\cite{Charles11-SparsityPenaltiesDynamicalSystemEstimation} is so closely related to ours. Nevertheless, using~\eqref{Eq:L1L1Noisy} has two important advantages: first, in practice, it is easier to obtain bounds on the measurement noise~$\sigma$ than it is to tune~$\beta_2$; second, and more importantly, the problem in~\eqref{Eq:L1L1Noisy} has well-characterized reconstruction guarantees~\cite{Mota14-CSSideInfo,Mota14-CSwithSideInfo-GlobalSIP}. It is exactly those guarantees that enable our scheme for computing of the number of measurements online. The work in~\cite{Charles11-SparsityPenaltiesDynamicalSystemEstimation}, as most prior work, assumes a fixed number of measurements for all signals.

\section{Preliminaries: Static signal \\reconstruction using    $\ell_1$-$\ell_1$ minimization}
\label{Sec:Preliminaries}
        
  This section reviews some results from~\cite{Mota14-CSSideInfo}, namely reconstruction guarantees
  for~\eqref{Eq:L1L1} in a static scenario, i.e., when we estimate just one signal, not a sequence. As mentioned before, $\beta = 1$ is an optimal choice: it not only minimizes the bounds in~\cite{Mota14-CSSideInfo}, but also leads to the best results in practice. This is the reason why we use~$\beta=1$ henceforth.

	\mypar{\boldmath{$\ell_1$}-\boldmath{$\ell_1$} minimization}	
  Let~$x^\star \in \mathbb{R}^{n}$ be a sparse vector, and assume we have~$m$ linear measurements of~$x^\star$: $y = Ax^\star$, where~$A \in \mathbb{R}^{m\times n}$. Denote the \textit{sparsity} of~$x^\star$ with $s := |\{i\,:\, x_i^\star \neq 0\}|$, where~$|\cdot|$ is the cardinality of a set. Assume we have access to a signal~$w \in \mathbb{R}^n$ similar to~$x^\star$ (in the sense that $\|x^\star - w\|_1$ is small) and suppose we attempt to reconstruct~$x^\star$ by solving the $\ell_1$-$\ell_1$ minimization problem~\eqref{Eq:L1L1} with~$\beta = 1$: 
  \begin{equation}\label{Eq:L1L1Simple}
		\begin{array}[t]{ll}
			\underset{x}{\text{minimize}} & \|x\|_1 + \|x - w\|_1 \\
      \text{subject to} & Ax = y\,.
    \end{array}
  \end{equation}
  The number of measurements that problem~\eqref{Eq:L1L1Simple} requires to reconstruct~$x^\star$ is a function of the ``quality'' of the side information~$w$. Quality in~\cite{Mota14-CSSideInfo} is measured in terms of the following parameters:       
  \begin{subequations}\label{Eq:QualityParameters}
  \begin{align}
			\xi &:= \big|\{i\,:\, w_i \neq x_i^\star = 0\}\big| - \big|\{i\,:\, w_i = x_i^\star \neq 0\}\big|\,,        
    \label{Eq:xi}
    \\
			\overline{h} 
    &:= 
			\big|
				\{i\,:\, x_i^\star > 0, \,\,x_i^\star > w_i\} \cup \{i\,:\, x_i^\star < 0, \,\,x_i^\star < w_i\}
      \big|\,.
    \label{Eq:hBar}
  \end{align}
  \end{subequations}
  Note that the number of components of~$w$ that contribute to~$\overline{h}$ are the ones defined on the support of~$x^\star$; thus, $0\leq \overline{h}\leq s$.  
  \begin{Theorem}[Th.\ 1 in \cite{Mota14-CSSideInfo}]
  \label{Thm:L1L1}
		Let~$x^\star, w \in \mathbb{R}^n$ be the vector to reconstruct and the side information, respectively. Assume~$\overline{h} > 0$ and that there exists at least one index~$i$ for which $x_i^\star = w_i = 0$. Let the entries of~$A \in \mathbb{R}^{m \times n}$ be i.i.d.\ Gaussian with zero mean and variance~$1/m$. If
    \begin{equation}\label{Eq:L1L1Bound}
			m \geq 2\overline{h}\log\Big(\frac{n}{s + \xi/2}\Big) + \frac{7}{5}\Big(s + \frac{\xi}{2}\Big) + 1\,,
    \end{equation}
    then, with probability at least $1 - \exp\big(-\frac{1}{2}(m - \sqrt{m})^2\big)$, $x^\star$ is the unique solution of~\eqref{Eq:L1L1Simple}.
  \end{Theorem}
  Theorem~\ref{Thm:L1L1} establishes that if the number of measurements is larger than~\eqref{Eq:L1L1Bound}
  then, with high probability, \eqref{Eq:L1L1Simple} reconstructs~$x^\star$ perfectly. The bound
  in~\eqref{Eq:L1L1Bound} is a function of the signal dimension~$n$ and sparsity~$s$, and of the
  quantities~$\xi$ and~$\overline{h}$, which depend on the signs of the entries of~$x^\star$ and~$w -
  x^\star$, but not on their magnitudes. When~$w$ approximates~$x^\star$ reasonably well, the bound
  in~\eqref{Eq:L1L1Bound} is much smaller than the one for basis pursuit\footnote{Recall that basis pursuit is~\eqref{Eq:L1L1} with~$\beta = 0$.} in~\cite{Chandrasekaran12-ConvexGeometryLinearInverseProblems}:
  \begin{equation}\label{Eq:ChandrasekaranBound}
		m \geq 2s\log\Big(\frac{n}{s}\Big) + \frac{7}{5}s + 1\,.
  \end{equation} 
  Namely, \cite{Chandrasekaran12-ConvexGeometryLinearInverseProblems} establishes that if~\eqref{Eq:ChandrasekaranBound} holds and if~$A \in \mathbb{R}^{m\times n}$ has i.i.d.\ Gaussian entries with zero mean and variance~$1/m$ then, with probability similar to the one in Theorem~\ref{Thm:L1L1}, $x^\star$ is the unique solution to basis pursuit.
  Indeed, if~$\overline{h}\ll s$ and~$\xi$ is larger than a small negative constant,
  then~\eqref{Eq:L1L1Bound} is much smaller than~\eqref{Eq:ChandrasekaranBound}. Note that, in practice,
  the quantities~$s$, $\xi$, and~$\overline{h}$ are unknown, since they depend on the unknown
  signal~$x^\star$. In the next section, we propose an online scheme to estimate them using past signals.
        
	\mypar{Noisy case} Theorem~\ref{Thm:L1L1} has a counterpart for noisy measurements, which we state informally; see~\cite{Mota14-CSSideInfo} for details. Let~$y = Ax^\star + \eta$, where~$\|\eta\|_2 \leq \sigma$. Let also~$A \in \mathbb{R}^{m\times n}$ be as in Theorem~\ref{Thm:L1L1} with
	\begin{equation}\label{Eq:L1L1BoundNoisy}
		m \geq \frac{1}{(1-\tau)^2}\bigg[2\overline{h}\log\Big(\frac{n}{s + \xi/2}\Big) + \frac{7}{5}\Big(s + \frac{\xi}{2}\Big) + \frac{3}{2}\bigg]\,,
  \end{equation}
  where~$0<\tau<1$. Let~$\hat{x}_{\text{noisy}}$ be any solution of
  \begin{equation}\label{Eq:L1L1Noisy}
		\begin{array}[t]{ll}
			\underset{x}{\text{minimize}} & \|x\|_1 + \beta \|x - w\|_1 \\
      \text{subject to} & \|Ax - y\|_2 \leq \sigma\,.
    \end{array}
  \end{equation}
  Then, with overwhelming probability, $\|\hat{x}_{\text{noisy}} - x^\star\|_2 \leq 2\sigma/\tau$, i.e.,  \eqref{Eq:L1L1Noisy} reconstructs~$x^\star$ stably. Our algorithm, described in the next section, adapts easily to the noisy scenario, but we provide reconstruction guarantees only for the noiseless case.
        
  \begin{algorithm}[h!]   
    \caption{\small Adaptive-Rate Sparse Signal Reconstruction}
    \algrenewcommand\algorithmicrequire{\textbf{Input:}}
    \label{Alg:AdaptiveRate}
    \begin{algorithmic}[1]
    \small
    \Require $0\leq\alpha\leq1$, a positive sequence $\{\delta_k\}$, and estimates $\hat{s}_1$ and $\hat{s}_2$ of the sparsity of $x[1]$ and~$x[2]$, respectively.
                
    \Statex
                
    \algrenewcommand\algorithmicrequire{\textbf{Part I: Initialization}}
    \Require
                
    \For{the first two time instants $k=1,2$}                                   \label{SubAlg:BegFirstFor}
                
        \State Set $m_k = 2\hat{s}_k\log(n/\hat{s}_k) + (7/5)\hat{s}_k + 1$       \label{SubAlg:step1}
        \State Generate Gaussian matrix $A_k \in \mathbb{R}^{m_k \times n}$     \label{SubAlg:step2}
        \State Acquire~$m_k$ measurements of $x[k]$: $y[k] = A_k\, x[k]$          \label{SubAlg:step3}  
        \State Find~$\hat{x}[k]$ such that                                        \label{SubAlg:step4}
        $$
					\hat{x}[k]
          \in
          \begin{array}[t]{cl}
						\underset{x}{\arg\min} & \|x\|_1 \\
            \text{s.t.} & A_k\, x = y[k]
          \end{array}
        $$
      \EndFor                                                                     \label{SubAlg:EndFirstFor}
                
      \State Set~$w[2] = f_2(\hat{x}[1])$ and compute                             \label{SubAlg:step5}
      \begin{align*}
					\hat{\xi}_2 
        &:= 
          \big|\{i\,:\, w_i[2] \neq \hat{x}_i[2] = 0\}\big| - \big|\{i\,:\, w_i[2] = \hat{x}_i[2] \neq 0\}\big|
        \\
          \hat{\overline{h}}_2 
        &:= 
          \big|
						\{i\,:\, \hat{x}_i[2] > 0, \,\,\hat{x}_i[2] > w_i[2]\} \cup \{i\,:\, \hat{x}_i[2] < 0, 
            \\& \qquad\qquad\qquad\qquad\qquad\qquad\qquad\qquad\qquad \hat{x}_i[2] < w_i[2]\}
          \big|           
          \,.
      \end{align*}

      \vspace{-0.3cm}
      \State Set 
        $\hat{\overline{m}}_2 = 2\hat{\overline{h}}_2\log\big(n/(\hat{s}_2 + \hat{\xi}_2/2)\big) + (7/5)\big(\hat{s}_2 + \hat{\xi}_2/2\big) + 1$
        \label{SubAlg:step6}
      \State Set $\phi_3 = \hat{\overline{m}}_2$
      \label{SubAlg:step7}
                    
      \Statex                                                         
                
      \algrenewcommand\algorithmicrequire{\textbf{Part II: Online estimation}}
      \Require
                                    
        \For{each time instant $k = 3,4,5,\ldots$}                                     \label{SubAlg:BegSecondFor}
    
            \State Set $m_k = (1 + \delta_k)\phi_k$                                        \label{SubAlg:ChooseMeas}   
                        
            \State Generate Gaussian matrix $A_k \in \mathbb{R}^{m_k \times n}$               \label{SubAlg:step10}                  
                        
            \State Acquire $m_k$ measurements of $x[k]$: $y[k] = A_k\, x[k]$             \label{SubAlg:step11}
                        
            \State Set $w[k] = f_k(\{\hat{x}[i]\}_{i=1}^{k-1})$ and find $\hat{x}[k]$ such that   \label{SubAlg:step12}                         
        $$
            \hat{x}[k]
          \in
          \begin{array}[t]{cl}
						\underset{x}{\arg\min} & \|x\|_1 + \big\|x - w[k]\big\|_1\\
            \text{s.t.} & A_k\, x = y[k]
          \end{array}
        $$
        \State Compute
					\begin{align*}
						\hat{s}_k &= |\{i\,:\, \hat{x}[k] \neq 0\}|
						\\
							\hat{\xi}_k 
            &= 
							\big|\{i\,:\, w_i[k] \neq \hat{x}_i[k] = 0\}\big| - \big|\{i\,:\, w_i[k] = \hat{x}_i[k] \neq 0\}\big|
            \\
							\hat{\overline{h}}_k 
            &= 
              \big|
								\{i\,:\, \hat{x}_i[k] > 0, \,\,\hat{x}_i[k] > w_i[k]\} \cup \{i\,:\, \hat{x}_i[k] < 0, 
                \\& \qquad\qquad\qquad\qquad\qquad\qquad\qquad\qquad\qquad \hat{x}_i[k] < w_i[k]\}
              \big|           
              \,.
          \end{align*}                                    
          \label{SubAlg:h_barxi}
                        
          \vspace{-0.4cm}
          %\State Set $\hat{s}_k = |\{i\,:\, \hat{x}[k] \neq 0\}|$
          %\label{SubAlg:s}
                        
          \State Set $\hat{\overline{m}}_k = 2\hat{\overline{h}}_k\log\big(n/(\hat{s}_k + \hat{\xi}_k/2)\big) + (7/5)\big(\hat{s}_k + \hat{\xi}_k/2\big) + 1$
          \label{SubAlg:mk}
                        
          \vspace{-0.21cm}
          \State Update $\phi_{k+1} = (1-\alpha)\phi_k + \alpha \,\hat{\overline{m}}_k$
          \label{SubAlg:PhiUpdate}
                        
          \vspace{0.1cm}
			\EndFor \label{SubAlg:EndSecondFor}
    \end{algorithmic}
  \end{algorithm}

\section{Online sparse signal estimation}
\label{Sec:DynamicalSignalReconstruction}

	Algorithm~\ref{Alg:AdaptiveRate} describes our online scheme for reconstructing a sparse sequence~$\{x[k]\}$ satisfying~\eqref{Eq:IntroStateSpace}. Although described for a noiseless measurement scenario, the algorithm adapts to the noisy scenario in a straightforward way, as discussed later. Such an adaptation is essential when using it on a real system, e.g., a single-pixel camera~\cite{Duarte08-SinglePixelImagingViaCompressiveSampling}.

\subsection{Algorithm description}          
        
 The algorithm consists of two parts: the initialization, where the first two signals~$x[1]$ and~$x[2]$ are reconstructed using basis pursuit, and the online estimation, where the remaining signals are reconstructed using $\ell_1$-$\ell_1$ minimization. 
%For a simple solver for the $\ell_1$-$\ell_1$ minimization problem, see \aref{Sec:AppSolver}.\footnote{There are many solvers available for basis pursuit. In our experiments, we used \text{spgl1}~\cite{Friedlander08-ProbingParetofrontierBasisPursuit-spgl1,vanDenBerg11-SparseOptimizationWithLeastSquaresConstraints-spgl1}.} For an alternative solver, see~\cite{TranDinh14-APrimalDualAlgorithmicFrameworkForConstrainedConvexMinimization}.
        
\mypar{Part I: Initialization}
In steps~\ref{SubAlg:BegFirstFor}-\ref{SubAlg:EndFirstFor}, we compute the number of measurements~$m_1$ and~$m_2$ according to the bound in~\eqref{Eq:ChandrasekaranBound}, and then reconstruct~$x[1]$ and~$x[2]$ via basis pursuit. The expressions for~$m_1$ and~$m_2$ in step~\ref{SubAlg:step1} require estimates~$\hat{s}_1$ and~$\hat{s}_2$ of the sparsity of~$x[1]$ and~$x[2]$, which are given as input to the algorithm. Henceforth, variables with hats refer to estimates. Steps~\ref{SubAlg:step5}-\ref{SubAlg:step7} initialize the estimator~$\phi_k$: during Part~II of the algorithm, $\phi_k$ should approximate the right-hand side of~\eqref{Eq:L1L1Bound} for~$x[k]$, i.e., with~$s=s_k$, $\overline{h} = \overline{h}_k$, and~$\xi = \xi_k$, where the subscript~$k$ indicates that these are parameters associated with~$x[k]$. 
        
\mypar{Part II: Online estimation}        
The loop in Part~II starts by computing the number of measurements as $m_k = (1+\delta_k)\phi_k$, where~$\delta_k$, an input to the algorithm, is a (positive) safeguard parameter. We take more measurements from~$x[k]$ than the ones prescribed by~$\phi_k$, because $\phi_k$ is only an approximation to the bound in~\eqref{Eq:L1L1Bound}, as explained next. After acquiring measurements from~$x[k]$, we reconstruct it as~$\hat{x}[k]$ via $\ell_1$-$\ell_1$ minimization with $w[k] = f_k(\{\hat{x}[i]\}_{i=1}^{k-1})$ (step~\ref{SubAlg:step12}). Next, in step~\ref{SubAlg:h_barxi}, we compute the sparsity~$\hat{s}_k$ and the quantities in~\eqref{Eq:QualityParameters}, $\hat{\xi}_k$ and~$\hat{\overline{h}}_k$, for~$\hat{x}[k]$. If the reconstruction of~$x[k]$ is perfect, i.e., $\hat{x}[k] = x[k]$, then all these quantities match their true values. In that case, $\hat{\overline{m}}_k$ in step~\ref{SubAlg:mk} will also match the true value of the bound in~\eqref{Eq:L1L1Bound}. Note, however, that the bound for~$x[k]$, $\hat{\overline{m}}_k$, is computed only after~$x[k]$ is reconstructed. Consequently, the number of measurements used in the acquisition of~$x[k]$, $k>2$, is a function of the bound~\eqref{Eq:L1L1Bound} for~$x[k-1]$. Since the bounds for~$x[k]$ and~$x[k-1]$ might differ, we take more measurements than the ones specified by~$\phi_k$ by a factor~$\delta_k$, as in step~\ref{SubAlg:ChooseMeas}. 
%Thus, Algorithm~\ref{Alg:AdaptiveRate} works well when the bound~\eqref{Eq:L1L1Bound} for~$x[k]$ is not much larger than the one for~$x[k-1]$. 
%And that bound depends only on the sparsity~$s_k$ and the quality parameters~$\overline{h}_k$ and~$\xi_k$. 
Also, we mitigate the effect of failed reconstructions by filtering~$\hat{\overline{m}}_k$ with an exponential moving average filter, in step~\ref{SubAlg:PhiUpdate}. Indeed, if reconstruction fails for some~$x[k]$, the resulting~$\hat{\overline{m}}_k$ might differ significantly from the true bound in~\eqref{Eq:L1L1Bound}. The role of the filter is to smooth out such variations. 

\mypar{Extension to the noisy case}
Algorithm~\ref{Alg:AdaptiveRate} can be easily extended to the scenario where the acquisition process is noisy, i.e., $y[k] = A_k x[k] + \eta_k$. Assume that~$\eta_k$ is arbitrary noise, but has bounded magnitude, i.e., we know~$\sigma_k$ such that~$\|\eta_k\|_2 \leq \sigma_k$. In that case, the constraint in the reconstruction problems in steps~\ref{SubAlg:step4} and~\ref{SubAlg:step12} should be replaced by~$\|A_k x - y[k]\|_2 \leq \sigma_k$. The other modification is in steps~\ref{SubAlg:step6} and~\ref{SubAlg:mk}, whose expressions for~$\hat{\overline{m}}_k$ are multiplied by~$1/(1-\tau)^2$ as in~\eqref{Eq:L1L1BoundNoisy}. Our reconstruction guarantees, however, hold only for the noiseless case.

\mypar{Remarks}
We will see in the next section that Algorithm~\ref{Alg:AdaptiveRate} works well when each~$\delta_k$ is chosen according to the prediction quality of~$f_k$: the worse the prediction quality, the larger~$\delta_k$ should be. In practice, it may be more convenient to make~$\delta_k$ constant, as we do in our experiments in \sref{Sec:ExperimentalResults}. Note that the conditions under which our algorithm performs well differ from the majority of prior work. For example, the algorithms in~\cite{Vaswani08-KalmanFilteredCS,Vaswani09-AnalyzingKalmanFilteredCS,Cevher08-CompressiveSensingForBackgroundSubtraction,Ziniel10-TrackingAndSmoothingTimeVaryingSparseSignalsBP,Warnell12-AdaptiveRateCompressiveSensingBackgroundSubtraction,Carmi10-MethodsForSparseSignalRecoveryUsingKalmanFiltering,Kanevsky10-KalmanFilteringForCompressedSensing,Chi13-PETRELS,Balzano10-OnlineIdentificationAndTrackingOfSubspacesFromHighlyIncompleteInformation,Balzano14-LocalConvergenceOfAnAlgorithmForSubspaceIdentificationFromPartialData-GROUSE,Warnell14-AdaptiveRateCompressiveSensingUsingSideInformation} work well when the sparsity pattern of~$x[k]$ varies slowly between consecutive time instants. Our algorithm, in contrast, works well when the quality parameters~$\xi_k$ and~$\overline{h}_k$ of the side information and also the sparsity~$s_k$ of~$x[k]$ vary slowly; in other words, when the quality of the prediction of~$f_k$ varies slowly.

\subsection{Reconstruction guarantees} 
        
The following result bounds the probability with which Algorithm~\ref{Alg:AdaptiveRate} with~$\alpha = 1$ perfectly reconstructs a finite-length sequence~$\{x[i]\}_{i=1}^k$. The idea is to rewrite the condition that~\eqref{Eq:L1L1Bound} applied to~$x[i-1]$ is $(1+\delta_i)$ times larger than~\eqref{Eq:L1L1Bound} applied to~$x[i]$. If that condition holds for the entire sequence then, using \thref{Thm:L1L1} and assuming that the matrices~$A_k$ are drawn independently, we can bound the probability of successful reconstruction. The proof is in \aref{Sec:AppProofLemmaDelta}.
\begin{Lemma}\label{lem:Delta}                	
    Let~$\alpha = 1$, $\underline{m} := \min\big\{m_1, m_2, \min_{i=3,\ldots,k} \hat{\overline{m}}_i\big\}$, and fix~$k > 2$. 
    Let also, for all~$i = 3,\ldots,k$,
    \begin{equation}\label{Eq:LemDelta}
        \delta_i \geq                           
        \frac{2\big[\overline{h}_i \log(\frac{n}{u_i}) - \overline{h}_{i-1} \log(\frac{n}{u_{i-1}})\big] + \frac{7}{5}(u_i - u_{i-1})}{2\overline{h}_{i-1}\log(\frac{n}{u_{i-1}}) + \frac{7}{5}u_{i-1} + 1}\,,
    \end{equation} 
    where~$u_i := s_i + \xi_i/2$. Assume $\hat{s}_q \geq s_q := |\{j: x_j[q] \neq 0\}|$, for~$q=1,2$, i.e., that the initial sparsity estimates~$\hat{s}_1$ and~$\hat{s}_2$ are not smaller than the true sparsity of~$x[1]$ and~$x[2]$. Assume also that the matrices~$\{A_i\}_{i=1}^k$ in Algorithm~\ref{Alg:AdaptiveRate} are drawn independently. Then, the probability (over the sequence of matrices~$\{A_i\}_{i=1}^k$) that Algorithm~\ref{Alg:AdaptiveRate} reconstructs~$x[i]$ perfectly in all time instants $1\leq i \leq k$ is at least
    \begin{equation}\label{Eq:LemProb}                      
        \Big(1 - \exp\Big[-\frac{1}{2}(\underline{m} - \sqrt{\underline{m}})^2\Big]\Big)^k\,.
    \end{equation}          
\end{Lemma}
When the conditions of Lemma~\ref{lem:Delta} hold, the probability of perfect reconstruction decreases with
the length~$k$ of the sequence, albeit at a very slow rate: for example, if~$\underline{m}$ is as small
as~$8$, then~\eqref{Eq:LemProb} equals~$0.9998$ for~$k=10^2$, and~$0.9845$ for~$k=10^4$. If~$\underline{m}$
is larger, these numbers are even closer to~$1$. 

\mypar{Interpretation of~\eqref{Eq:LemDelta}}
As shown in the proof, condition~\eqref{Eq:LemDelta} is equivalent to~$(1 + \delta_i)\overline{m}_{i-1} \geq \overline{m}_i$, where~$\overline{m}_i$ is~\eqref{Eq:L1L1Bound} applied to~$x[i]$. To get more insight about this condition, rewrite it as
\begin{equation}\label{Eq:ConditionOnDeltaRenewed}
		\delta_i 
   \geq                                    
     \frac{\overline{h}_i - \overline{h}_{i-1} + c_1(n)}{\overline{h}_{i-1} + c_2(n)}\,,            
\end{equation}
where
\begin{align*}
		c_1(n) &:= 
			\frac{
				2\overline{h}_{i-1}\log u_{i-1} - 2\overline{h}_i\log u_i + \frac{7}{5}(u_i - u_{i-1})                                    
        }{
        2\log n
        }
     %\label{Eq:Defc1}
    \\
    c_2(n) &:=
			\frac{
				\frac{7}{5}u_{i-1} + 1 - 2\overline{h}_{i-1}\log u_{i-1}                                 
        }{
        2\log n
        }\,.
    %\label{Eq:Defc2}
\end{align*}
Suppose~$\{x[i]\}$ and~$\{\epsilon[i]\}$ are signals for which~$n \gg u_i, \overline{h}_i$. In that case, $c_1(n), c_2(n) \simeq 0$, and condition~\eqref{Eq:ConditionOnDeltaRenewed} tells us that the oversampling factor~$\delta_i$ should be larger than the relative variation of~$\overline{h}_i$ from time~$i-1$ to time~$i$. In general, the magnitude of~$c_1(n)$ and~$c_2(n)$ can be significant, since they approach zero at a relatively slow rate, $o(1/\log n)$. Hence, those terms should not be ignored. 
%Note that~\eqref{Eq:ConditionOnDeltaRenewed} is also a consequence of~$\overline{h}$ being the most important term in~\eqref{Eq:L1L1Bound}.	

\mypar{Remarks on the noisy case}
There is an inherent difficulty in establishing a counterpart of \lref{lem:Delta} for the noisy measurement scenario: namely, the quality parameters~$\xi$ and~$\overline{h}$ in~\eqref{Eq:QualityParameters} are not continuous functions of~$x$. So, no matter how close a reconstructed signal is from the original one, their quality parameters can differ arbitrarily. And, for the noisy measurement case, we can never guarantee that the reconstructed and the original signals are equal; at most, if~\eqref{Eq:L1L1BoundNoisy} holds, they are within a distance~$2\sigma/\tau$, for $0<\tau<1$. 

So far, we have considered~$\{x[k]\}$ and~$\{\epsilon[k]\}$ to be deterministic sequences. In the next section, we will model~$\{\epsilon[k]\}$ (and thus~$\{x[k]\}$) as a Laplacian stochastic process. 
  
\section{Compressive Video Background subtraction}
\label{Sec:ApplicationToBackgroundSubtraction}
	
  We now consider the application of our algorithm to compressive video background subtraction. We start by modeling the problem of compressive background subtraction as the estimation of a sequence of sparse signals satisfying~\eqref{Eq:IntroStateSpace}. Our background subtraction system, based on Algorithm~\ref{Alg:AdaptiveRate}, is then introduced. Finally, we establish reconstruction guarantees for our scheme when~$\epsilon[k]$ in~\eqref{Eq:IntroStateSpaceModel} is Laplacian noise. 
  
  \subsection{Model}
  
  Let~$\{Z[k]\}_{k\geq 1}$ be a sequence of images, each with resolution~$N_1 \times N_2$, and let~$z[k] \in \mathbb{R}^n$ with~$n := N_1 \cdot N_2$ be the (column-major) vectorization of the $k$th image. At time instant~$k$, we collect~$m_k$ linear measurements of~$Z[k]$: $u[k] = A_k z[k]$, where~$A_k \in \mathbb{R}^{m_k \times n}$ is a measurement matrix. We decompose each image~$Z[k]$ as~$Z[k] = X[k] + B$, where~$X[k]$ is the $k$th foreground image, typically sparse, and~$B$ is the background image, assumed known and to be the same in all the images. Let~$x[k]$ and~$b$ be vectorizations of~$X[k]$ and~$B$, respectively.    
  %Given that~$m_k$ is typically smaller than~$n$, it may seem impossible to reconstruct~$Z[k]$ by using only the measurements~$u[k]$. Yet, as noticed in~\cite{Cevher08-CompressiveSensingForBackgroundSubtraction}, in many scenarios, images consist of a background, which typically is known, and a foreground, which typically is unknown. In our setting, we write this as~$Z[k] = X[k] + B$, where~$X[k]$ is the $k$th foreground image and~$B$ is the background image. We assume~$B$ is known. Since the foreground of an image usually occupies a small area, it now seems plausible to be able to reconstruct~$Z[k]$ from~$u[k]$, even when~$m_k \ll n$. To see how, let~$x[k]$ and~$b$ be vectorizations of~$X[k]$ and~$B$, respectively. 
  Because the background image is known, we take measurements from it using~$A_k$: $u^b[k] = A_k b$. Then, as suggested in~\cite{Cevher08-CompressiveSensingForBackgroundSubtraction}, we subtract~$u^b[k]$ to~$u[k]$:
  \begin{equation}\label{Eq:BackgroundSubtrMeasurements}
		y[k]:= u[k] - u^b[k] = A_k (z[k] - b) = A_k x[k]\,.
  \end{equation}
  This equation tells us that, although we cannot measure the foreground image~$x[k]$ directly, we can still construct a
  vector measurements, $y[k]$, as if we would. Given that~$x[k]$ is usually sparse, the theory of compressed sensing tells us that it can be reconstructed by solving, for example, basis pursuit~\cite{Candes06-RobustUncertaintyPrinciplesExactSignalReconstructionHighlyIncomplete,Donoho06-CompressedSensing}. Specifically, if~$x[k]$ has sparsity~$s_k$ and the entries of~$A_k$ are realizations of i.i.d.\ zero-mean Gaussian random variables with variance~$1/m_k$, then~$2s_k\log(n/s_k) + (7/5)s_k + 1$ measurements suffice to reconstruct~$x[k]$ perfectly~\cite{Chandrasekaran12-ConvexGeometryLinearInverseProblems} [cf.\ \eqref{Eq:ChandrasekaranBound}].
	%After reconstructing~$x[k]$, we can obtain~$z[k]$ by adding it to~$b$: $z[k] = x[k] + b$. We thus see that this approach not only allows reconstructing each~$Z[k]$ from a limited number of measurements, but also performs background subtraction as a by-product.
	
	Notice that~\eqref{Eq:BackgroundSubtrMeasurements} is exactly the equation of measurements in~\eqref{Eq:IntroStateSpaceObservations}. Regarding equation~\eqref{Eq:IntroStateSpaceModel}, we will use it to model the estimation of the foreground of each frame, $x[k]$, from previous foregrounds, $\{x[i]\}_{i=1}^{k-1}$. We use a motion-compensated extrapolation technique, as explained in \ssref{SubSec:Extrapolation}. This technique is known to produce image estimates with an error well modeled as Laplacian and, thus, each~$\|\epsilon[k]\|_1$ is expected to be small. This perfectly aligns with the way we integrate side information in our reconstruction scheme: namely, the second term in the objective of the optimization problem in step~\ref{SubAlg:step12} of Algorithm~\ref{Alg:AdaptiveRate} is nothing but~$\|\epsilon[k]\|_1$.

  %where both~$\{x[k]\}$ and~$\{\epsilon[k]\}$ are unknown sequences of sparse vectors in~$\mathbb{R}^n$. Note that while the sparsity of the vectors~$x[k]$ is justified by the small area occupied by the foreground~\cite{Cevher08-CompressiveSensingForBackgroundSubtraction}, the sparsity of the vectors~$\epsilon[k]$ is justified by motion estimation models~\cite{Deligiannis12-SideInformationDependentCorrelationChannelEstimationHashBasedDistributedVideoCoding,Deligiannis14-MaximumLikelihoodLaplacianCorrelationChannelEstimationLayeredWynerZivCoding}. We also assume that the entries of each matrix~$A_k \in \mathbb{R}^{m_k \times n}$ are i.i.d.\ zero-mean Gaussian random variables with variance~$1/m_k$. 
  %Although Gaussian matrices are hard to implement in practical systems, they have optimal performance. There are, however, other more practical matrices with a similar performance~\cite{Berinde08-SparseRecoveryUsingSparseMatrices,Liutkus14-ImagingWithNature-CompressiveImagingUsingMultiplyScatteringMedium}.                 

  \subsection{Our background subtraction scheme}
		
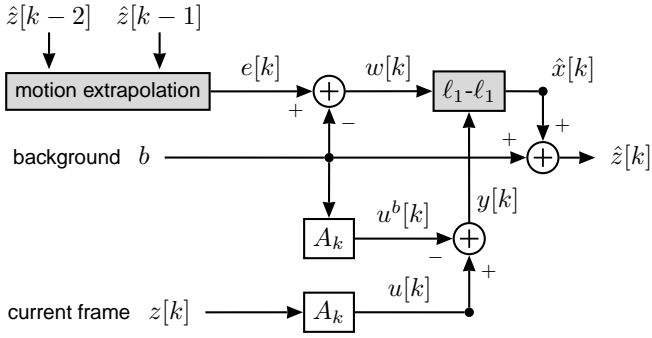
\begin{figure}[t]
		\centering
		\psscalebox{0.98}{
		\begin{pspicture}(8.2,4.7)
			\pssignal(0.3,4.5){zpp}{$\hat{z}[k-2]$}
			\pssignal(1.8,4.5){zp}{$\hat{z}[k-1]$}
			
			\pssignal(0.74,2.6){b}{{\footnotesize \textsf{background}}\,\, $b$}			
			\pssignal(0.98,0.5){z}{{\footnotesize \textsf{current frame}}\,\, $z[k]$}
			
			\pssignal(8.2,2.6){zh}{$\hat{z}[k]$}
			
			\psblock[FillColor=black!15!white](1.1,3.5){ME}{\footnotesize \textsf{motion extrapolation}}
			\psblock(4.1,1.5){Ak1}{$A_k$}
			\psblock(4.1,0.5){Ak2}{$A_k$}
			\psblock[FillColor=black!15!white](6.0,3.5){L1L1}{$\ell_1$-$\ell_1$}
			
			\dotnode(4.1,2.6){dot1}
			\dotnode(7.0,3.5){dot2}
			\nput[labelsep=3pt]{45}{dot2}{$\hat{x}[k]$}
			\dotnode(6.0,0.5){dot3}
			
			\pscircleop(4.1,3.5){op1}
			\nput{210}{op1}{$\scriptstyle +$}
			\nput{300}{op1}{$\scriptstyle -$}
			
			\pscircleop(6.0,1.5){op2}
			\nput{210}{op2}{$\scriptstyle -$}
			\nput{300}{op2}{$\scriptstyle +$}
			\nput{58}{op2}{$y[k]$}
			
			\pscircleop(7.0,2.6){op3}
			\nput{60}{op3}{$\scriptstyle +$}
			\nput{150}{op3}{$\scriptstyle +$}
			
			%-------------
			\psset{style=Arrow,arrowinset=0.05,arrowsize=5pt,labelsep=3pt}
			
			\psline(0.3,4.3)(0.3,3.8)
			\psline(1.8,4.3)(1.8,3.8)
			
			\ncline{ME}{op1}
			\naput{$e[k]$}
			\ncline{-}{b}{dot1}
			\ncline{dot1}{Ak1}
			\ncline{z}{Ak2}
			\ncline{dot1}{op1}
			\ncline{op1}{L1L1}
			\naput{$w[k]$}
			\ncline{-}{L1L1}{dot2}
			\ncline{dot2}{op3}
			\ncline{op2}{L1L1}			
			\ncline{-}{Ak2}{dot3}
			\naput{$u[k]$}
			\ncline{dot3}{op2}
			\ncline{Ak1}{op2}
			\naput{$u^b[k]$}
			\ncline{op3}{zh}
			
			\ncline{dot1}{op3}
			
			%\psgrid
		\end{pspicture}
		}
	  \caption{Block diagram of Algorithm~\ref{Alg:AdaptiveRate} when applied to background subtraction. The main blocks are highlighted.}
		\label{Fig:BlockDiagram}
	\end{figure}

	Fig.~\ref{Fig:BlockDiagram} shows the block diagram of our compressive background subtraction scheme and, essentially, translates Algorithm~\ref{Alg:AdaptiveRate} into a diagram. The scheme does not apply to the reconstruction of the first two frames, which are reconstructed as in~\cite{Cevher08-CompressiveSensingForBackgroundSubtraction}, i.e., by solving basis pursuit. This corresponds to Part I of Algorithm~\ref{Alg:AdaptiveRate}. The scheme in Fig.~\ref{Fig:BlockDiagram} depicts Part II of Algorithm~\ref{Alg:AdaptiveRate}. The motion extrapolation module constructs a motion-compensated prediction~$e[k]$ of the current frame, $z[k]$, by using the two past (reconstructed) frames, $\hat{z}[k-2]$ and~$\hat{z}[k-1]$. Motion estimation is performed in the image domain ($z[k]$) rather than in the foreground domain ($x[k]$), as the former contains more texture, thereby yielding a more accurate motion field.	Next, the background frame~$b$ is subtracted from~$e[k]$ to obtain a prediction of the foreground~$x[k]$, i.e., the side information~$w[k]$. These two operations are modeled in Algorithm~\ref{Alg:AdaptiveRate} with the function~$f_k$, which takes a set of past reconstructed signals (in our case, $\hat{x}[k-2]$ and~$\hat{x}[k-1]$, to which we add~$b$, obtaining~$\hat{z}[k-2]$ and~$\hat{z}[k-1]$, respectively), and outputs the side information~$w[k]$. This is one of the inputs of the $\ell_{1}$-$\ell_{1}$ block, which solves the optimization problem~\eqref{Eq:L1L1Simple}. To obtain the other input, i.e., the set of foreground measurements~$y[k]$, we proceed as specified in equation~\eqref{Eq:BackgroundSubtrMeasurements}: we take measurements~$u[k] = A_k z[k]$ of the current frame and, using the same matrix, we take measurements of the background~$u[k] = A_k b$. Subtracting them we obtain~$y[k] = u[k] - u^b[k]$. The output of the $\ell_1$-$\ell_1$ module is the estimated foreground~$\hat{x}[k]$, from which we obtain the estimate of the current frame as $\hat{z}[k] = \hat{x}[k] + b$.
	
%By subtracting the vectorized background frame~$b$ from the vectorized motion-compensated prediction~$e[k]$ yields the foreground prediction frame~$w[k]$, which acts as the side information signal in the $\ell_{1}$-$\ell_{1}$ minimization problem defined in~\eqref{Eq:L1L1}. As stated in Section~\ref{Sec:Model}, the compressive measurements, defining the constrains in~\eqref{Eq:L1L1}, are obtained as~$y[k]=u^{b}[k]-u[k]$, where~$u^{b}[k]$ and~$u[k]$ denote the measurements taken from the background and the current frame, respectively. The required number of measurements (and thus, the size of the sensing matrix~$A_k$) is estimated online  as described in Algorithm~\ref{Alg:AdaptiveRate}. Finally, by solving  the $\ell_{1}$-$\ell_{1}$ minimization problem~\eqref{Eq:L1L1}, we obtain the vectorized foreground image~$\hat{x}[k]$ and, by adding the background signal~$b$, the reconstructed current frame~$\hat{z}[k]$.          

	\subsection{Motion-compensated extrapolation}
	\label{SubSec:Extrapolation}
		
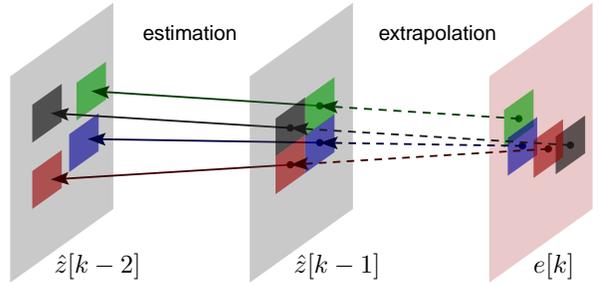
\begin{figure}[t]
		\centering
		\psscalebox{0.98}{
		\begin{pspicture}(8.2,3.8)
					
			\psset{blendmode=2}
		
			\def\frameskew#1{
				%\pspolygon*[linecolor=black!15!white](0,0)(1.2,1.2)(1.2,4.5)(0,3.3)
				\pspolygon*[linecolor=#1,opacity=0.3](0,0)(1.4,1.0)(1.4,3.8)(0,2.8)
			}
			
			\def\littlebk[#1]#2{
				\pspolygon*[linestyle=#1,linecolor=#2,opacity=0.6](0,0)(0.4,0.2857)(0.4,0.8)(0,0.5143)
			}			
			
			\def\ball#1{				
					\psellipse*[linecolor=#1](0,0)(0.05,0.05)
			}
			
			\psset{arrowsize=5.0pt,arrowinset=0.18,linewidth=0.7pt}
			
			% Frame 1
			\rput(0,0){\frameskew{black!70!white}}
			\rput(0.3,1.0){\littlebk[solid]{red}}
			\rput(0.9,2.2){\littlebk[solid]{green}}
			\rput(0.8,1.5){\littlebk[solid]{blue}}
			\rput(0.3,1.9){\littlebk[solid]{black}}
			
			\psline[linecolor=green!40!black]{->}(4.2,2.4)(1.1,2.6)
			\psline[linecolor=red!40!black]{->}(3.8,1.6)(0.55,1.4)
			\psline[linecolor=blue!40!black]{->}(4.2,1.9)(1.05,1.95)
			\psline[linecolor=black!90!white]{->}(3.8,2.1)(0.55,2.3)
			
			% Frame 2
			\rput(3.25,0){\frameskew{black!70!white}}
			\rput(4.0,1.4857){\littlebk[solid]{blue}}
			\rput(3.6,1.2){\littlebk[solid]{red}}
			\rput(4.0,1.99){\littlebk[solid]{green}}
			\rput(3.6,1.7143){\littlebk[solid]{black}}
			
			\rput(4.2,2.4){\ball{darkgreen}}
			\rput(3.8,1.6){\ball{darkred}}
			\rput(3.8,2.1){\ball{black!90!white}}
			\rput(4.2,1.9){\ball{darkblue}}

			% Frame 3
			\rput(6.5,0){\frameskew{red!70!white}}
						
			\rput(6.70,1.83){\littlebk[solid]{green}}
			\rput(7.40,1.47){\littlebk[solid]{black}}
			\rput(6.75,1.45){\littlebk[solid]{blue}}
			\rput(7.1,1.42){\littlebk[solid]{red}}
			
			\rput(6.9,2.2258){\ball{darkgreen}}
			\rput(7.3,1.8154){\ball{darkred}}
			\rput(6.95,1.8563){\ball{darkblue}}
			\rput(7.6,1.8662){\ball{black!90!white}}
			
			\psset{linestyle=dashed,dash=3pt 2.5pt}
			\psline[linecolor=green!40!black]{->}(6.9,2.2258)(4.2,2.4)
			\psline[linecolor=red!40!black]{->}(7.3,1.8154)(3.8,1.6)
			\psline[linecolor=blue!40!black]{->}(6.95,1.8563)(4.2,1.9)
			\psline[linecolor=black!90!white]{->}(7.6,1.8662)(3.8,2.1)
			
			% Text and math
			\rput[lt](0.60,0.4){$\hat{z}[k-2]$}
			\rput[lt](3.85,0.4){$\hat{z}[k-1]$}
			\rput[lt](7.10,0.4){$e[k]$}
			
			\rput[lt](1.8,3.5){\footnotesize \textsf{estimation}}
			\rput[lt](5.0,3.5){\footnotesize \textsf{extrapolation}}
								
			%\psgrid
		\end{pspicture}
		}
		\vspace{0.1cm}
	  \caption{Scheme of motion-compensated extrapolation. We use the motion between matching blocks in $\hat{z}[k-2]$ and $\hat{z}[k-1]$ to create an estimate~$e[k]$ of frame~$z[k]$.
		}
			%Motion is estimated from reconstructed frames $\hat{z}[k-2]$ and~$\hat{z}[k-1]$ to create an estimated frame~$e[k]$.}
		\label{Fig:MCEschema}
	\end{figure}

	To obtain an accurate predition~$e[k]$, we use a motion-compensated extrapolation technique similar to what is used in distributed video coding for generating decoder-based motion-compensated predictions~\cite{Girod05-DistributedVideoCoding,Deligiannis12-SideInformationDependentCorrelationChannelEstimationHashBasedDistributedVideoCoding,natario2006extrapolating}.
	Our technique is illustrated in \fref{Fig:MCEschema}. In the first stage, we perform forward block-based motion estimation between the reconstructed frames~$\hat{z}[k-2]$ and~$\hat{z}[k-1]$. The block matching algorithm is performed with half-pel accuracy and considers a block size of~$\gamma\times{\gamma}$ pixels and a search range of~$\rho$ pixels. The required interpolation for half-pel motion estimation is performed using the 6-tap filter of H.264/AVC~\cite{wiegand2003overview}. In addition, we use the $\ell_{1}$-norm (or sum of absolute differences: SAD) as error metric. The resulting motion vectors are then spatially smoothed by applying a weighted vector-median filter~\cite{alparone1996adaptively}. The filtering improves the spatial coherence of the resulting motion field by removing outliers (i.e., motion vectors that are far from the true motion field). Assuming linear motion between~$\hat{z}[k-2]$ and~$\hat{z}[k-1]$, and~$\hat{z}[k-1]$ and~$\hat{z}[k]$, we linearly project the motion vectors between~$\hat{z}[k-2]$ and~$\hat{z}[k-1]$ to obtain~$e[k]$, our estimate of~$z[k]$; see \fref{Fig:MCEschema}. During motion compensation, pixels in~$e[k]$ that belong to overlapping prediction blocks are estimated as the average of their corresponding motion-compensated pixel predictors in~$\hat{z}[k-1]$. Pixels in uncovered areas (i.e., no motion-compensated predictor is available) are estimated by taking averaging the three neighbor pixel values in~$e[k]$ (up, left and up-left pixel positions, following a raster scan of the frame) and the corresponding pixel in $\hat{z}[k-1]$. 
	      	      
	%This operation makes the function modeling the extrapolation, i.e., $f_k$, nonlinear. 
	
  \subsection{Reconstruction guarantees for Laplacian modeling noise}
        
  It is well known that the noise produced by a motion-compensated prediction module, as the one just described, is Laplacian \cite{Girod05-DistributedVideoCoding,Deligiannis14-MaximumLikelihoodLaplacianCorrelationChannelEstimationLayeredWynerZivCoding}. In our model, that corresponds to each~$\epsilon[k]$ in~\eqref{Eq:IntroStateSpaceModel} being Laplacian. We assume each~$\epsilon[k]$ is independent from the matrix of measurements~$A_k$. 
  
  \mypar{Model for~\boldmath{$\epsilon[k]$}}
  As in \cite{Girod05-DistributedVideoCoding,Deligiannis14-MaximumLikelihoodLaplacianCorrelationChannelEstimationLayeredWynerZivCoding,fan2010transform} (and references therein), we assume that~$\epsilon[k]$ is independent from~$\epsilon[l]$, for $k\neq l$, and that the entries of each~$\epsilon[k]$ are independent and have zero-mean. The probability distribution of~$\epsilon[k]$ is then
  \begin{align}
		\mathbb{P}(\epsilon[k] \leq  u) 
    &=
    \mathbb{P}(\epsilon_1[k] \leq  u_1, \,\,\epsilon_2[k] \leq u_2,\ldots, \,\,\epsilon_n[k] \leq  u_n) 
    \notag
    \\
    &=
    \prod_{j=1}^n \mathbb{P}(\epsilon_j[k] \leq  u_j)
    \notag
    \\
    &=                      
		\prod_{j=1}^n \int_{-\infty}^{u_j} \frac{\lambda_j}{2} \exp\big[-\lambda_j |\epsilon_j| \      \big]\,d\epsilon_j\,,
    \label{Eq:ProbDistributionEpsilon}
  \end{align}
  where~$u \in \mathbb{R}^n$, and~$\lambda_j\geq 0$ is the parameter of the distribution of~$\epsilon_j[k]$. The entries of~$\epsilon[k]$, although independent, are not identically distributed,
  since they have possibly different parameters~$\lambda_j$. The variance~$\sigma_j^2$ of each
  component~$\epsilon_j[k]$ is given by~$\sigma_j^2 = 2/\lambda_j^2$. 
  
  \mypar{Resulting model for~\boldmath{$x[k]$}}
  The sequence~$\{\epsilon[k]\}$ being stochastic implies that~$\{x[k]\}$ is also stochastic. Indeed, if each~$f_k$ in~\eqref{Eq:IntroStateSpace} is measurable, then~$\{x[k]\}_{k\geq 2}$ is a
  sequence of random variables. Given the independence across time and across components
  of the sequence~$\{\epsilon[k]\}$, the distribution of~$x[k]$ given~$\{x[i]\}_{i=1}^{k-1}$ is also
  Laplacian, yet not necessarily with zero-mean. That is, for~$u \in \mathbb{R}^n$ and~$k\geq 2$,
  \begin{align}
			&\mathbb{P}\Big(x[k] \leq u \,\,\big|\,\, \{x[i]\}_{i=1}^{k-1}\Big)
		\notag
    \\
    &=
			\mathbb{P}\Big(f_k(\{x[i]\}_{i=1}^{k-1}) + \epsilon[k] \leq u \,\,\big|\,\, \{x[i]\}_{i=1}^{k-1}\Big)
    \notag
    \\
    &=
			\mathbb{P}\Big(\epsilon[k] \leq u - f_k(\{x[i]\}_{i=1}^{k-1})\,\,\big|\,\, \{x[i]\}_{i=1}^{k-1}\Big)
    \notag
    \\
    &=
			\prod_{j=1}^n \int_{-\infty}^{u_j - [f_k( \{x[i]\}_{i=1}^{k-1})]_j} \frac{\lambda_j}{2} \exp\big[-\lambda_j |\epsilon_j| \big]\,d\epsilon_j
    \notag
    \\
    &=
			\prod_{j=1}^n \int_{-\infty}^{u_j} \frac{\lambda_j}{2} \exp\Big[\!-\lambda_j \big|z_j - [ f_k(\{x[i]\}_{i=1}^{k-1})]_j\big| \Big]\,dz_j
		\label{Eq:ProbDistributionX}
  \end{align}
  where~$[f_k(\{x[i]\}_{i=1}^{k-1})]_j$ is the $j$th component of~$f_k(\{x[i]\}_{i=1}^{k-1})$. 
  %In the last step, we changed the integration variable from~$\epsilon_j$ to $z_j = \epsilon_j + [f_k(\{x[i]\}_{i=1}^{k-1})]_j$. 
  In words, the distribution of each component of~$x[k]$ conditioned on all past realizations~$x[i]$, $1\leq i < k$, is Laplacian with mean~$[f_k(\{x[i]\}_{i=1}^{k-1})]_j$ and parameter~$\lambda_j$. Furthermore, it is independent from the other components. 
  
  \mypar{Reconstruction guarantees}
  Note that~$\{x[k]\}$ and $\{\epsilon[k]\}$ being stochastic processes implies that the quantities in~\eqref{Eq:QualityParameters}, which we will denote with~$\xi_k$ and~$\overline{h}_k$ for signal~$x[k]$, are random variables. Hence, at each time~$k$, the conditions of Theorem~\ref{Thm:L1L1}, namely that $\overline{h}_k > 0$ and that there is at least one index~$i$ such that~$x_i[k] = w_i[k] = 0$, become events, and may or may not hold. We now impose conditions on the variances~$\sigma_j^2 = 2/\lambda_j$ that guarantee the conditions of Theorem~\ref{Thm:L1L1} are satisfied and, thus, that $\ell_1$-$\ell_1$ minimization reconstructs~$x[k]$ perfectly, with high probability. Given a set~$S \in \{1,\ldots,n\}$, we use~$S^c$ to denote its complement in~$\{1,\ldots,n\}$. 
  \begin{Theorem}
		\label{Thm:LaplacianModelingNoise}
		Let~$w \in \mathbb{R}^n$ be given. Let~$\epsilon$ have distribution~\eqref{Eq:ProbDistributionEpsilon}, where the variance of component~$\epsilon_j$ is $\sigma_j^2 = 2/\lambda_j^2$. Define~$x^\star := w + \epsilon$, and the sets $\Sigma := \{j\,:\, \sigma_j^2 \neq 0\}$ and $\mathcal{W} := \{j\,:\, w_j \neq 0\}$.
% 		\begin{align*}		
% 			\Sigma &:= \{j\,:\, \sigma_j^2 \neq 0\}
% 			\\
% 			\mathcal{W} &:= \{j\,:\, w_j \neq 0\}\,.
% 		\end{align*}
		Assume~$\Sigma^c \cap \mathcal{W}^c \neq \emptyset$, that is, there exists~$j$ such that~$\sigma_j^2 = 0$ and~$w_j = 0$. Assume~$A \in \mathbb{R}^{m \times n}$ is generated as in \thref{Thm:L1L1} with a number of measurements 
		\begin{multline}\label{Eq:lemConditionsLaplacianModelNumMeas}
			m \geq 2(\mu + t)\log\bigg(\frac{n}{\big|\Sigma\big| + \frac{1}{2}\big|\Sigma^c \cap \mathcal{W}\big|}\bigg) 
			\\
			+ \frac{7}{5}\Big(\big|\Sigma\big| + \frac{1}{2}\big|\Sigma^c \cap \mathcal{W}\big|\Big) + 1\,,
		\end{multline}
		for some~$t>1$, where
		$
			\mu := \frac{1}{2}\sum_{j \in \Sigma} \big[1 + \exp\big(-\sqrt{2}|w_j|/\sigma_j\big)\big]
		$.
		Let~$\hat{x}$ denote the solution of $\ell_1$-$\ell_1$ minimization~\eqref{Eq:L1L1Simple}. Then,		
		\begin{multline}\label{Eq:lemConditionsLaplacianProb}
			\mathbb{P}\big(\hat{x} = x^\star\big)
			\geq 
			\bigg[1 - \exp\Big(-\frac{(m - \sqrt{m})^2}{2}\Big)\bigg]
			\times\\\times
			\bigg[1 - \exp\Big(-\frac{2\mu^2}{|\Sigma|}\Big) - \exp\Big(-\frac{2(t-1)^2}{|\Sigma|}\Big)\bigg]\,.
		\end{multline}		
	\end{Theorem}
	The proof is in \aref{Sec:AppProofLemmaConditions}.	By assuming each component~$\epsilon_j$ is Laplacian with parameter~$\lambda_j = \sqrt{2}/\sigma_j$ (independent from the other components), \thref{Thm:LaplacianModelingNoise} establishes a lower bound on the number of measurements that guarantee perfect reconstruction of~$x^\star$ with probability as in~\eqref{Eq:lemConditionsLaplacianProb}. Note that all the quantities in~\eqref{Eq:lemConditionsLaplacianModelNumMeas} are deterministic. This contrasts with the direct application of \thref{Thm:L1L1} to the problem, since the right-hand side of~\eqref{Eq:L1L1Bound} is a function of the random variables~$s$, $\overline{h}$, and~$\xi$. The assumption $\Sigma^c \cap \mathcal{W}^c \neq \emptyset$ implies~$\Sigma^c \neq \emptyset$, which means that some components of~$\epsilon$ have zero variance and, hence, are equal to zero with probability~$1$. 
	%Because~$\Sigma^c \cap \mathcal{W} \neq \emptyset$, at least one of those components coincides with a zero component of~$w$, i.e., $\epsilon_j = w_j = 0$ for some~$j$. 
	Note that, provided the variances~$\sigma_j^2$ are known, all the quantities in~\eqref{Eq:lemConditionsLaplacianModelNumMeas}, and consequently in~\eqref{Eq:lemConditionsLaplacianProb}, are known. 
			
	The proof of \thref{Thm:LaplacianModelingNoise} uses the fact that the sparsity of~$x^\star$ is $s = |\Sigma| + |\Sigma^c \cap \mathcal{W}|/2$ with probability~$1$. This implies that the bound in~\eqref{Eq:lemConditionsLaplacianModelNumMeas} is always smaller than the one for basis pursuit in~\eqref{Eq:ChandrasekaranBound} whenever~$\mu + t < s = |\Sigma| + |\Sigma^c \cap \mathcal{W}|$. Since~$\mu \leq |\Sigma|$, this holds if~$t < |\Sigma^c \cap \mathcal{W}|/2$.
	%That is partly because knowing the cardinality of the sets~$\Sigma$ and~$\Sigma^c \cap \mathcal{W}$ implies knowing the sparsity~$s$ of~$x^\star$, a fact used in the proof. Arguably, $\mu$ and~$t$ are the parameters that most influence the bound in~\eqref{Eq:lemConditionsLaplacianModelNumMeas}. While~$t > 1$ is free, $\mu$ depends on the cardinality of~$\Sigma$ and, for~$j \in \Sigma$, depends on the magnitude of~$w_j$ and standard deviation~$\sigma_j$ of~$\epsilon_j$. Specifically, $\mu$ (and thus~$m$) increases when~$|\Sigma|$ increases, $|w_j|$ decreases, and~$\sigma_j$ increases. The sizes of~$\Sigma$ and~$\Sigma \cap \mathcal{W}$ have a smaller influence but, when~$n$ is large enough so that the first term in~\eqref{Eq:lemConditionsLaplacianModelNumMeas} becomes dominant, the larger these sets, the smaller the bound in~\eqref{Eq:lemConditionsLaplacianModelNumMeas}. In turn, the probability of successful reconstruction in~\eqref{Eq:lemConditionsLaplacianProb} increases with~$m$, $\mu$, and~$t$, and decreases with~$|\Sigma|$.
	
	We state without proof a consequence of~\thref{Thm:LaplacianModelingNoise} that is obtained by reasoning as in \lref{lem:Delta}:
	\begin{Corollary}
		\label{Cor:LaplacianModellingNoiseAlg}
		Let~$\{\epsilon[k]\}$ be a stochastic process where~$\epsilon[k]$ has distribution~\eqref{Eq:ProbDistributionEpsilon} and each~$\epsilon[k]$ is independent from~$\epsilon[l]$, $k\neq l$. Assume that~$\{x[k]\}$ is generated as in~\eqref{Eq:IntroStateSpaceModel} and consider Algorithm~\ref{Alg:AdaptiveRate} with~$\alpha = 1$ at iteration~$k>2$. Assume that~$\epsilon[k]$ and~$A_k$ are independent. Assume also, for~$i=3,\ldots,k$, that
		\begin{multline}\label{Eq:CorDelta}
			\delta_i \geq 
				\bigg\{ 
					2\Big[ (\mu_i + t_i)\log\Big(\frac{n}{u_i}\Big) - (\mu_{i-1} + t_{i-1})\log\Big(\frac{n}{u_{i-1}}\Big)\Big] \\ + \frac{7}{5}(u_i - u_{i-1})
				\bigg\}
				\Big/
				\bigg\{	
					2(\mu_{i-1} + t_{i-1})\log\Big(\frac{n}{u_{i-1}}\Big) + \frac{7}{5}u_{i-1} + 1
				\bigg\}
					\,,
		\end{multline}
		where~$u_i := |\Sigma_i| + |\Sigma_i^c \cap \mathcal{W}_i|/2$, and the quantities~$\mu_i$, $t_i$, $\Sigma_i$, and~$\mathcal{W}_i$ are defined as in \thref{Thm:LaplacianModelingNoise} for signal~$x[i]$. Assume the initial sparsity estimates satisfy $\hat{s}_1 \geq s_1$ and $\hat{s}_2 \geq s_2$ with probability~$1$, where~$s_1$ and~$s_2$ are the sparsity of~$x[1]$ and~$x[2]$. Then, the probability over the sequences $\{A_i\}_{i=1}^k$ and $\{\epsilon[i]\}_{i=1}^k$ that Algorithm~\ref{Alg:AdaptiveRate} reconstructs $x[i]$ perfectly in all time instants $1\leq i \leq k$ is at least
		\begin{multline*}%\label{Eq:CorProb}
			\prod_{i=1}^k 
			\bigg[1 - \exp\Big(-\frac{(m_i - \sqrt{m_i})^2}{2}\Big)\bigg]		
			\bigg[1 - \exp\Big(-\frac{2\mu_i^2}{|\Sigma_i|}\Big) \\- \exp\Big(-\frac{2(t_i-1)^2}{|\Sigma_i|}\Big)\bigg]\,.
		\end{multline*} 
	\end{Corollary}
	\coref{Cor:LaplacianModellingNoiseAlg} establishes reconstruction guarantees of Algorithm~\ref{Alg:AdaptiveRate} when the modeling noise~$\epsilon[k]$ in~\eqref{Eq:IntroStateSpaceModel} is Laplacian. In contrast with \lref{lem:Delta}, the bound in~\eqref{Eq:CorDelta} is a function of known parameters, but it requires the variances~$\sigma_j^2[i]$ of each~$\epsilon_j[i]$, which can be estimated from the past frame in a block-based way~\cite{Deligiannis12-SideInformationDependentCorrelationChannelEstimationHashBasedDistributedVideoCoding,Deligiannis14-MaximumLikelihoodLaplacianCorrelationChannelEstimationLayeredWynerZivCoding}. For some insight on~\eqref{Eq:CorDelta}, assume~$u_i \simeq u_{i-1}$, $t_i = t_{i-1}$, and that~$n$ is large enough so that terms not depending on it are negligible. Then, \eqref{Eq:CorDelta} becomes $\delta_i \gtrsim (\mu_i - \mu_{i-1})/(\mu_{i-1} + t_{i-1})$, and we can select
	\begin{align}
		\delta_i &= 2\kappa - 1 \simeq \frac{\kappa - \frac{1}{2}}{\frac{1}{|\Sigma_{i-1}|} + \frac{1}{2}} = \frac{|\Sigma_i| - \frac{1}{2}|\Sigma_{i-1}|}{1 + \frac{1}{2}|\Sigma_{i-1}|} 				
		\geq \frac{\mu_i - \mu_{i-1}}{\mu_{i-1} + t_{i-1}}\,,
		\label{Eq:DeltaLaplacianApprox}
	\end{align}
	where~$\kappa := |\Sigma_i|/|\Sigma_{i-1}|$, and the inequality is due to~$|\Sigma_i|/2\leq \mu_i \leq |\Sigma_i|$. The approximation in~\eqref{Eq:DeltaLaplacianApprox} holds if~$|\Sigma_{i-1}| \gg 2$, which is often the case in practice. The expression in~\eqref{Eq:DeltaLaplacianApprox} tells us that, for large~$n$, $\delta_i$ is mostly determined by the ratio~$\kappa$: if~$\kappa > 1$ (resp. $<1$), then we should select~$\delta_i > 1$ (resp. $<1$). We observe that, in practice, \eqref{Eq:CorDelta} and \eqref{Eq:DeltaLaplacianApprox} give conservative estimates for~$\delta_i$. We will see in the next section that selecting a small, constant~$\delta_i$ (namely~$0.1$) leads to excellent results without compromising perfect reconstruction.

	\begin{figure*}
		\raggedleft
		
		\begin{pspicture}(18,0.2)
			%\rput[l](0.0,2.6){\small \sf Frame \#}
			\rput[l](0.0,-1.4){\sf Original}
			
			\rput[b](3.380,0.1){\sf background}
			\rput[b](6.580,0.16){\sf frame $\mathsf{4}$}
			\rput[b](9.750,0.16){\sf frame $\mathsf{100}$}
			\rput[b](12.92,0.16){\sf frame $\mathsf{150}$}
			\rput[b](16.05,0.16){\sf frame $\mathsf{250}$}

			\rput[lb](0.0,-4.8){\sf Estimated}
			
			\rput[lb](0.0,-7.9){\sf Reconstructed} 
			
			\rput[lb](0.0,-11.0){\sf Reconstructed foreground} 
			\rput[lb](0.0,-11.5){\sf (binarized)} 
			
			%\psgrid
		\end{pspicture}
		
% 		\def\widthHall{3.035cm}
% 		\includegraphics[width=\widthHall]{figures/Hall_original_f1.eps}
% 		\includegraphics[width=\widthHall]{figures/Hall_original_f4.eps}
% 		\includegraphics[width=\widthHall]{figures/Hall_original_f100.eps}
% 		\includegraphics[width=\widthHall]{figures/Hall_original_f150.eps}
% 		\includegraphics[width=\widthHall]{figures/Hall_original_f250.eps}
% 		\hspace{0.35cm}
% 		
% 		\vspace{0.1cm}
% 		\includegraphics[width=\widthHall]{figures/Hall_estimated_image_f4.eps}
% 		\includegraphics[width=\widthHall]{figures/Hall_estimated_image_f100.eps}
% 		\includegraphics[width=\widthHall]{figures/Hall_estimated_image_f150.eps}
% 		\includegraphics[width=\widthHall]{figures/Hall_estimated_image_f250.eps}
% 		\hspace{0.35cm}
% 		
% 		\vspace{0.1cm}
% 		\includegraphics[width=\widthHall]{figures/Hall_reconstructed_image_f4.eps}
% 		\includegraphics[width=\widthHall]{figures/Hall_reconstructed_image_f100.eps}
% 		\includegraphics[width=\widthHall]{figures/Hall_reconstructed_image_f150.eps}
% 		\includegraphics[width=\widthHall]{figures/Hall_reconstructed_image_f250.eps}
% 		\hspace{0.35cm}
% 		
% 		\vspace{0.1cm}
% 		\includegraphics[width=\widthHall]{figures/Hall_reconstructed_foreground_f4.eps}
% 		\includegraphics[width=\widthHall]{figures/Hall_reconstructed_foreground_f100.eps}
% 		\includegraphics[width=\widthHall]{figures/Hall_reconstructed_foreground_f150.eps}
% 		\includegraphics[width=\widthHall]{figures/Hall_reconstructed_foreground_f250.eps}
% 		\hspace{0.35cm}

		% =====================================================================================
		% Volkan
		
		\def\widthHall{3.035cm}
		\includegraphics[width=\widthHall]{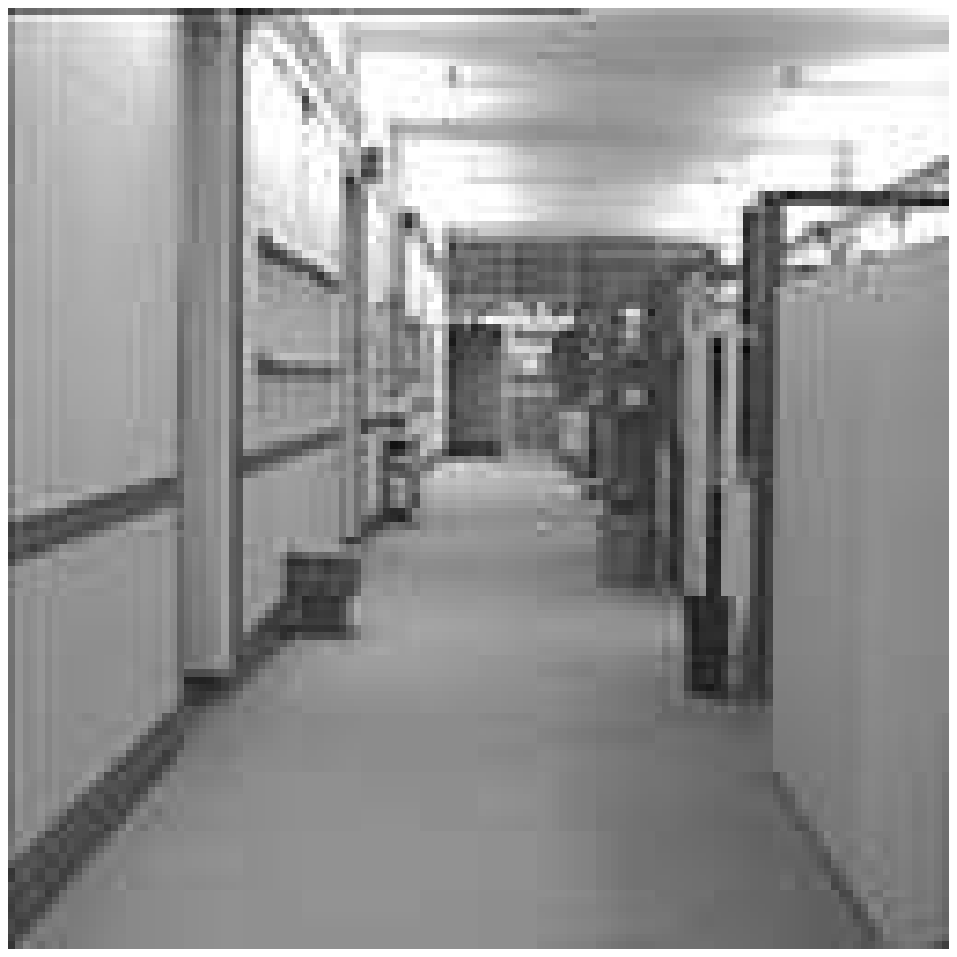}
		\includegraphics[width=\widthHall]{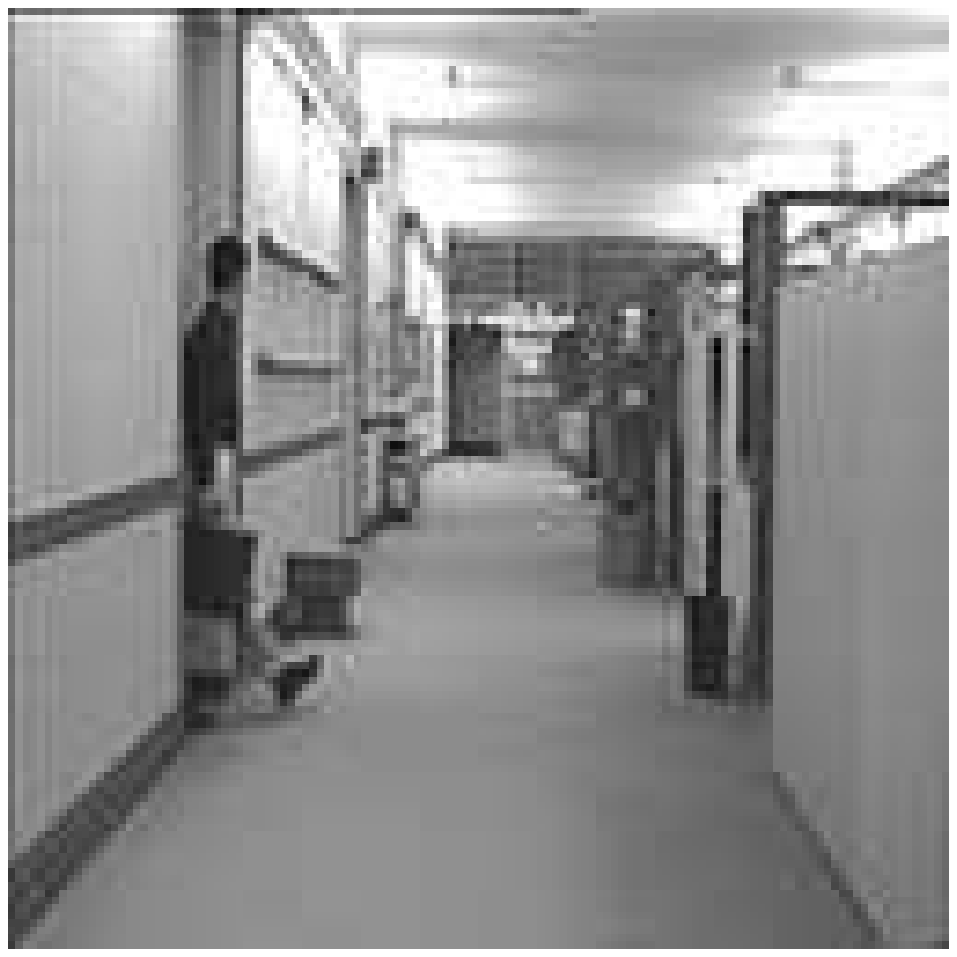}
		\includegraphics[width=\widthHall]{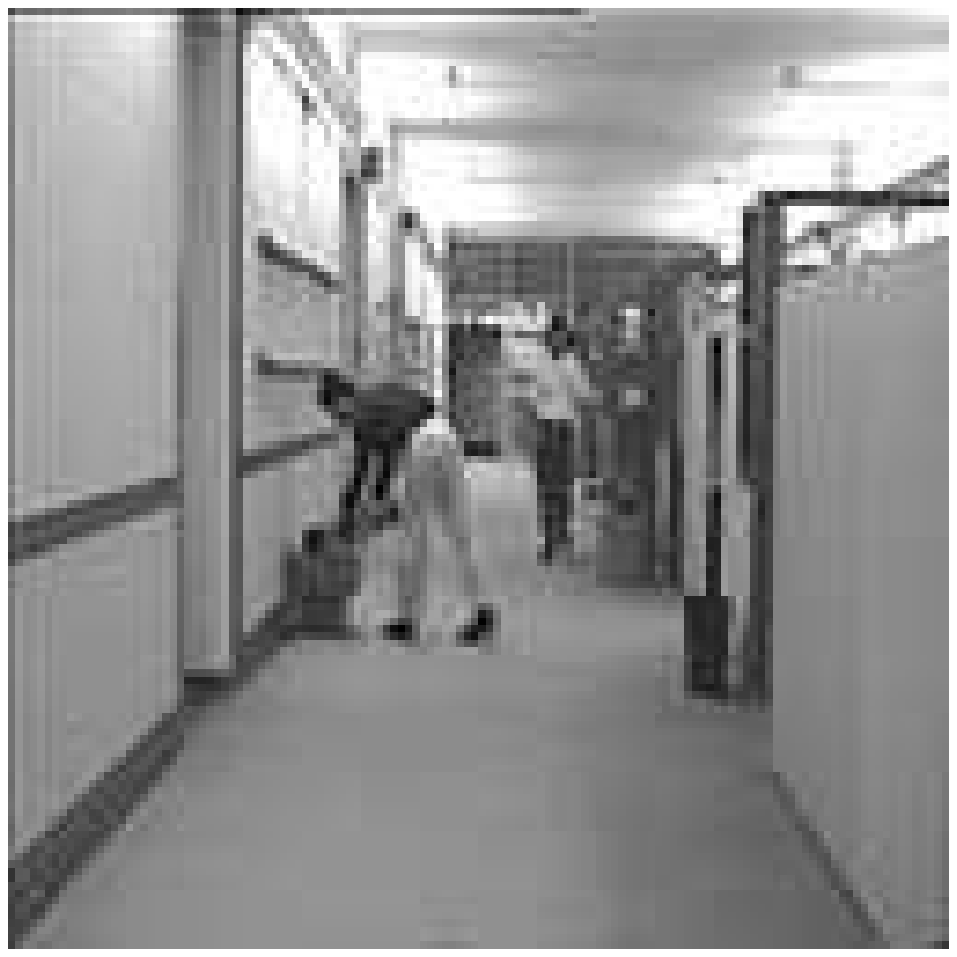}
		\includegraphics[width=\widthHall]{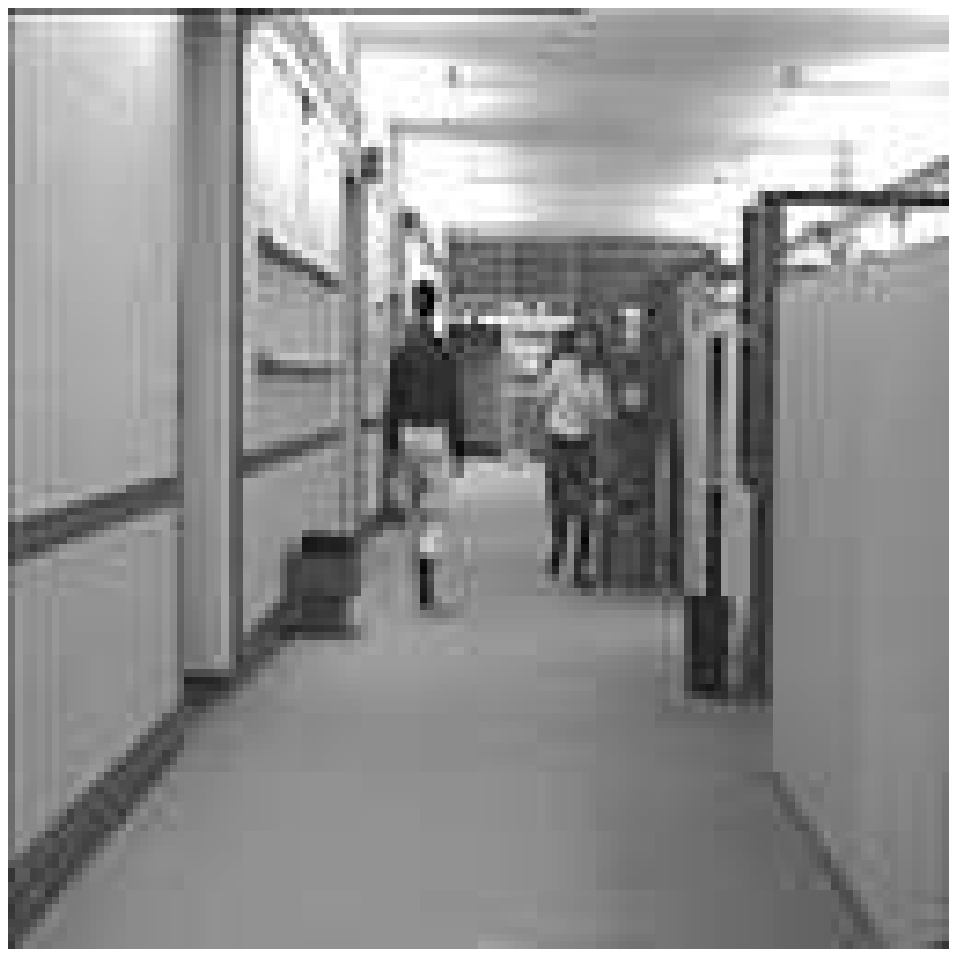}
		\includegraphics[width=\widthHall]{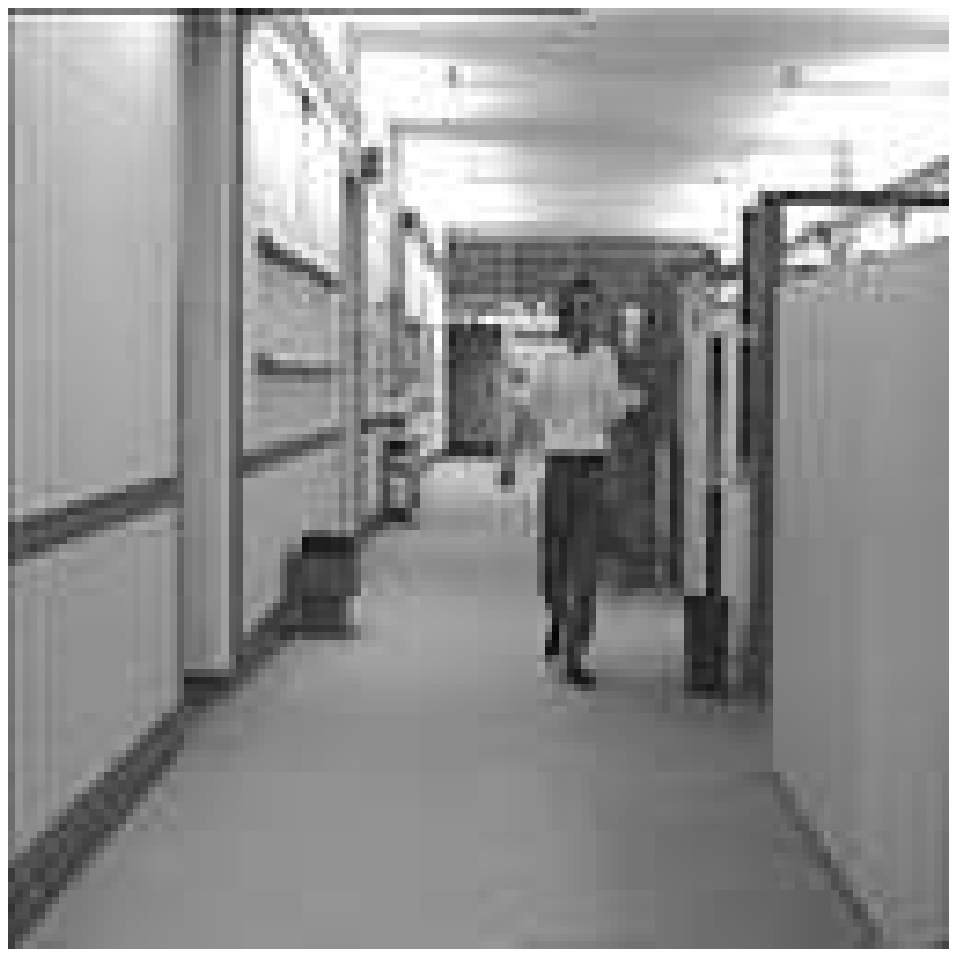}
		\hspace{0.35cm}
		
		\vspace{0.1cm}
		\includegraphics[width=\widthHall]{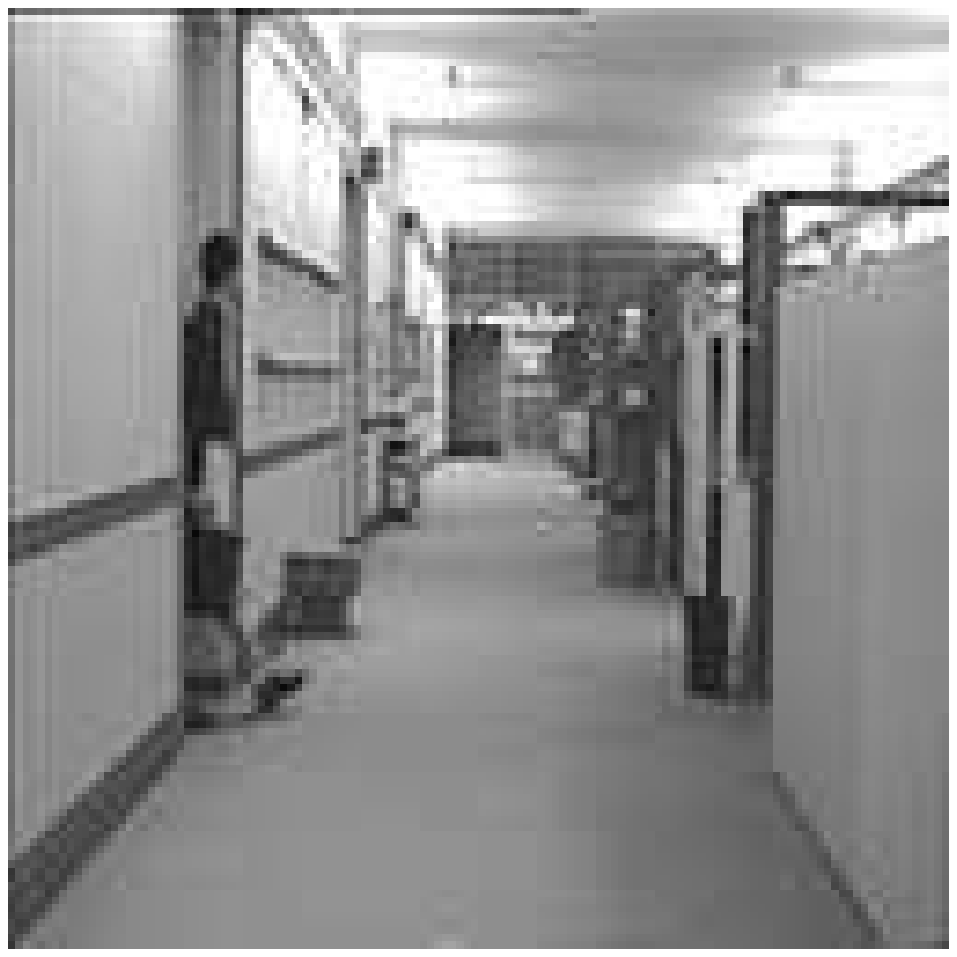}
		\includegraphics[width=\widthHall]{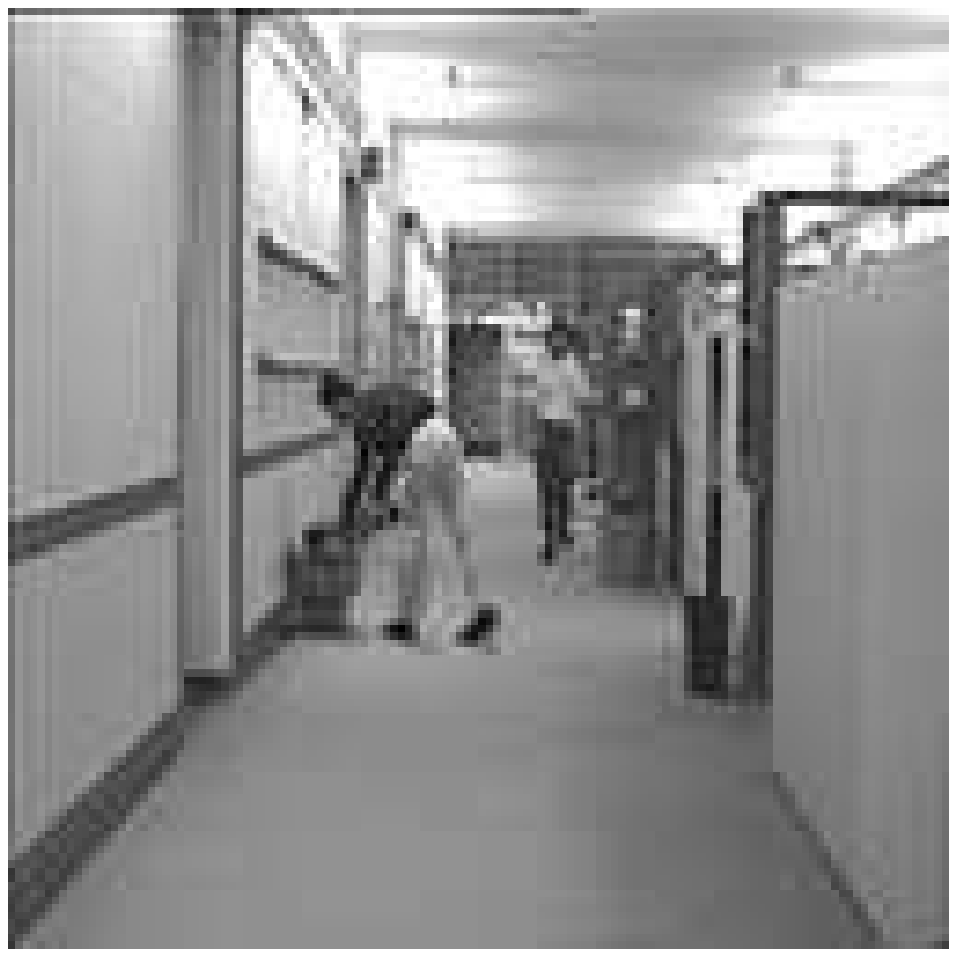}
		\includegraphics[width=\widthHall]{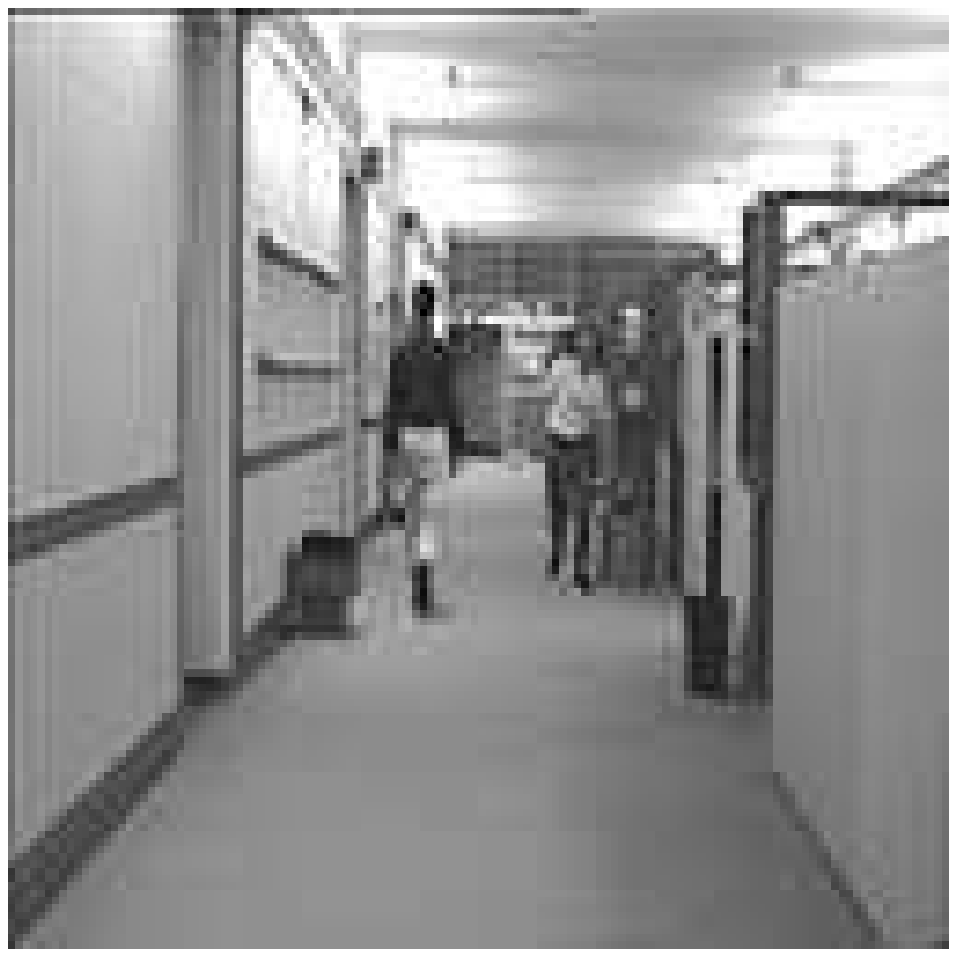}
		\includegraphics[width=\widthHall]{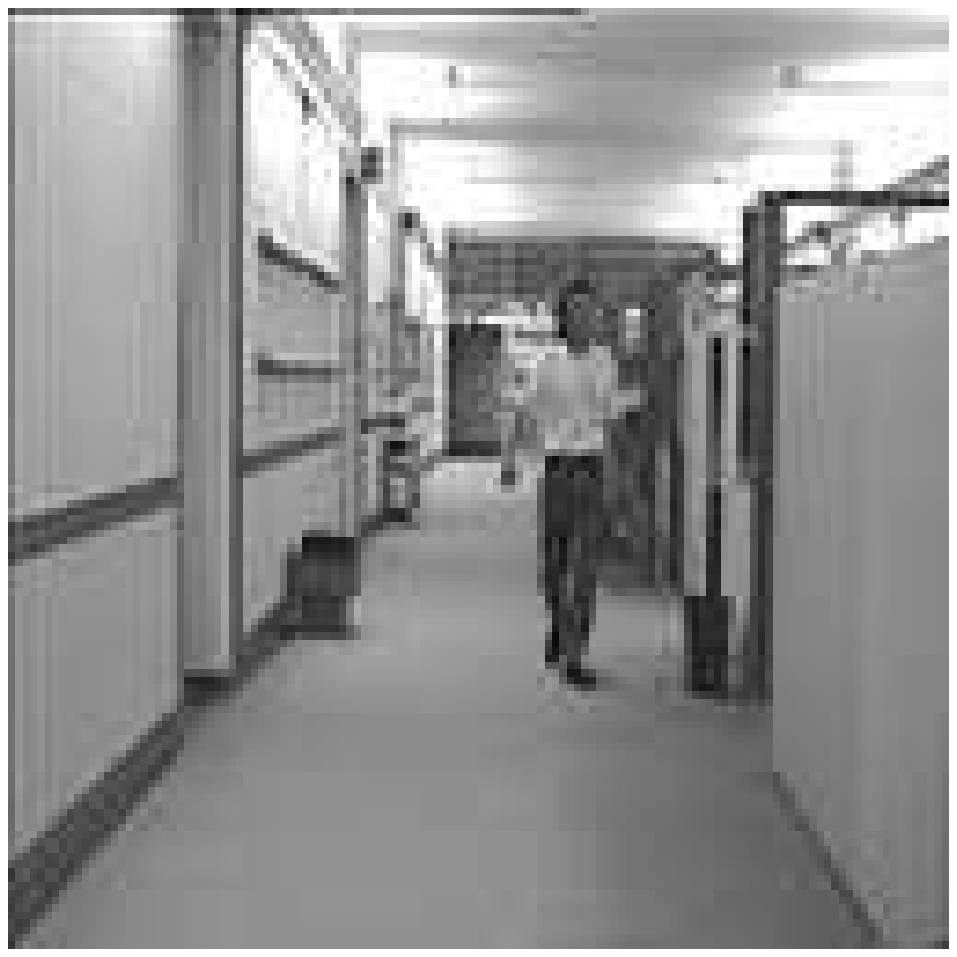}
		\hspace{0.35cm}
		
		\vspace{0.1cm}
		\includegraphics[width=\widthHall]{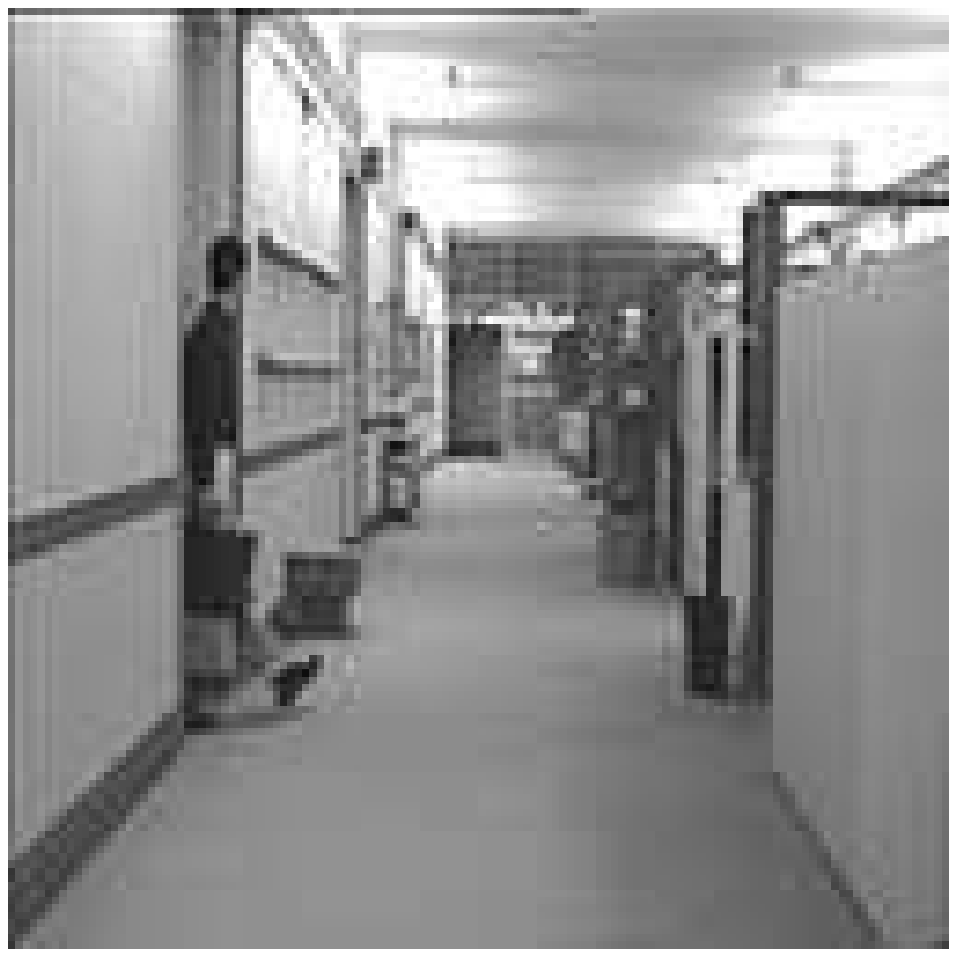}
		\includegraphics[width=\widthHall]{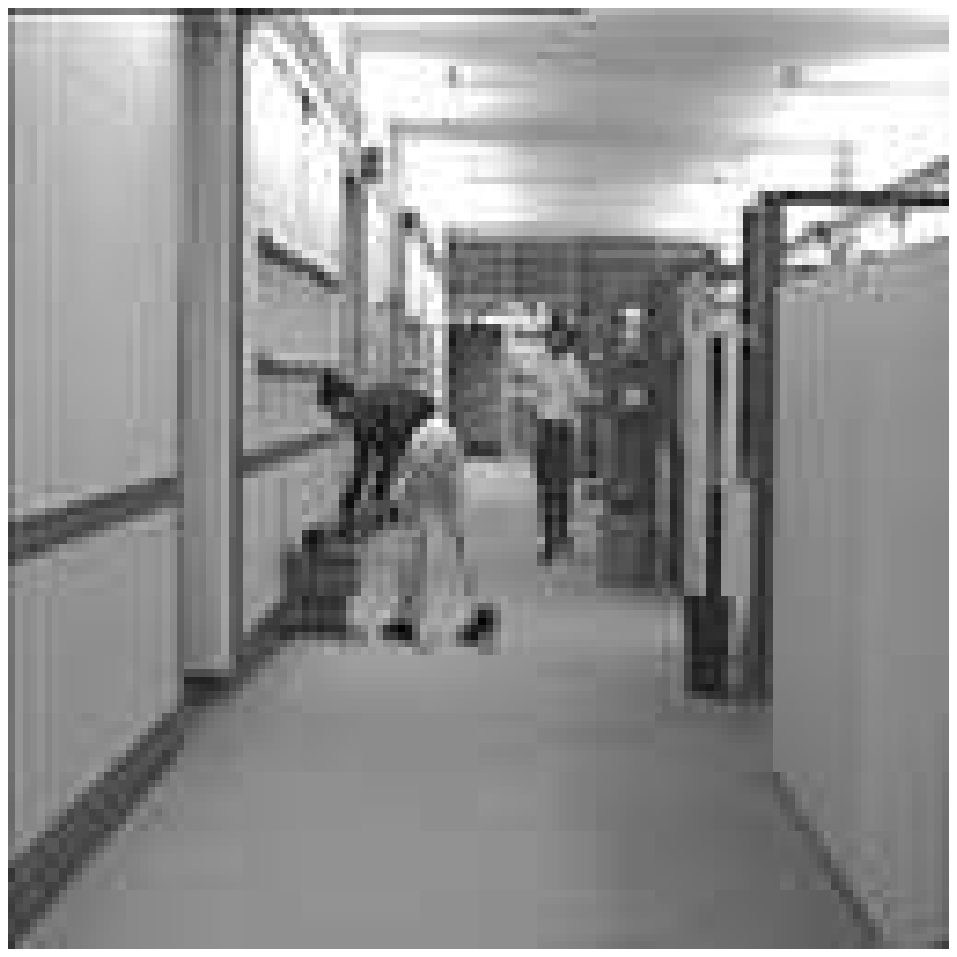}
		\includegraphics[width=\widthHall]{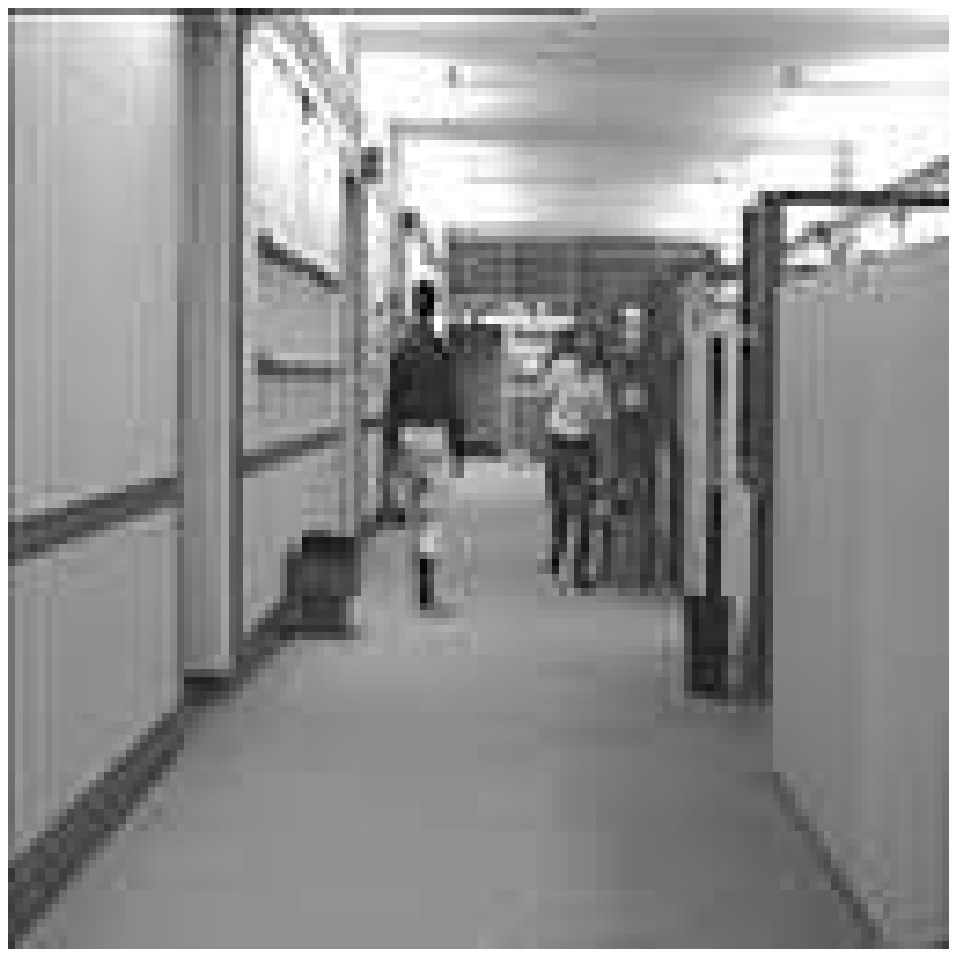}
		\includegraphics[width=\widthHall]{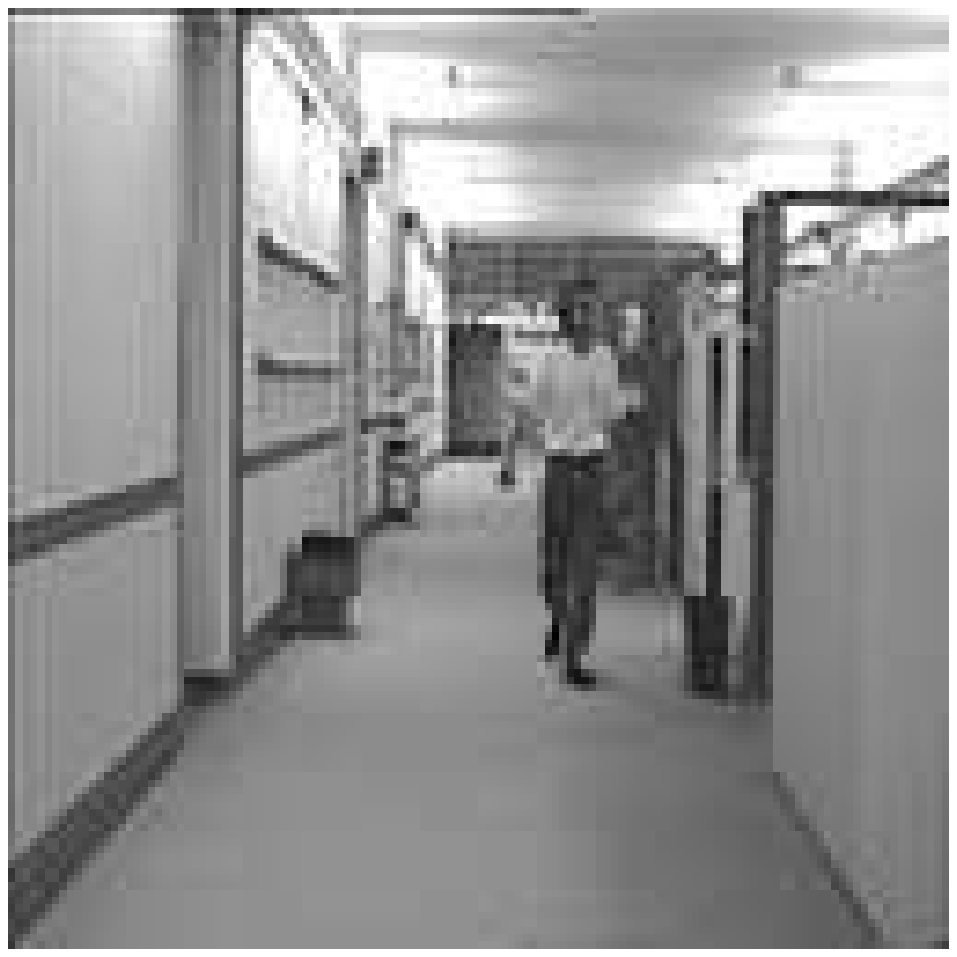}
		\hspace{0.35cm}
		
		\vspace{0.1cm}
		\includegraphics[width=\widthHall]{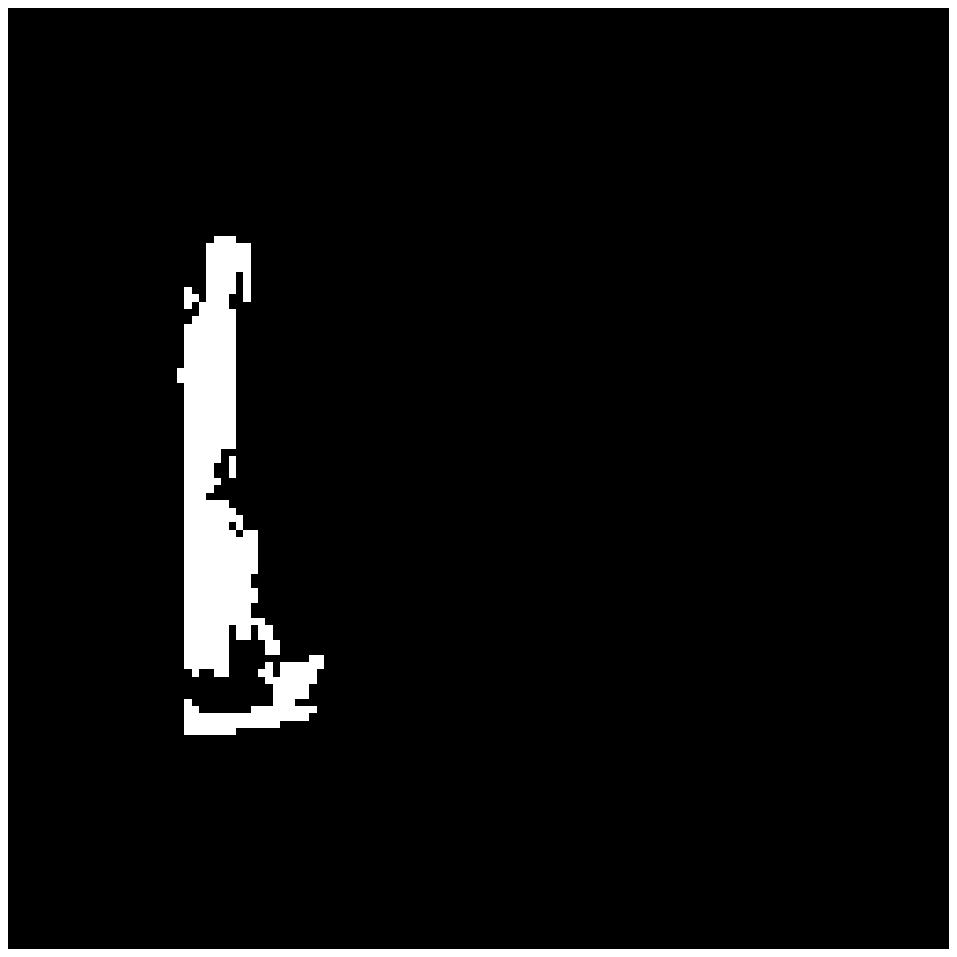}
		\includegraphics[width=\widthHall]{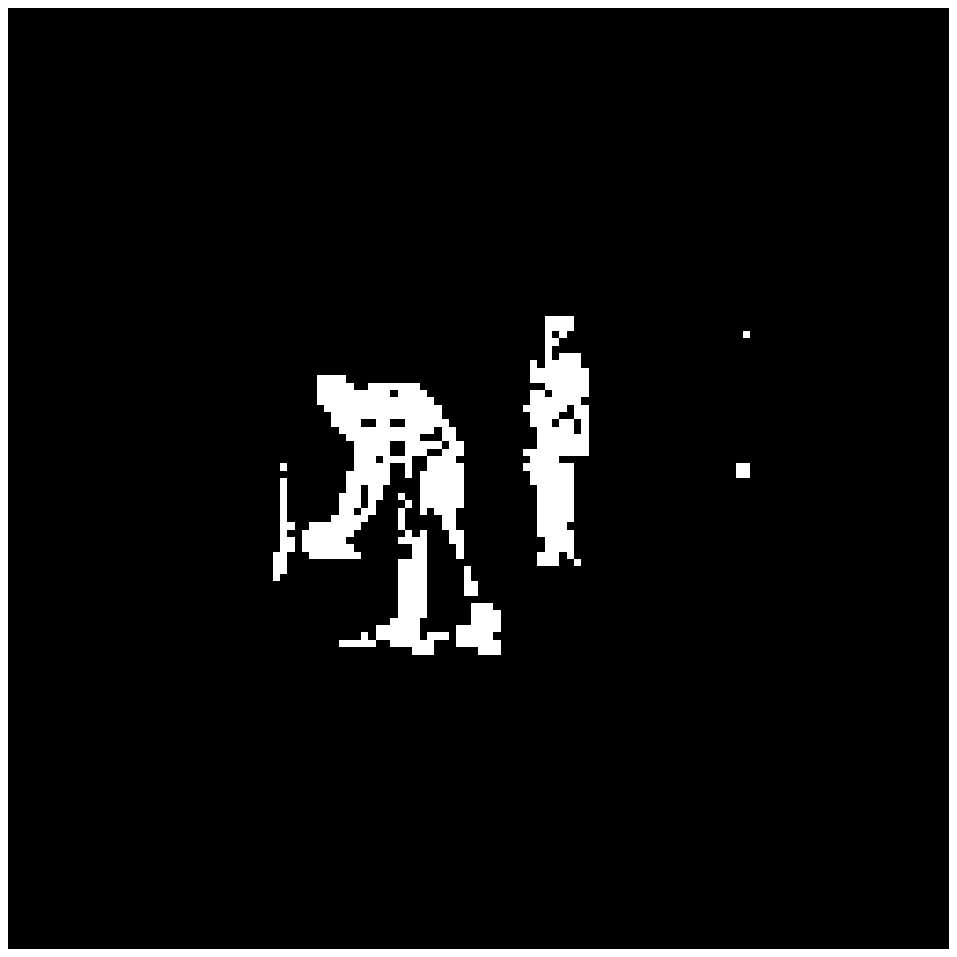}
		\includegraphics[width=\widthHall]{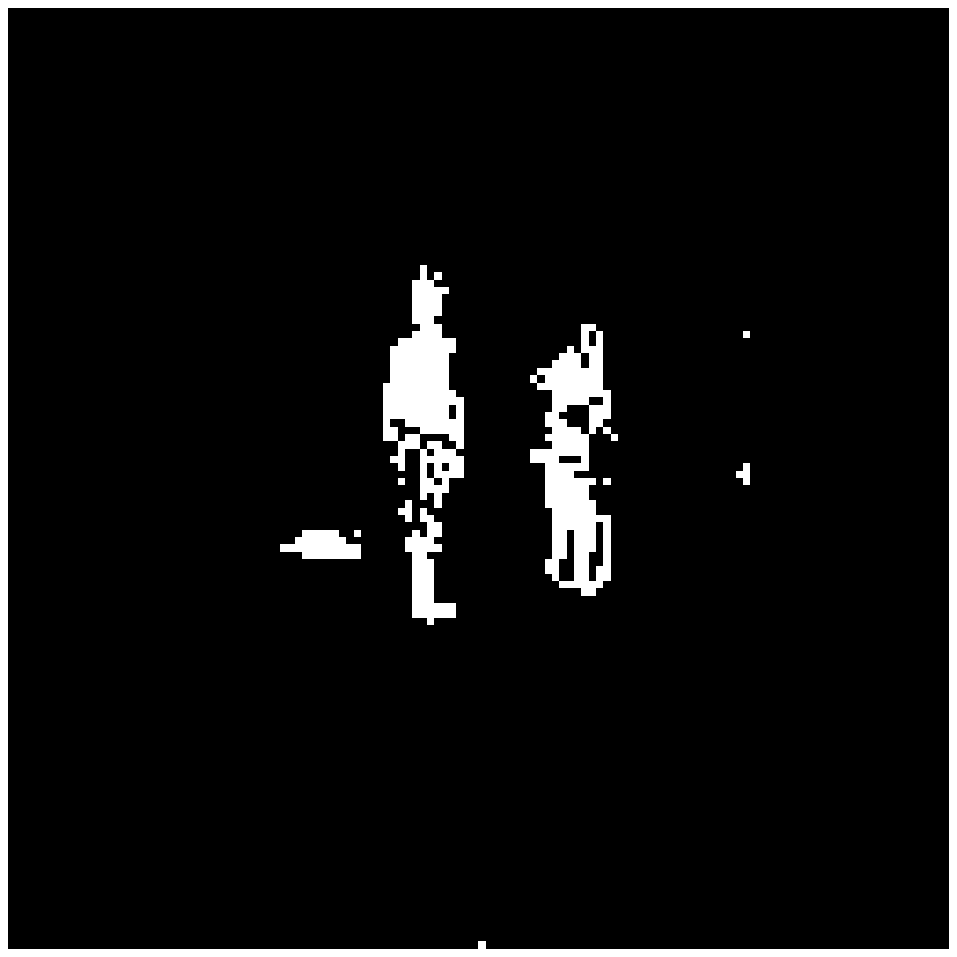}
		\includegraphics[width=\widthHall]{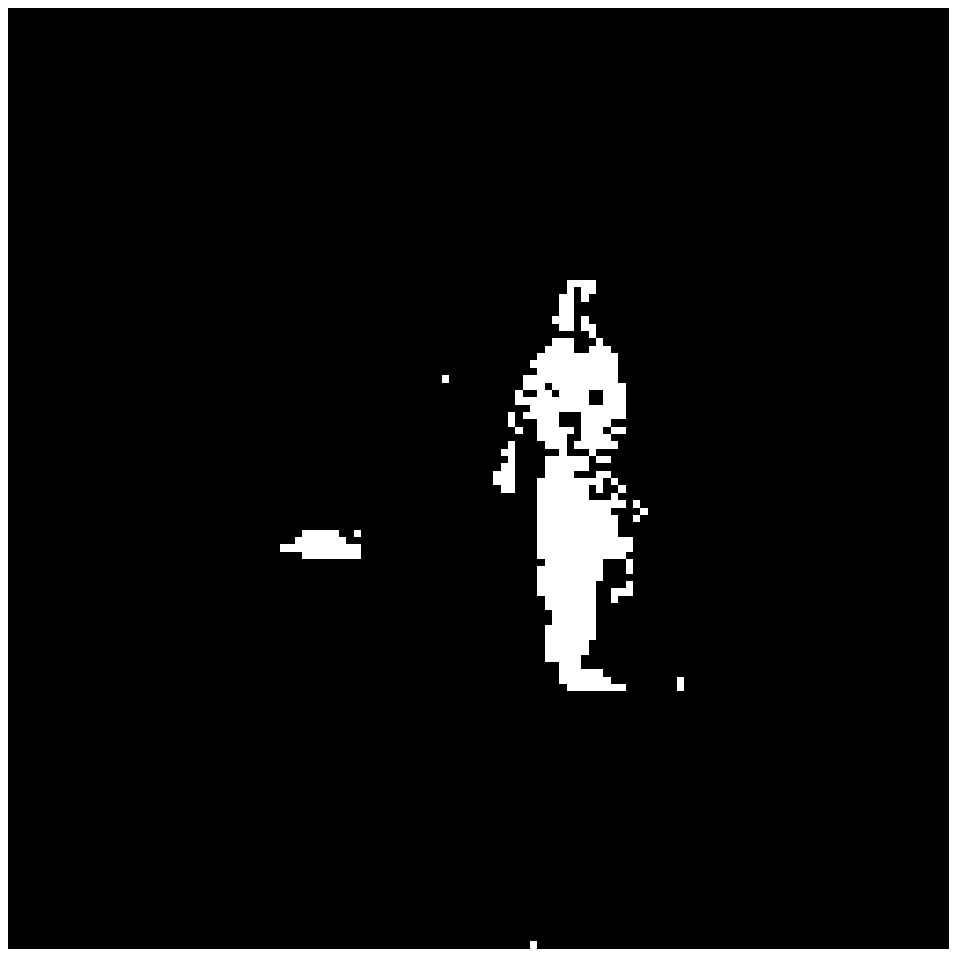}
		\hspace{0.35cm}
		% =====================================================================================
				
		\caption{
		    Hall sequence. The top panel shows the background and~$4$ different frames of the original
            sequence, which consists of~$282$ frames. The remaining panels show the estimated
            frames~$e[k]$, the reconstructed frame~$\hat{z}[k]$, and the reconstructed
            foreground~$\hat{x}[k]$ (binarized for better visualization).
		}
		\label{Fig:HallSequence}
		
	\end{figure*}

	\begin{figure*}
	\centering
	
	\subfigure[Quantities related to the number of measurements.]{\label{SubFig:HallMeasurements}
% 	\readdata{\dataCSOracle}{figures/HallCSmeasurements_oracle.dat}
% 	\readdata{\dataLLMinOracle}{figures/HallL1L1meas_oracle.dat}
% 	\readdata{\dataLLmeasEstim}{figures/HallL1L1meas_estim.dat}
% 	\readdata{\dataMeasurements}{figures/HallMeasurements.dat}
	
	% =======================================================================
	% Volkan
	\readdata{\dataCSOracle}{figures/VolkanAlg/HallCSmeasurements_oracle.dat}
	\readdata{\dataLLMinOracle}{figures/VolkanAlg/HallL1L1meas_oracle.dat}
	\readdata{\dataLLmeasEstim}{figures/VolkanAlg/HallL1L1meas_estim.dat}
	\readdata{\dataMeasurements}{figures/VolkanAlg/HallMeasurements.dat}
	% =======================================================================
			
	\psscalebox{0.96}{
	\begin{pspicture}(8.8,5.1)
					
		\def\xMax{300}                                % Maximum value of x
		\def\xMin{0}                                  % Minimum value of x
		\def\xNumTicks{6}                             % Number of x ticks						
		\def\yMax{6000}                               % Maximum value of y
		\def\yMin{0}                                  % Minimum value of y
		\def\yNumTicks{6}                             % Number of y ticks						
		\def\xIncrement{50}                           % =(xMax-xMin)/xNumTicks						
		\def\yIncrement{1000}                         % =(yMax-yMin)/yNumTicks			
					
		\def\xOrig{0.60}                              % Origin of the plot: X
		\def\yOrig{0.80}                              % Origin of the plot: Y
		\def\SizeX{8.00}                              % Size of plot: horizontal
		\def\SizeY{3.80}                              % Size of plot: vertical									
		\def\xTickIncr{1.33}                          % = SizeX/xNumTicks
		\def\yTickIncr{0.64}                          % = SizeY/yNumTicks

		\input{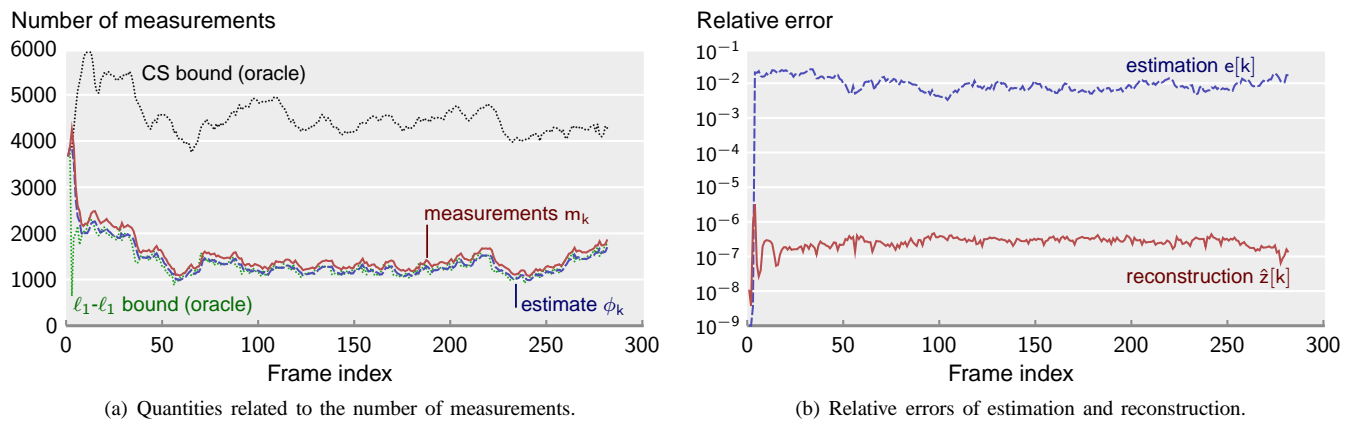}
		
		% ===========================================================================================
		% Display raw data
		
		\psset{xunit=\xScale\psunit,yunit=\yScale\psunit,linewidth=0.8pt}		
		
		\dataplot[origin={\xDataOrig,\yDataOrig},showpoints=false,linestyle=dotted,dotsep=0.6pt,linecolor=black!90!white]{\dataCSOracle}
		\dataplot[origin={\xDataOrig,\yDataOrig},showpoints=false,linestyle=dotted,dotsep=0.6pt,linecolor=green!90!white]{\dataLLMinOracle}	
		\dataplot[origin={\xDataOrig,\yDataOrig},showpoints=false,linestyle=dashed,dash=4pt 1pt,linecolor=blue!70!white]{\dataLLmeasEstim}
		\dataplot[origin={\xDataOrig,\yDataOrig},showpoints=false,linestyle=solid,dotsep=0.8pt,linecolor=red!70!white]{\dataMeasurements}
				
		\psset{xunit=\psunit,yunit=\psunit}
		% ===========================================================================================
		
		\rput[lb](-0.16,4.92){\small \textbf{\sf Number of measurements}}
		\rput[ct](4.25,0.14){\small \textbf{\sf Frame index}}
			
		\rput[lt](5.56,2.45){\footnotesize {\sf \color{red!70!black}{measurements $\mathsf{m_k}$}}}
		\psline[linewidth=0.7pt,linecolor=red!70!black](5.6,1.74)(5.6,2.2)
		
		\rput[lb](6.9,0.95){\footnotesize {\sf \color{blue!70!black}{estimate $\mathsf{\phi_k}$}}}
		\psline[linewidth=0.7pt,linecolor=blue!70!black](6.83,1.05)(6.83,1.38)
		
		\rput[lb](0.7,0.95){\footnotesize {\sf \color{green!70!black}{$\mathsf{\ell_1}$-$\mathsf{\ell_1}$ bound (oracle)}}}
		
		\rput[lb](1.65,4.15){\footnotesize {\sf CS bound (oracle)}}
		%\psgrid
	\end{pspicture}
	}
	\hfill
	}
	\subfigure[Relative errors of estimation and reconstruction.]{\label{SubFig:HallErrors}
	%\readdata{\estimError}{figures/Hall_EstimError.dat}
	%\readdata{\recError}{figures/Hall_RecError.dat}	
	% ========================================================
	% Volkan
 	\readdata{\estimError}{figures/VolkanAlg/Hall_EstimError.dat}
 	\readdata{\recError}{figures/VolkanAlg/Hall_RecError.dat}	
	% ========================================================
			
	\psscalebox{0.96}{
	\begin{pspicture}(8.8,5.1)
					
		\def\xMax{300}                                % Maximum value of x
		\def\xMin{0}                                  % Minimum value of x
		\def\xNumTicks{6}                             % Number of x ticks						
		\def\yMax{8}                                  % Maximum value of y
		\def\yMin{0}                                  % Minimum value of y
		\def\yNumTicks{8}                             % Number of y ticks						
		\def\xIncrement{50}                           % =(xMax-xMin)/xNumTicks						
		\def\yIncrement{1}                            % =(yMax-yMin)/yNumTicks			
					
		\def\xOrig{0.60}                              % Origin of the plot: X
		\def\yOrig{0.80}                              % Origin of the plot: Y
		\def\SizeX{8.00}                              % Size of plot: horizontal
		\def\SizeY{3.80}                              % Size of plot: vertical									
		\def\xTickIncr{1.33}                          % = SizeX/xNumTicks
		\def\yTickIncr{0.48}                          % = SizeY/yNumTicks

		%\input{auxFiles/prettyDataPlots}
		
		% ===========================================================================================
		% PrettyDataPlots code: because of log, it has to be tweaked
		
		% ==============================================================================
		% Parameters

    % Colors
    \definecolor{colorXAxis}{gray}{0.55}           % x-Axis color
    \definecolor{colorBackground}{gray}{0.93}      % Background color (0.91)

    % Distance between x-labels and the x-axis
    \def \distXLabels{0.15}

    % Distance between y-labels and the y-axis
    \def \distYLabels{0.12}

    % Width of the x-axis ticks
    \def \xTickWidth{0.08}
    % ==============================================================================

    \fpAdd{\xNumTicks}{1}{\xNumTicksPOne}         % = xNumTicks + 1
    \fpAdd{\yNumTicks}{1}{\yNumTicksPOne}

    \FPadd \xEndPoint \xOrig \SizeX                % xEndPoint = xOrig + SizeX
    \FPadd \yEndPoint \yOrig \SizeY                % yEndPoint = yOrig + SizeY

    \FPset \unit 1

    \FPsub \xRange \xMax \xMin			
    \FPdiv \xScale \SizeX \xRange 
    \FPdiv \xMultByOrigin \unit \xScale
    \FPmul \xDataOrig  \xMultByOrigin \xOrig

    \FPsub \yRange \yMax \yMin
    \FPdiv \yScale \SizeY \yRange 
    \FPdiv \yMultByOrigin \unit \yScale
    \FPmul \yDataOrig  \yMultByOrigin \yOrig

    \fpSub{\yOrig}{\distXLabels}{\xPosLabels}
    \fpSub{\xOrig}{\distYLabels}{\yPosLabels}

    \fpAdd{\yOrig}{\xTickWidth}{\xTickTop}

    % Background
    \psframe*[linecolor=colorBackground](\xOrig,\yOrig)(\xEndPoint,\yEndPoint)
                                        
    % x-labels and x-ticks: dx-position, iB-label			
    \multido{\nx=\xOrig+\xTickIncr, \iB=\xMin+\xIncrement}{\xNumTicksPOne}{	
        \rput[t](\nx,\xPosLabels){\small $\mathsf{\iB}$}
        \psline[linewidth=0.8pt,linecolor=colorXAxis](\nx,\yOrig)(\nx,\xTickTop)	
    }						

    % y-labels and lines			
    \multido{\ny=\yOrig+\yTickIncr, \nB=\yMin+\yIncrement}{\yNumTicksPOne}{	
        %\rput[r](\yPosLabels,\ny){\small $\mathsf{\nB}$}
        \psline[linecolor=white,linewidth=0.8pt](\xOrig,\ny)(\xEndPoint,\ny)
    }						

    % x-axis			
    \psline[linewidth=1.2pt,linecolor=colorXAxis]{-}(\xOrig,\yOrig)(\xEndPoint,\yOrig)
    % ===========================================================================================
		
		% Y-labels
		\rput[r](0.52,0.81){\footnotesize $\mathsf{10^{-9}}$}
		\rput[r](0.52,1.28){\footnotesize $\mathsf{10^{-8}}$}
		\rput[r](0.52,1.76){\footnotesize $\mathsf{10^{-7}}$}
		\rput[r](0.52,2.23){\footnotesize $\mathsf{10^{-6}}$}
		\rput[r](0.52,2.74){\footnotesize $\mathsf{10^{-5}}$}
		\rput[r](0.52,3.20){\footnotesize $\mathsf{10^{-4}}$}
		\rput[r](0.52,3.69){\footnotesize $\mathsf{10^{-3}}$}
		\rput[r](0.52,4.18){\footnotesize $\mathsf{10^{-2}}$}		
		\rput[r](0.52,4.62){\footnotesize $\mathsf{10^{-1}}$}
		            
		% ===========================================================================================
		% Display raw data
		
		\psset{xunit=\xScale\psunit,yunit=\yScale\psunit,linewidth=0.8pt}		
		
		\dataplot[origin={\xDataOrig,\yDataOrig},showpoints=false,linestyle=dashed,dash=4pt 1pt,linecolor=blue!70!white]{\estimError}
		
		\dataplot[origin={\xDataOrig,\yDataOrig},showpoints=false,linestyle=solid,dotsep=0.8pt,linecolor=red!70!white]{\recError}
				
		\psset{xunit=\psunit,yunit=\psunit}
		% ===========================================================================================
		
		\rput[lb](-0.1,4.92){\small \textbf{\sf Relative error}}
		\rput[ct](4.25,0.14){\small \textbf{\sf Frame index}}
			 		
 		%\rput[lb](5.85,2.75){\footnotesize {\sf \color{red!70!black}{reconstruction $\mathsf{\hat{z}[k]}$}}}
 		\rput[lt](5.85,1.60){\footnotesize {\sf \color{red!70!black}{reconstruction $\mathsf{\hat{z}[k]}$}}} % Volkan
 		\rput[lb](5.85,4.25){\footnotesize {\sf \color{blue!70!black}{estimation $\mathsf{e[k]}$}}}

 		%\psgrid
	\end{pspicture}
	}		
	}
	
	\caption{
		Results for the Hall sequence. \text{(a)} Number of measurements~$m_k$ taken from each frame (solid red line) and estimate~$\phi_k$ (dashed blue line); the dotted lines are the right-hand side of~\eqref{Eq:L1L1Bound} and~\eqref{Eq:ChandrasekaranBound} (green and black, respectively). 
		\text{(b)} Relative error of estimation $\|e[k] - z[k]\|_2/\|z[k]\|_2$, and reconstruction $\|\hat{z}[k] - z[k]\|_2/\|z[k]\|_2$. Figure \text{(b)} is illustrative, since the reconstruction error is mostly determined by the precision of the solver for~\eqref{Eq:L1L1Simple}.
	}
	\label{Fig:ResultsHall}
	\end{figure*}
% =================================

	\begin{figure*}
		\raggedleft
		
		\begin{pspicture}(18,0.2)
			%\rput[l](0.0,2.6){\small \sf Frame \#}
			\rput[l](0.0,-1.4){\sf Original}
			
			\rput[b](3.380,0.1){\sf background}
			\rput[b](6.580,0.16){\sf frame $\mathsf{5}$}
			\rput[b](9.750,0.16){\sf frame $\mathsf{75}$}
			\rput[b](12.92,0.16){\sf frame $\mathsf{100}$}
			\rput[b](16.05,0.16){\sf frame $\mathsf{170}$}

			\rput[lb](0.0,-4.8){\sf Estimated}
			
			\rput[lb](0.0,-7.9){\sf Reconstructed} 
			
			\rput[lb](0.0,-11.0){\sf Reconstructed foreground} 
			\rput[lb](0.0,-11.5){\sf (binarized)} 
			
			%\psgrid
		\end{pspicture}
		
% 		\def\widthPETS{3.035cm}
% 		\includegraphics[width=\widthPETS]{figures/PETS_original_f1.eps}
% 		\includegraphics[width=\widthPETS]{figures/PETS_original_f5.eps}
% 		\includegraphics[width=\widthPETS]{figures/PETS_original_f75.eps}
% 		\includegraphics[width=\widthPETS]{figures/PETS_original_f100.eps}
% 		\includegraphics[width=\widthPETS]{figures/PETS_original_f170.eps}
% 		\hspace{0.35cm}
% 		
% 		\vspace{0.1cm}
% 		\includegraphics[width=\widthPETS]{figures/PETS_estimated_image_f5.eps}
% 		\includegraphics[width=\widthPETS]{figures/PETS_estimated_image_f75.eps}
% 		\includegraphics[width=\widthPETS]{figures/PETS_estimated_image_f100.eps}
% 		\includegraphics[width=\widthPETS]{figures/PETS_estimated_image_f170.eps}
% 		\hspace{0.35cm}
% 		
% 		\vspace{0.1cm}
% 		\includegraphics[width=\widthPETS]{figures/PETS_reconstructed_image_f5.eps}
% 		\includegraphics[width=\widthPETS]{figures/PETS_reconstructed_image_f75.eps}
% 		\includegraphics[width=\widthPETS]{figures/PETS_reconstructed_image_f100.eps}
% 		\includegraphics[width=\widthPETS]{figures/PETS_reconstructed_image_f170.eps}
% 		\hspace{0.35cm}
% 		
% 		\vspace{0.1cm}
% 		\includegraphics[width=\widthPETS]{figures/PETS_reconstructed_foreground_f5.eps}
% 		\includegraphics[width=\widthPETS]{figures/PETS_reconstructed_foreground_f75.eps}
% 		\includegraphics[width=\widthPETS]{figures/PETS_reconstructed_foreground_f100.eps}
% 		\includegraphics[width=\widthPETS]{figures/PETS_reconstructed_foreground_f170.eps}
% 		\hspace{0.35cm}
		
		% =================================================================================
		% Volkan
		\def\widthPETS{3.035cm}
		\includegraphics[width=\widthPETS]{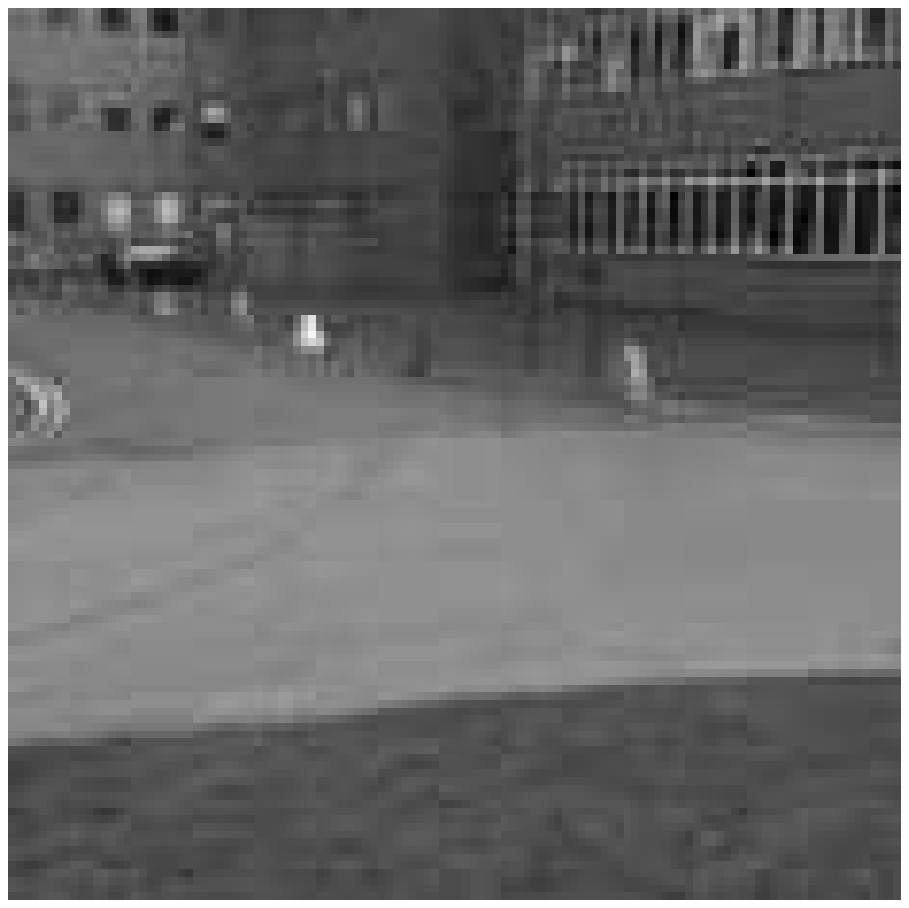}
		\includegraphics[width=\widthPETS]{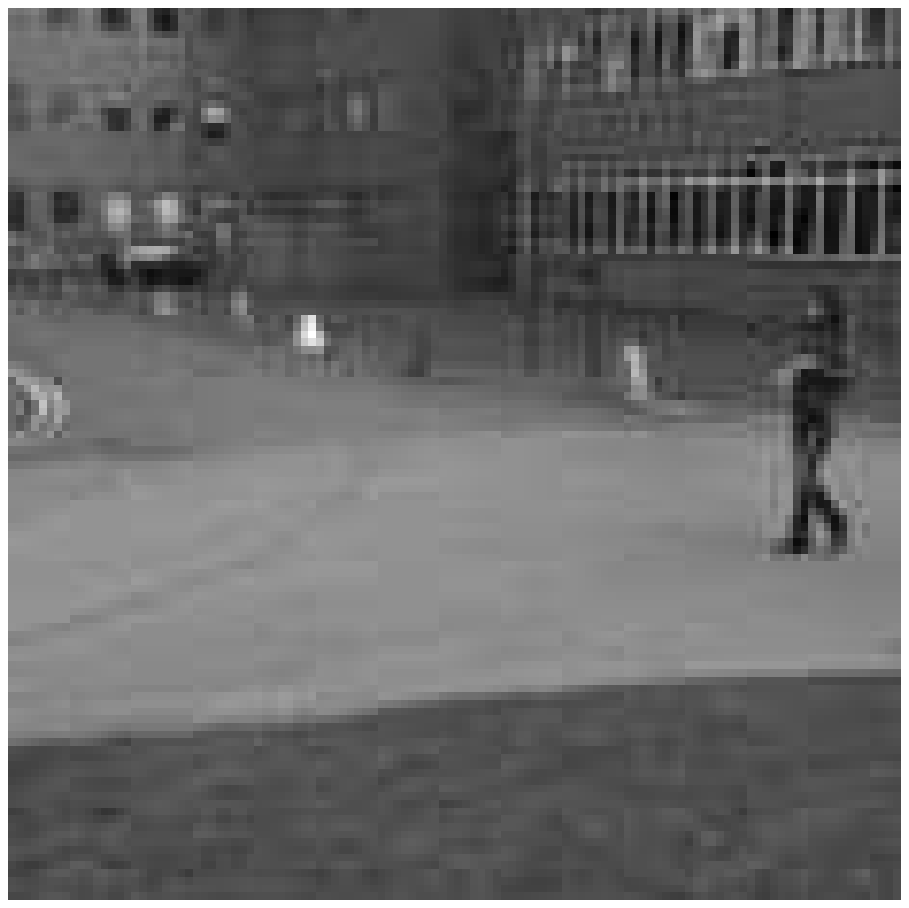}
		\includegraphics[width=\widthPETS]{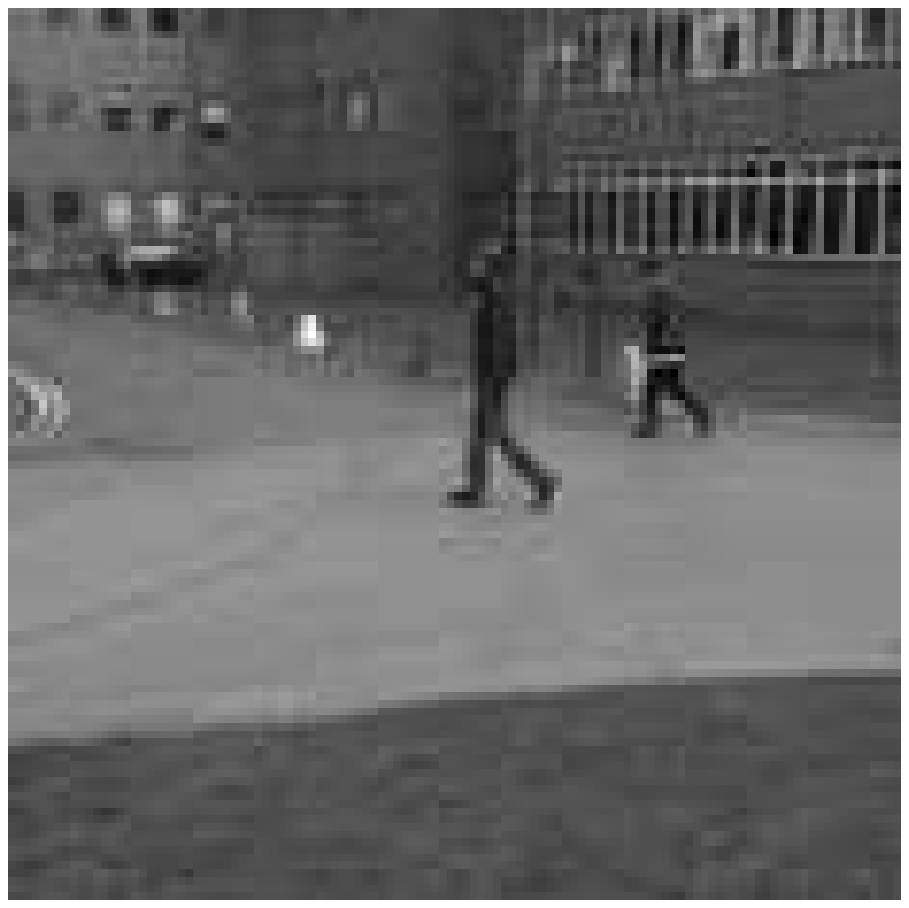}
		\includegraphics[width=\widthPETS]{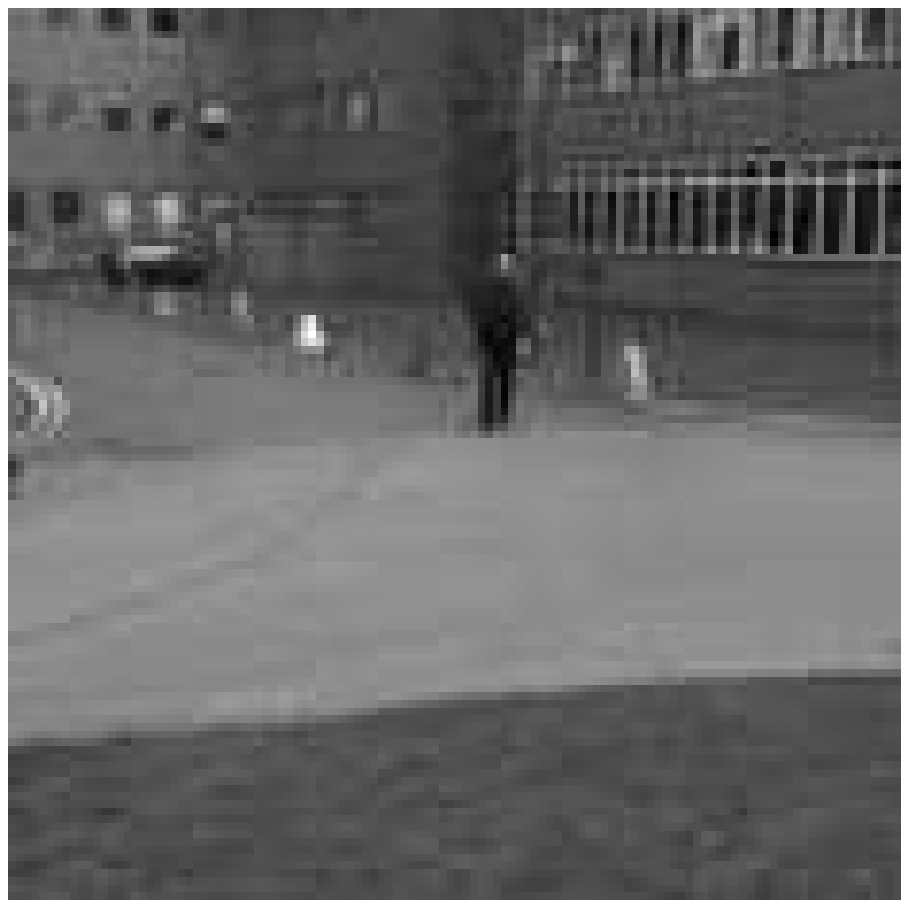}
		\includegraphics[width=\widthPETS]{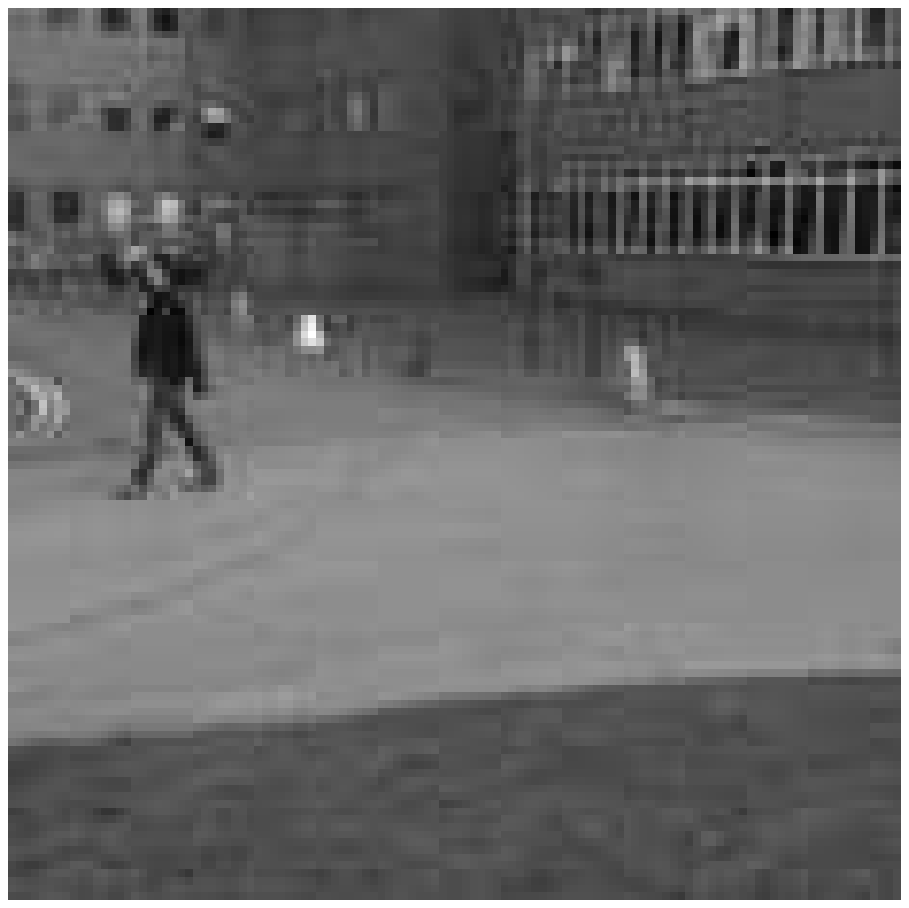}
		\hspace{0.35cm}
		
		\vspace{0.1cm}
		\includegraphics[width=\widthPETS]{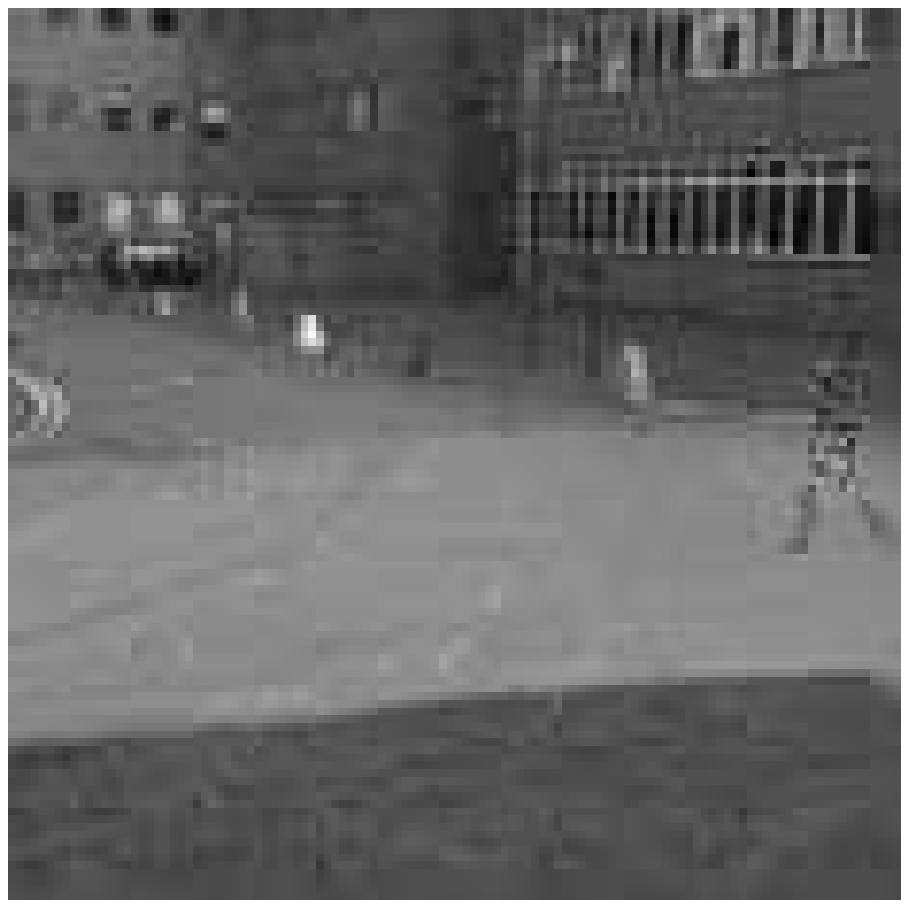}
		\includegraphics[width=\widthPETS]{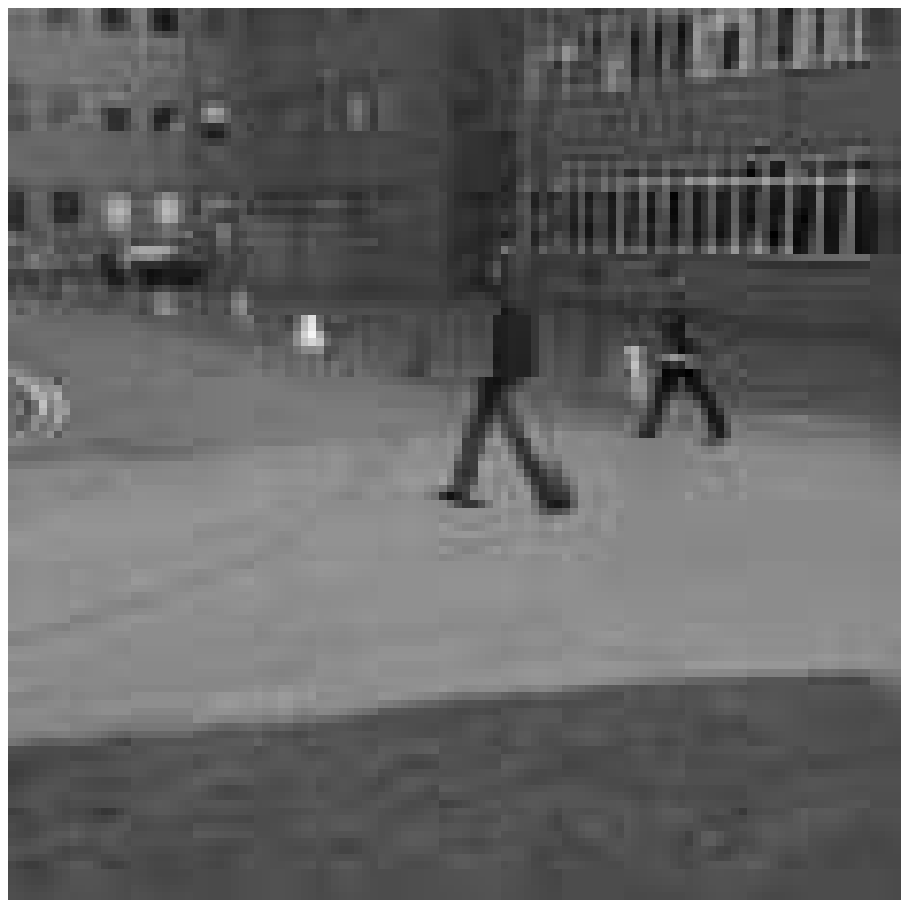}
		\includegraphics[width=\widthPETS]{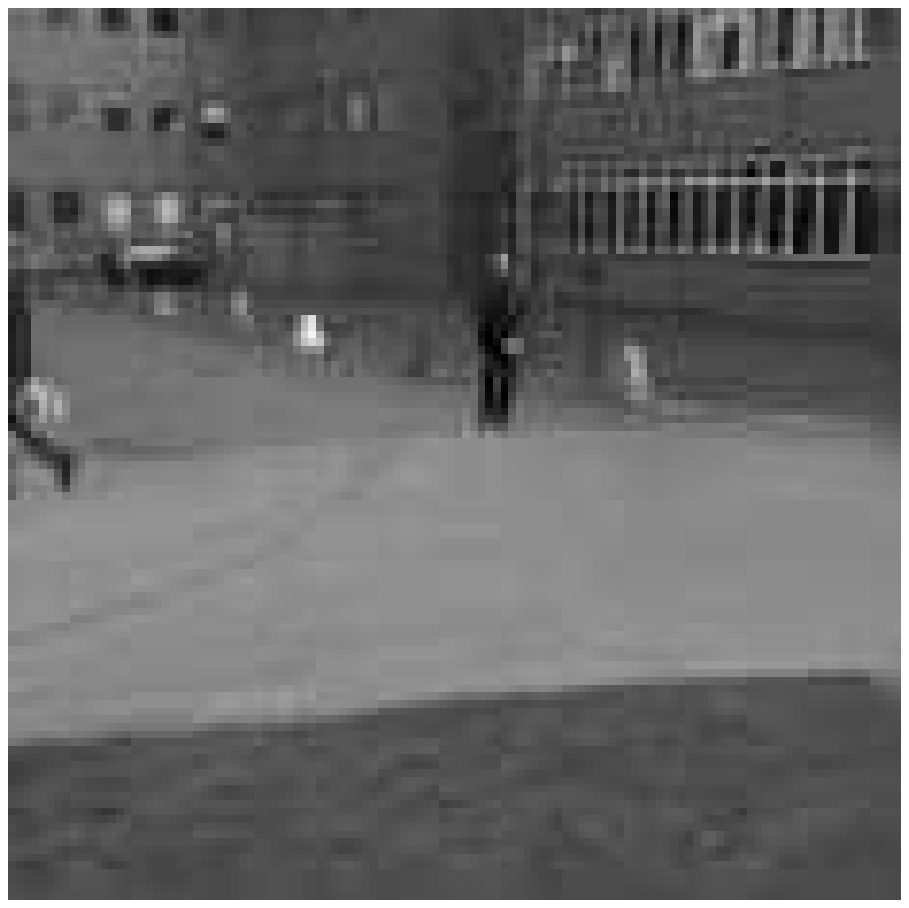}
		\includegraphics[width=\widthPETS]{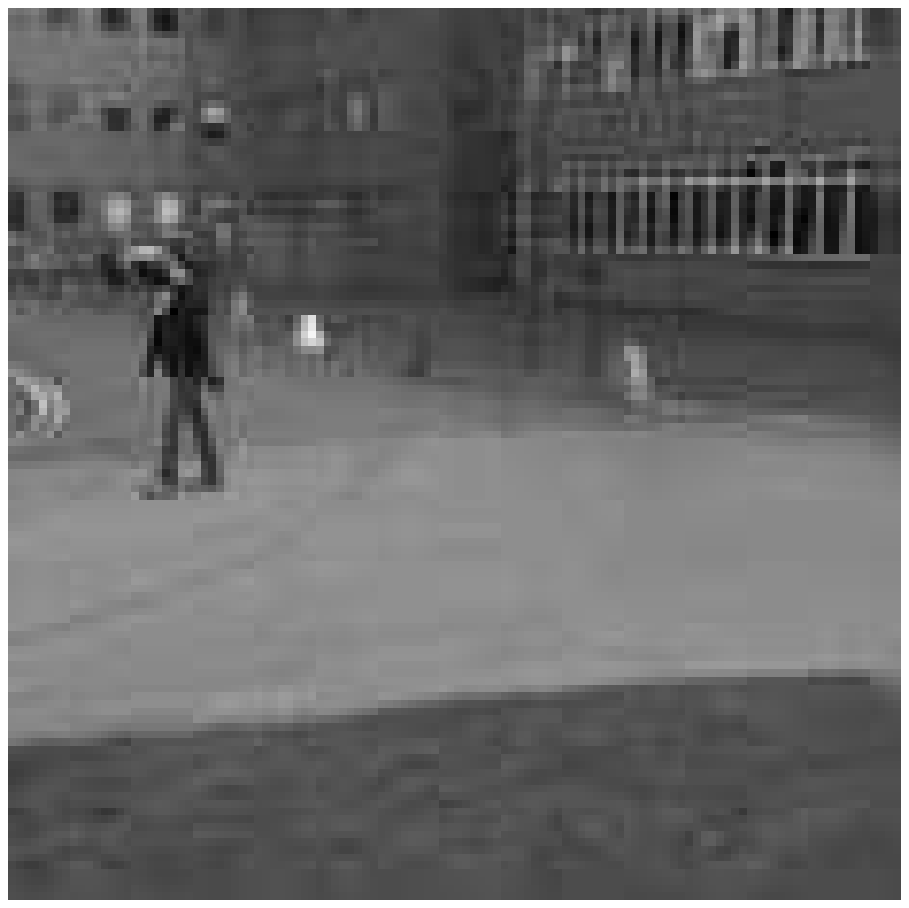}
		\hspace{0.35cm}
		
		\vspace{0.1cm}
		\includegraphics[width=\widthPETS]{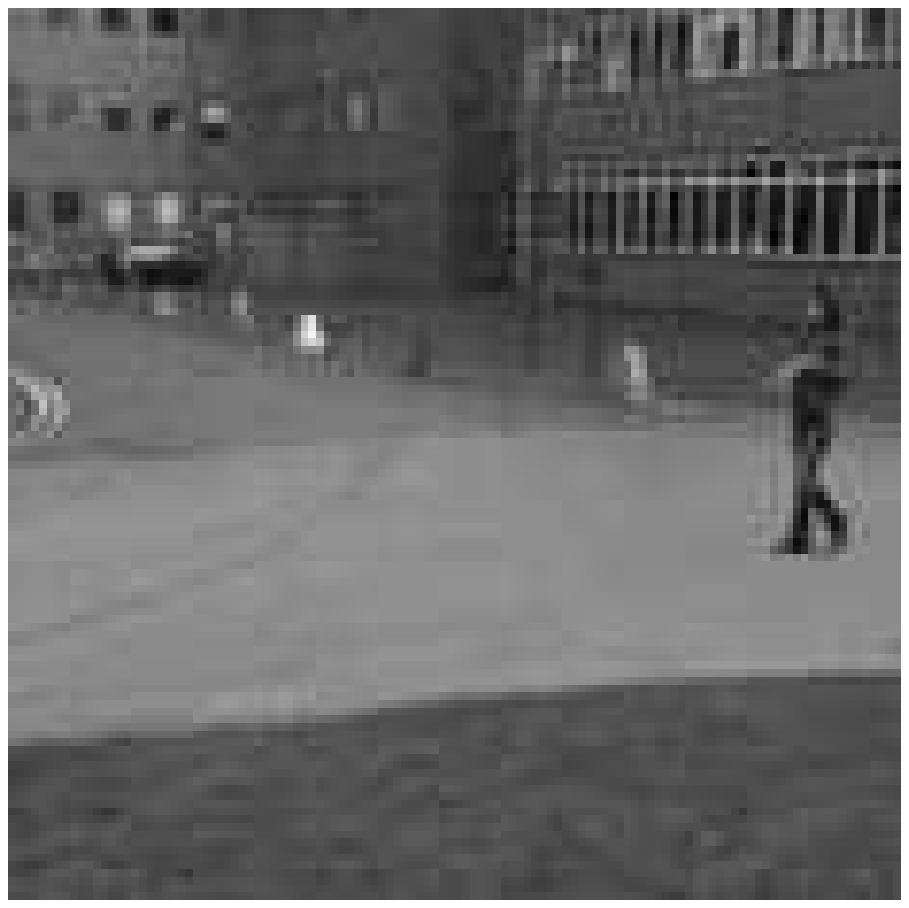}
		\includegraphics[width=\widthPETS]{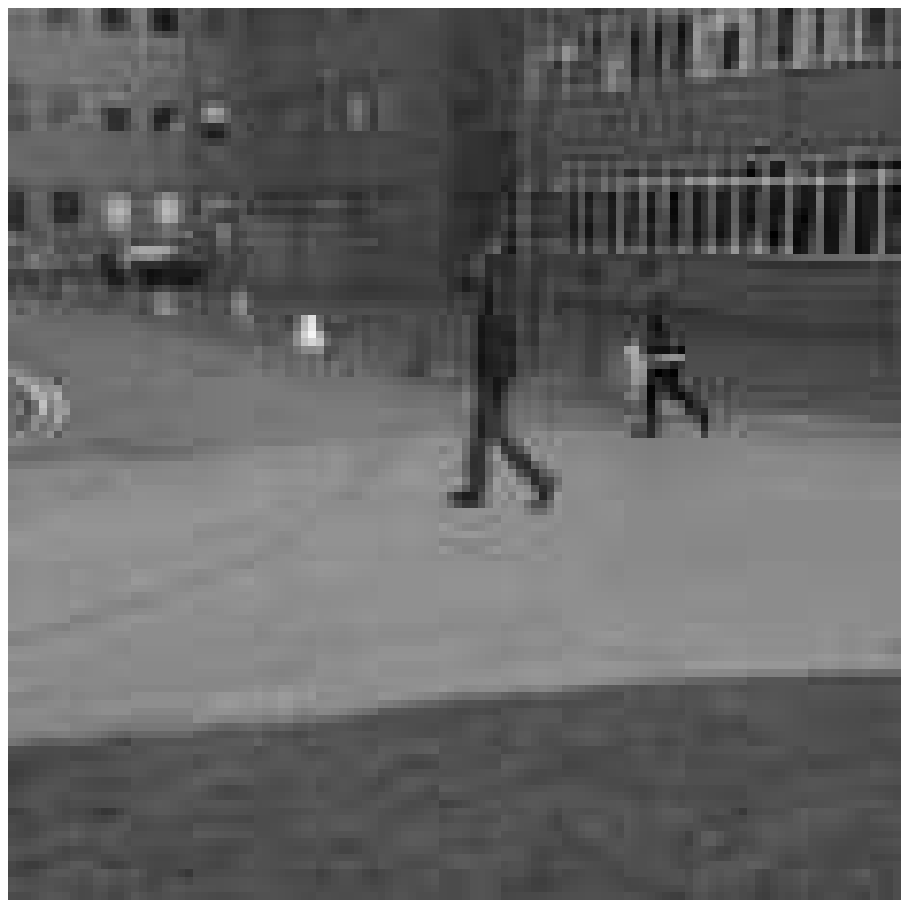}
		\includegraphics[width=\widthPETS]{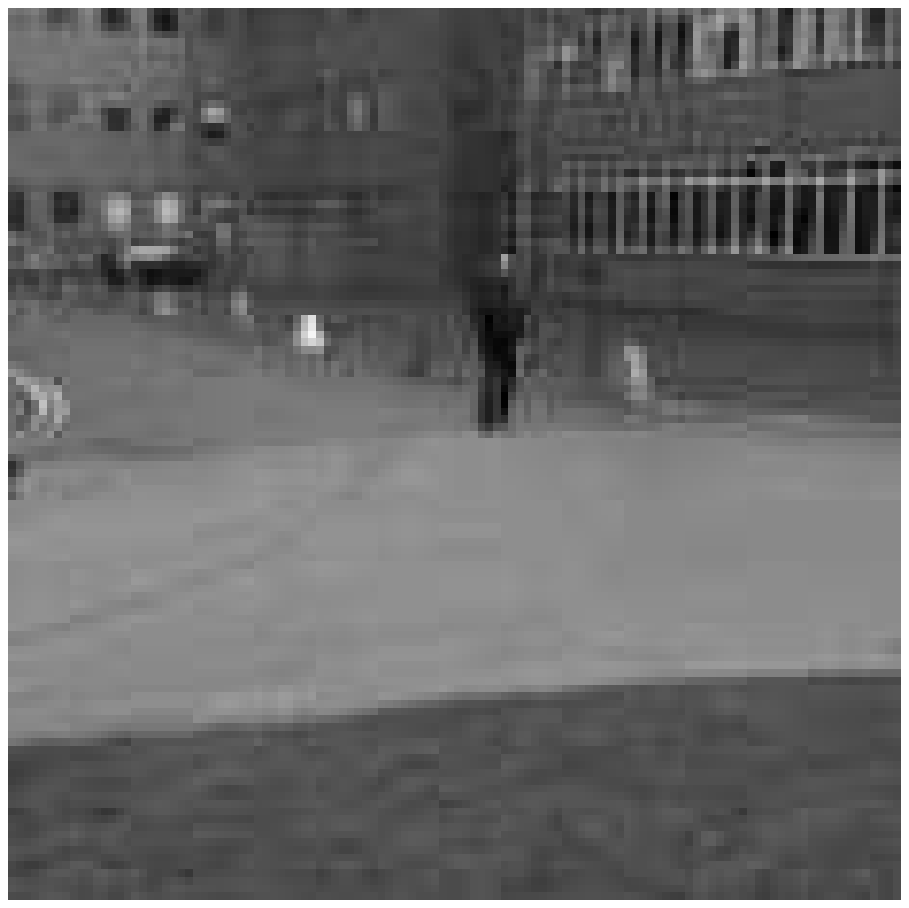}
		\includegraphics[width=\widthPETS]{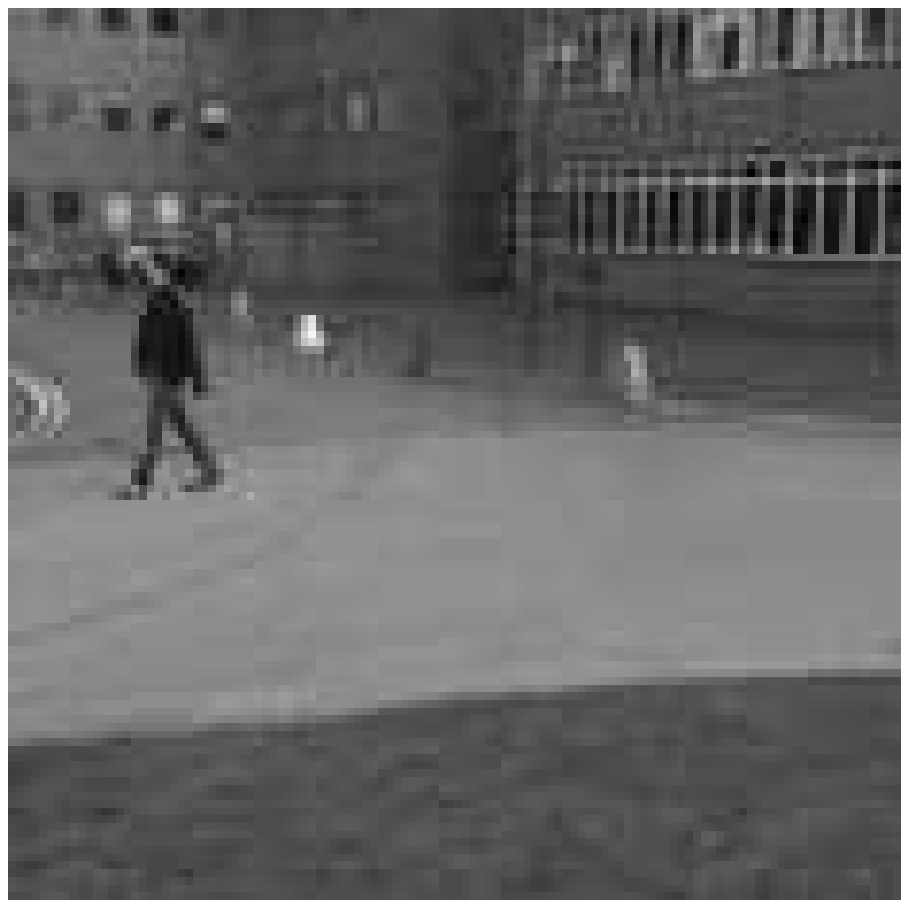}
		\hspace{0.35cm}
		
		\vspace{0.1cm}
		\includegraphics[width=\widthPETS]{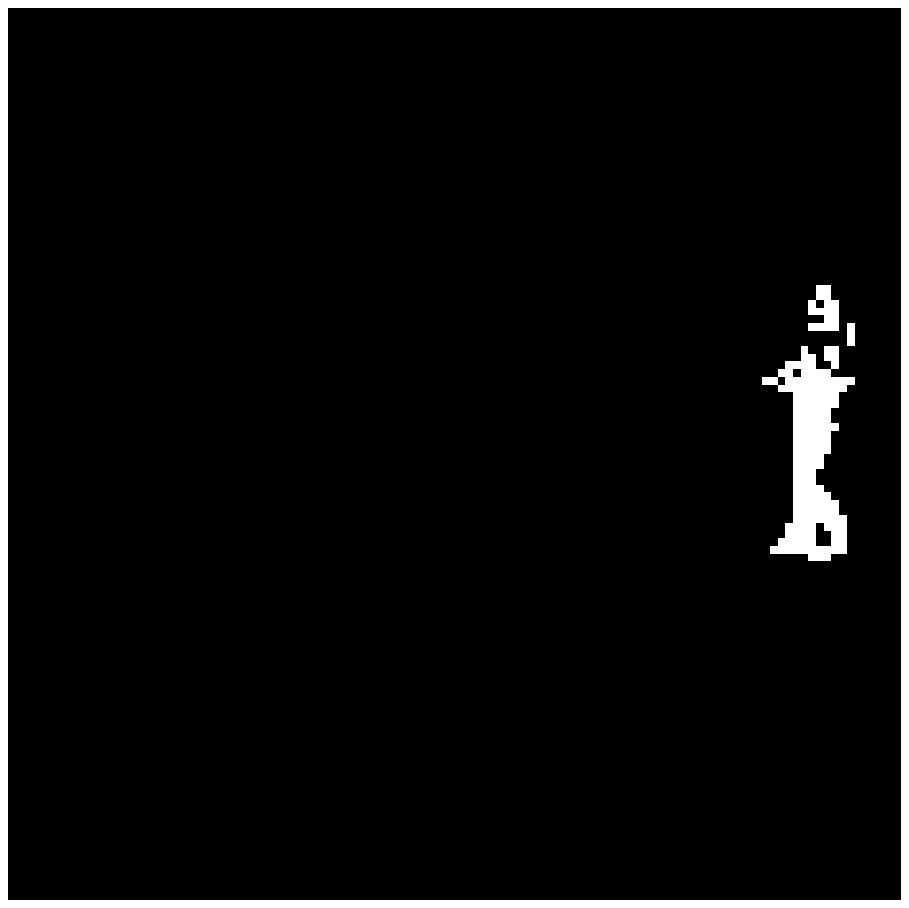}
		\includegraphics[width=\widthPETS]{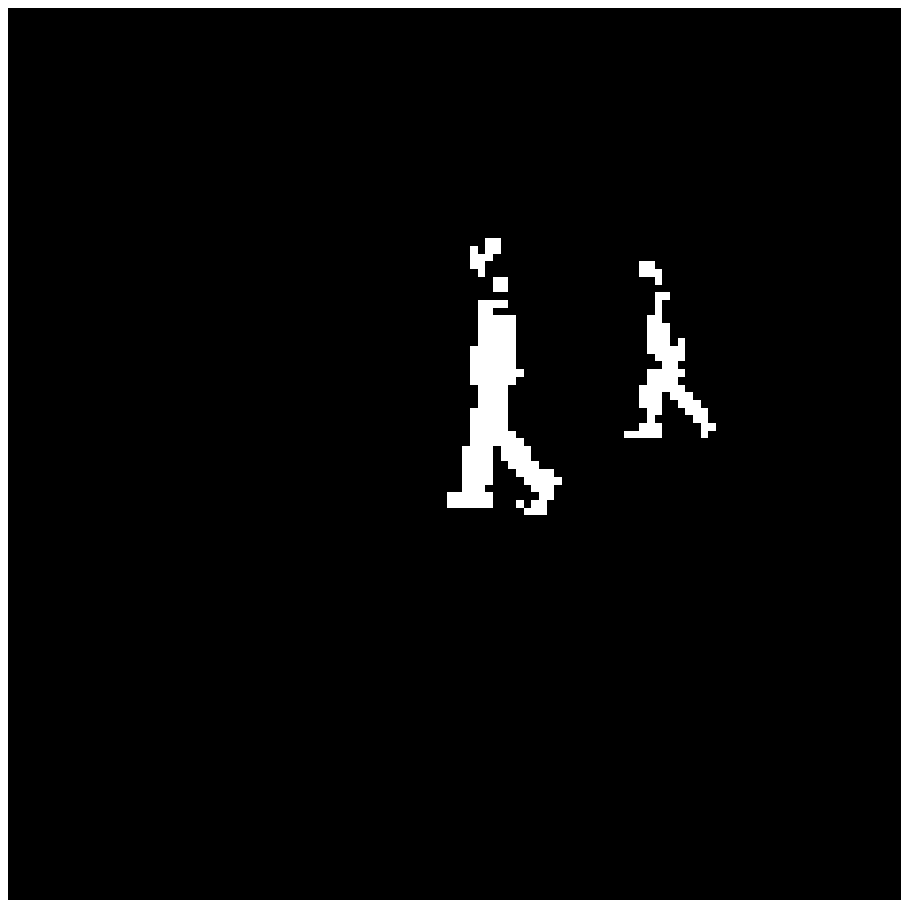}
		\includegraphics[width=\widthPETS]{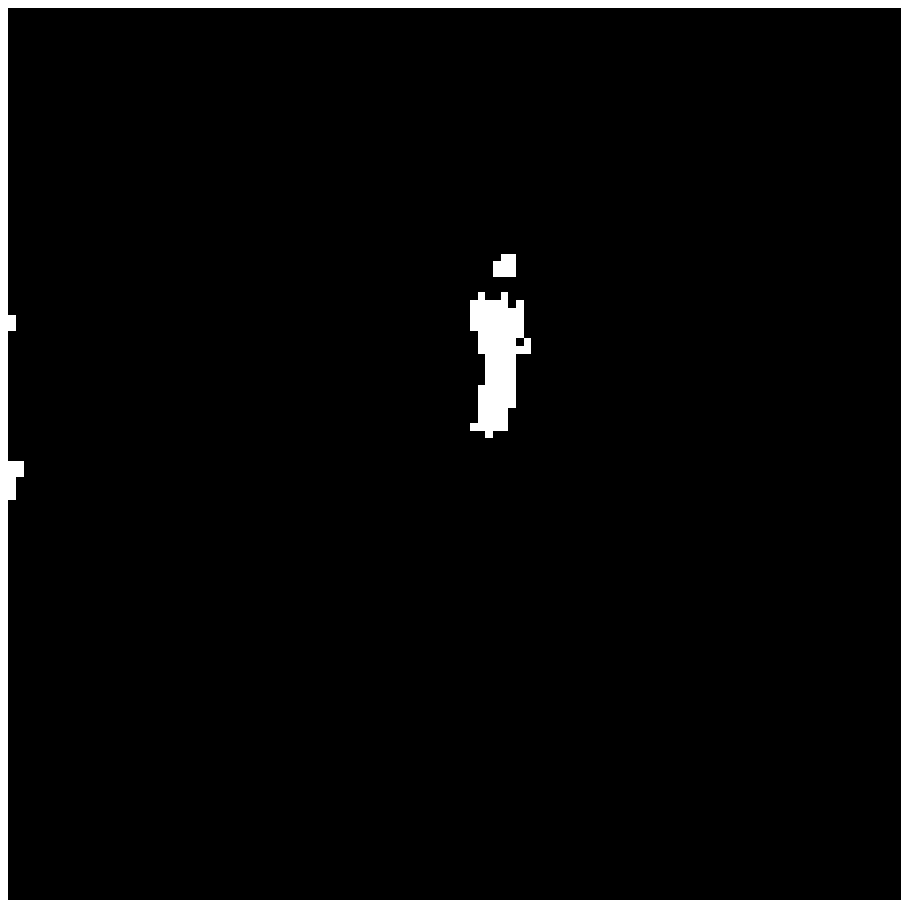}
		\includegraphics[width=\widthPETS]{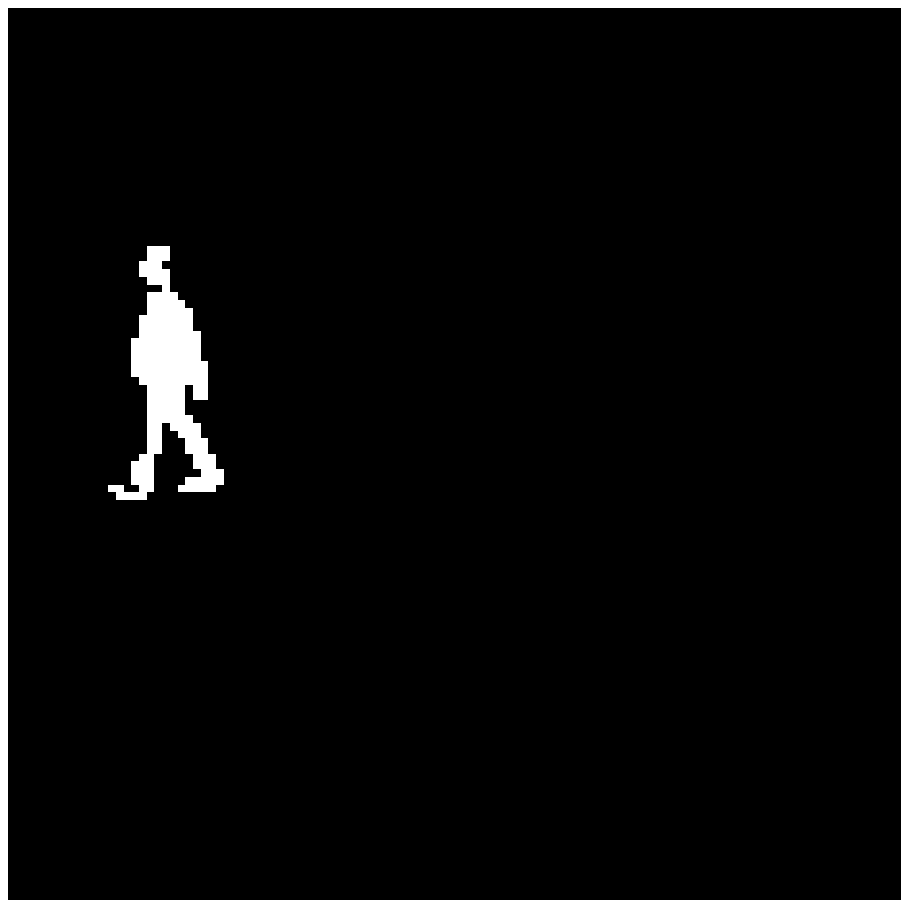}
		\hspace{0.35cm}
		% =================================================================================

		\caption{
      PETS sequence. The top panel shows the background and~$4$ different frames of the original sequence, which consists of~$171$ frames. The remaining panels show the estimated frames~$e[k]$, the reconstructed frame~$\hat{z}[k]$, and the reconstructed foreground~$\hat{x}[k]$ (binarized for better visualization).                              
		}
		\label{Fig:PETSSequence}
		
	\end{figure*}

	\begin{figure*}
	\centering
	
	\subfigure[Quantities related to the number of measurements.]{\label{SubFig:PETSMeasurements}
% 	\readdata{\dataCSOraclePETS}{figures/PETSCSmeasurements_oracle.dat}
%  	\readdata{\dataLLMinOracle}{figures/PETSL1L1meas_oracle.dat}
%  	\readdata{\dataLLmeasEstim}{figures/PETSL1L1meas_estim.dat}
%  	\readdata{\dataMeasurements}{figures/PETSMeasurements.dat}
 	
 	% ===================================================================
 	% Volkan Alg
 	\readdata{\dataCSOraclePETS}{figures/VolkanAlg/PETSCSmeasurements_oracle.dat}
 	\readdata{\dataLLMinOracle}{figures/VolkanAlg/PETSL1L1meas_oracle.dat}
 	\readdata{\dataLLmeasEstim}{figures/VolkanAlg/PETSL1L1meas_estim.dat}
 	\readdata{\dataMeasurements}{figures/VolkanAlg/PETSMeasurements.dat}
 	% ===================================================================
				
	\psscalebox{0.96}{
	\begin{pspicture}(8.8,5.1)
					
		\def\xMax{180}                                % Maximum value of x
		\def\xMin{0}                                  % Minimum value of x
		\def\xNumTicks{9}                             % Number of x ticks						
		\def\yMax{3500}                               % Maximum value of y
		\def\yMin{0}                                  % Minimum value of y
		\def\yNumTicks{7}                             % Number of y ticks						
		\def\xIncrement{20}                           % =(xMax-xMin)/xNumTicks						
		\def\yIncrement{500}                         % =(yMax-yMin)/yNumTicks			
					
		\def\xOrig{0.60}                              % Origin of the plot: X
		\def\yOrig{0.80}                              % Origin of the plot: Y
		\def\SizeX{8.00}                              % Size of plot: horizontal
		\def\SizeY{3.80}                              % Size of plot: vertical									
		\def\xTickIncr{0.89}                          % = SizeX/xNumTicks
		\def\yTickIncr{0.55}                          % = SizeY/yNumTicks

		\input{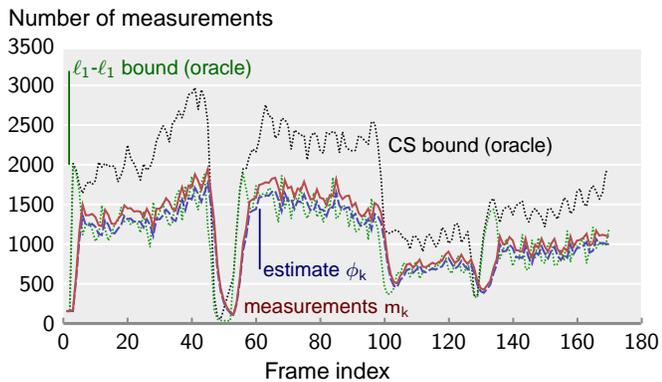}
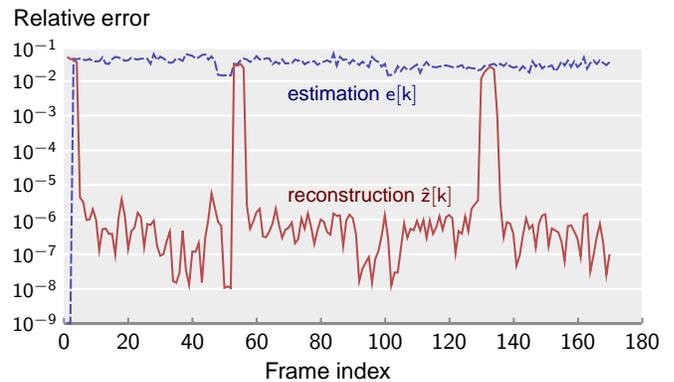
		
		% ===========================================================================================
		% Display raw data
		
		\psset{xunit=\xScale\psunit,yunit=\yScale\psunit,linewidth=0.8pt}		
		
		\dataplot[origin={\xDataOrig,\yDataOrig},showpoints=false,linestyle=dotted,dotsep=0.6pt,linecolor=black!90!white]{\dataCSOraclePETS}
 		\dataplot[origin={\xDataOrig,\yDataOrig},showpoints=false,linestyle=dotted,dotsep=0.6pt,linecolor=green!90!white]{\dataLLMinOracle}	
 		\dataplot[origin={\xDataOrig,\yDataOrig},showpoints=false,linestyle=dashed,dash=4pt 1pt,linecolor=blue!70!white]{\dataLLmeasEstim}
 		\dataplot[origin={\xDataOrig,\yDataOrig},showpoints=false,linestyle=solid,dotsep=0.8pt,linecolor=red!70!white]{\dataMeasurements}
 				
		\psset{xunit=\psunit,yunit=\psunit}
		% ===========================================================================================
		
		\rput[lb](-0.16,4.92){\small \textbf{\sf Number of measurements}}
		\rput[ct](4.25,0.14){\small \textbf{\sf Frame index}}
 			
 		\rput[lb](3.1,0.89){\footnotesize {\sf \color{red!70!black}{measurements $\mathsf{m_k}$}}}
 		
		\rput[lb](3.37,1.4){\footnotesize {\sf \color{blue!70!black}{estimate $\mathsf{\phi_k}$}}}
		\psline[linewidth=0.7pt,linecolor=blue!70!black](3.32,1.55)(3.32,2.39)
 		
 		\rput[lb](0.73,4.2){\footnotesize {\sf \color{green!70!black}{$\mathsf{\ell_1}$-$\mathsf{\ell_1}$ bound (oracle)}}}
 		\psline[linewidth=0.7pt,linecolor=green!70!black](0.68,4.3)(0.68,3.0)
 		
 		\rput[l](5.1,3.25){\footnotesize {\sf CS bound (oracle)}}
 		
		%\psgrid
	\end{pspicture}
	}
	}
	\hfill
	\subfigure[Relative errors of estimation and reconstruction.]{\label{SubFig:PETSErrors}
% 	\readdata{\estimError}{figures/PETSEstimError.dat}
% 	\readdata{\recError}{figures/PETSRecError.dat}
	% ===================================================
	% Volkan
 	\readdata{\estimError}{figures/VolkanAlg/PETSEstimError.dat}
 	\readdata{\recError}{figures/VolkanAlg/PETSRecError.dat}	
	% ===================================================
				
	\psscalebox{0.96}{
	\begin{pspicture}(8.8,5.1)
					
		\def\xMax{180}                                % Maximum value of x
		\def\xMin{0}                                  % Minimum value of x
		\def\xNumTicks{9}                             % Number of x ticks						
		\def\yMax{8}                                  % Maximum value of y
		\def\yMin{0}                                  % Minimum value of y
		\def\yNumTicks{8}                             % Number of y ticks						
		\def\xIncrement{20}                           % =(xMax-xMin)/xNumTicks						
		\def\yIncrement{1}                            % =(yMax-yMin)/yNumTicks			
					
		\def\xOrig{0.60}                              % Origin of the plot: X
		\def\yOrig{0.80}                              % Origin of the plot: Y
		\def\SizeX{8.00}                              % Size of plot: horizontal
		\def\SizeY{3.80}                              % Size of plot: vertical									
		\def\xTickIncr{0.89}                          % = SizeX/xNumTicks
		\def\yTickIncr{0.48}                          % = SizeY/yNumTicks

		%\input{auxFiles/prettyDataPlots}
		
		% ===========================================================================================
		% PrettyDataPlots code: because of log, it has to be tweaked
		
		% ==============================================================================
		% Parameters

    % Colors
    \definecolor{colorXAxis}{gray}{0.55}           % x-Axis color
    \definecolor{colorBackground}{gray}{0.93}      % Background color (0.91)

    % Distance between x-labels and the x-axis
    \def \distXLabels{0.15}

    % Distance between y-labels and the y-axis
    \def \distYLabels{0.12}

    % Width of the x-axis ticks
    \def \xTickWidth{0.08}
    % ==============================================================================

    \fpAdd{\xNumTicks}{1}{\xNumTicksPOne}         % = xNumTicks + 1
    \fpAdd{\yNumTicks}{1}{\yNumTicksPOne}

    \FPadd \xEndPoint \xOrig \SizeX                % xEndPoint = xOrig + SizeX
    \FPadd \yEndPoint \yOrig \SizeY                % yEndPoint = yOrig + SizeY

    \FPset \unit 1

    \FPsub \xRange \xMax \xMin			
    \FPdiv \xScale \SizeX \xRange 
    \FPdiv \xMultByOrigin \unit \xScale
    \FPmul \xDataOrig  \xMultByOrigin \xOrig

    \FPsub \yRange \yMax \yMin
    \FPdiv \yScale \SizeY \yRange 
    \FPdiv \yMultByOrigin \unit \yScale
    \FPmul \yDataOrig  \yMultByOrigin \yOrig

    \fpSub{\yOrig}{\distXLabels}{\xPosLabels}
    \fpSub{\xOrig}{\distYLabels}{\yPosLabels}

    \fpAdd{\yOrig}{\xTickWidth}{\xTickTop}

    % Background
    \psframe*[linecolor=colorBackground](\xOrig,\yOrig)(\xEndPoint,\yEndPoint)
                                        
    % x-labels and x-ticks: dx-position, iB-label			
    \multido{\nx=\xOrig+\xTickIncr, \iB=\xMin+\xIncrement}{\xNumTicksPOne}{	
        \rput[t](\nx,\xPosLabels){\small $\mathsf{\iB}$}
        \psline[linewidth=0.8pt,linecolor=colorXAxis](\nx,\yOrig)(\nx,\xTickTop)	
    }						

    % y-labels and lines			
    \multido{\ny=\yOrig+\yTickIncr, \nB=\yMin+\yIncrement}{\yNumTicksPOne}{	
        %\rput[r](\yPosLabels,\ny){\small $\mathsf{\nB}$}
        \psline[linecolor=white,linewidth=0.8pt](\xOrig,\ny)(\xEndPoint,\ny)
    }						

    % x-axis			
    \psline[linewidth=1.2pt,linecolor=colorXAxis]{-}(\xOrig,\yOrig)(\xEndPoint,\yOrig)
    % ===========================================================================================
		
		% Y-labels
		\rput[r](0.52,0.81){\footnotesize $\mathsf{10^{-9}}$}
		\rput[r](0.52,1.28){\footnotesize $\mathsf{10^{-8}}$}
		\rput[r](0.52,1.76){\footnotesize $\mathsf{10^{-7}}$}
		\rput[r](0.52,2.23){\footnotesize $\mathsf{10^{-6}}$}
		\rput[r](0.52,2.74){\footnotesize $\mathsf{10^{-5}}$}
		\rput[r](0.52,3.20){\footnotesize $\mathsf{10^{-4}}$}
		\rput[r](0.52,3.69){\footnotesize $\mathsf{10^{-3}}$}
		\rput[r](0.52,4.18){\footnotesize $\mathsf{10^{-2}}$}		
		\rput[r](0.52,4.62){\footnotesize $\mathsf{10^{-1}}$}
		            
		% ===========================================================================================
		% Display raw data
		
		\psset{xunit=\xScale\psunit,yunit=\yScale\psunit,linewidth=0.8pt}		
		
		\dataplot[origin={\xDataOrig,\yDataOrig},showpoints=false,linestyle=dashed,dash=4pt 1pt,linecolor=blue!70!white]{\estimError}
		
		\dataplot[origin={\xDataOrig,\yDataOrig},showpoints=false,linestyle=solid,dotsep=0.8pt,linecolor=red!70!white]{\recError}
				
		\psset{xunit=\psunit,yunit=\psunit}
		% ===========================================================================================
		
		\rput[lb](-0.1,4.92){\small \textbf{\sf Relative error}}
		\rput[ct](4.25,0.14){\small \textbf{\sf Frame index}}
			
 		\rput[lb](3.7,2.41){\footnotesize {\sf \color{red!70!black}{reconstruction $\mathsf{\hat{z}[k]}$}}}
 		\rput[lt](3.7,4.11){\footnotesize {\sf \color{blue!70!black}{estimation $\mathsf{e[k]}$}}}

 		%\psgrid
	\end{pspicture}
	}		
	}
	
	\caption{
		Results for the PETS sequence. The displayed quantities are the same as in Fig.~\ref{Fig:ResultsHall}.
  }
	\label{Fig:ResultsPETS}
	\end{figure*}

\section{Experimental results}
\label{Sec:ExperimentalResults}

We applied the scheme described in the previous section to two sequences of images: the Hall monitor sequence\footnote{Obtained from \url{http://trace.eas.asu.edu/yuv/}} and a sequence from the PETS 2009 database.\footnote{Obtained from \url{http://garrettwarnell.com/ARCS-1.0.zip}} The Hall monitor sequence has~$282$ frames\footnote{The original sequence has~$300$ frames, but we removed the first~$18$, since they contain practically no foreground.} with two people walking in an office; the top panel of Fig.~\ref{Fig:HallSequence} shows the background image and frames~$4$, $100$, $150$, and~$250$. The PETS sequence is a sequence of~$171$ frames with several people walking on a street; the top panel of Fig.~\ref{Fig:PETSSequence} shows the background image and frames~$5$, $75$, $100$, and~$170$. Since the background of both sequences is static and the foreground in each frame is sparse, we can apply our scheme to simultaneously reconstruct and perform background subtraction on each frame. The remaining panels of Figs.~\ref{Fig:HallSequence} and~\ref{Fig:PETSSequence} show the estimated frames, the reconstructed frames, and the reconstructed foregrounds, binarized for better visualization (note that the foreground pixels are dark). 

\mypar{Experimental setup}
In both sequences, we set the oversampling parameters as~$\delta := \delta_k = 0.1$, for all~$k$, and the filter parameter as~$\alpha = 0.5$. While for the Hall sequence we used the true sparsity of the first two foregrounds, i.e., $\hat{s}_1 = s_1= 417$ and~$\hat{s}_2 = s_2 = 446$, for the PETS sequence we set these input parameters to values much smaller than their true values: $10 = \hat{s}_1 \ll s_1= 194$ and~$10 = \hat{s}_2 \ll s_2 = 211$. In spite of this poor initialization, the algorithm was able to quickly adapt, as we will see. For memory reasons, we downsampled each frame of the Hall sequence to $128 \times 128$ and each frame of the PETS sequence to~$116\times 116$. We also removed camera noise from each frame, i.e., isolated pixels, by preprocessing the full sequences. For the motion estimation, we used a block size of $\gamma\times \gamma =  8 \times 8$, and a search limit of~$6$. Finally, we mention that, after computing the side information~$w[k]$ for frame~$k$, we amplified the magnitude of its components by~$30\%$. This, according to the theory in~\cite{Mota14-CSSideInfo}, improves the quality of the side information. To solve basis pursuit in the reconstruction of the first two frames we used \textsc{spgl1}~\cite{Friedlander08-ProbingParetofrontierBasisPursuit-spgl1,vanDenBerg11-SparseOptimizationWithLeastSquaresConstraints-spgl1}.\footnote{Available at \url{http://www.math.ucdavis.edu/~mpf/spgl1/}} To solve $\ell_1$-$\ell_1$ minimization problem~\eqref{Eq:L1L1Simple} in the reconstruction of the remaining frames we used \textsc{decopt}~\cite{TranDinh14-ConstrainedConvexMinimizationViaModelBasedExcessiveGap,TranDinh14-APrimalDualAlgorithmicFrameworkForConstrainedConvexMinimization}.\footnote{Available at \url{http://lions.epfl.ch/decopt/}}

\mypar{Results}
We benchmark Algorithm~\ref{Alg:AdaptiveRate} with the CS (oracle) bound given by~\eqref{Eq:ChandrasekaranBound}. Note that the prior state-of-the-art in compressive background subtraction, \cite{Warnell12-AdaptiveRateCompressiveSensingBackgroundSubtraction,Warnell14-AdaptiveRateCompressiveSensingUsingSideInformation},  requires always more measurements than the ones given by~\eqref{Eq:ChandrasekaranBound}. The results of the experiments are in \fref{Fig:ResultsHall} for the Hall sequence and in \fref{Fig:ResultsPETS} for the PETS sequence. Figs.~\ref{SubFig:HallMeasurements} and~\ref{SubFig:PETSMeasurements} show the number of measurements~$m_k$ Algorithm~\ref{Alg:AdaptiveRate} took from each frame and the corresponding estimate~$\phi_k$ of~\eqref{Eq:L1L1Bound}. These figures also show the bounds~\eqref{Eq:L1L1Bound} and~\eqref{Eq:ChandrasekaranBound} as if an oracle told us the true values of~$s_k$, $\overline{h}_k$, and~$\xi_k$. We can see that~$m_k$ and~$\phi_k$ are always below the CS bound~\eqref{Eq:ChandrasekaranBound}, except at a few frames in Fig.~\ref{SubFig:PETSMeasurements} (PETS sequence). In those frames, there is no foreground and thus the number of required measurements approaches zero. Since there are no such frames in the Hall sequence, all quantities in Fig.~\ref{SubFig:HallMeasurements} do not exhibit such large fluctuations. Fig.~\ref{SubFig:HallMeasurements} clearly shows the advantage of our algorithm with respect to using standard CS (i.e., basis pursuit)~\cite{Cevher08-CompressiveSensingForBackgroundSubtraction,Warnell12-AdaptiveRateCompressiveSensingBackgroundSubtraction,Warnell14-AdaptiveRateCompressiveSensingUsingSideInformation}, even if CS reconstruction is performed using the knowledge of the true foreground sparsity: our algorithm required an average of~$33\%$ of the measurements that standard (oracle) CS required. Recall that the performance of the prior state-of-the-art algorithm~\cite{Warnell12-AdaptiveRateCompressiveSensingBackgroundSubtraction,Warnell14-AdaptiveRateCompressiveSensingUsingSideInformation} is always above the CS bound line. In Figs.~\ref{SubFig:HallMeasurements} and~\ref{SubFig:PETSMeasurements}, the estimate~$\phi_k$ is very close to the oracle bound~\eqref{Eq:L1L1Bound} and, for most frames, the number of measurements~$m_k$ is larger than~\eqref{Eq:L1L1Bound}, even though the oversampling factor~$\delta=0.1$ is quite small. In fact, $m_k$ was smaller than~\eqref{Eq:L1L1Bound} in less than~$7\%$ (resp.\ $30\%$) of the frames for the Hall (resp.\ PETS) sequence. Yet, the corresponding frames were reconstructed with a relatively small error, as shown in Figs.~\ref{SubFig:HallErrors} and~\ref{SubFig:PETSErrors}, and the algorithm quickly adapted. 

Figs.~\ref{SubFig:HallErrors} and~\ref{SubFig:PETSErrors} show the relative errors of the estimated image~$e[k]$ and the reconstruction image~$\hat{z}[k]$, i.e., $\|e[k] - z[k]\|_2/\|z[k]\|_2$ and~$\|\hat{z}[k] - z[k]\|_2/\|z[k]\|_2$. It can be seen that the estimation errors were approximately constant, around $0.01$ for the Hall sequence and around~$0.93$ for the PETS sequence. The reconstruction error is essentially determined by the precision of the solver for~\eqref{Eq:L1L1Simple}. It varied between $3.8\times 10^{-9}$ and~$3.5\times 10^{-6}$ for the Hall sequence [Fig.~\ref{SubFig:HallErrors}]. For the PETS sequence [Fig.~\ref{SubFig:PETSErrors}], it was always below~$10^{-5}$ except at three instances, where the reconstruction error approached the estimation error. These correspond to the frames with no foreground (making the bounds in~\eqref{Eq:L1L1Bound} and~\eqref{Eq:ChandrasekaranBound} approach zero) and to the initial frames, where the number of measurements was much smaller than~\eqref{Eq:L1L1Bound}. In spite of these ``ill-conditioned'' frames, our algorithm was able to quickly adapt in the next frames, and follow the $\ell_1$-$\ell_1$ bound curve closely.
%We also point out that we do not show an (explicit) comparison with the methods in~\cite{Warnell12-AdaptiveRateCompressiveSensingBackgroundSubtraction,Warnell14-AdaptiveRateCompressiveSensingUsingSideInformation}, since they require always more measurements than the CS oracle bound.

\begin{figure*}
	\centering
	
	\subfigure[Quantities related to the number of measurements.]{\label{SubFig:HallMeasurementsNoisy}
	
	% =======================================================================
	% Volkan
	\readdata{\dataCSOracle}{figures/noisy/HallCSmeasurements_oracle.dat}
	\readdata{\dataLLMinOracle}{figures/noisy/HallL1L1meas_oracle.dat}
	\readdata{\dataLLmeasEstim}{figures/noisy/HallL1L1meas_estim.dat}
	\readdata{\dataMeasurements}{figures/noisy/HallMeasurements.dat}
	% =======================================================================
			
	\psscalebox{0.96}{
	\begin{pspicture}(8.8,5.1)
					
		\def\xMax{300}                                % Maximum value of x
		\def\xMin{0}                                  % Minimum value of x
		\def\xNumTicks{6}                             % Number of x ticks						
		\def\yMax{6000}                               % Maximum value of y
		\def\yMin{0}                                  % Minimum value of y
		\def\yNumTicks{6}                             % Number of y ticks						
		\def\xIncrement{50}                           % =(xMax-xMin)/xNumTicks						
		\def\yIncrement{1000}                         % =(yMax-yMin)/yNumTicks			
					
		\def\xOrig{0.60}                              % Origin of the plot: X
		\def\yOrig{0.80}                              % Origin of the plot: Y
		\def\SizeX{8.00}                              % Size of plot: horizontal
		\def\SizeY{3.80}                              % Size of plot: vertical									
		\def\xTickIncr{1.33}                          % = SizeX/xNumTicks
		\def\yTickIncr{0.64}                          % = SizeY/yNumTicks

		\input{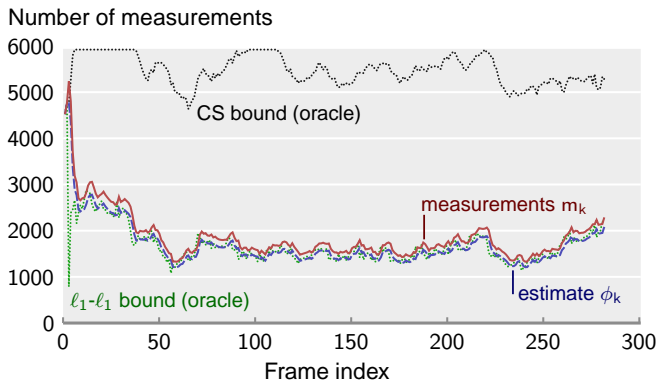}
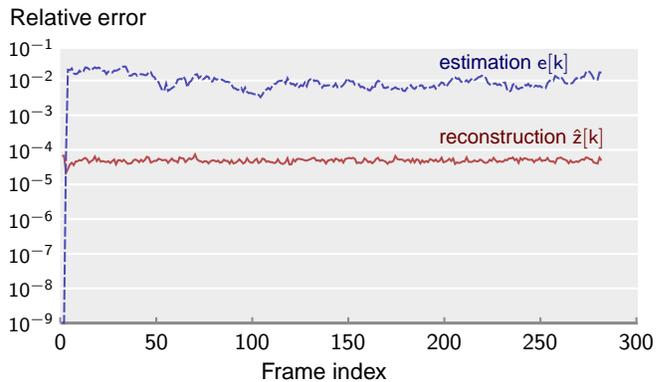
		
		% ===========================================================================================
		% Display raw data
		
		\psset{xunit=\xScale\psunit,yunit=\yScale\psunit,linewidth=0.8pt}		
		
		\dataplot[origin={\xDataOrig,\yDataOrig},showpoints=false,linestyle=dotted,dotsep=0.6pt,linecolor=black!90!white]{\dataCSOracle}
		\dataplot[origin={\xDataOrig,\yDataOrig},showpoints=false,linestyle=dotted,dotsep=0.6pt,linecolor=green!90!white]{\dataLLMinOracle}	
		\dataplot[origin={\xDataOrig,\yDataOrig},showpoints=false,linestyle=dashed,dash=4pt 1pt,linecolor=blue!70!white]{\dataLLmeasEstim}
		\dataplot[origin={\xDataOrig,\yDataOrig},showpoints=false,linestyle=solid,dotsep=0.8pt,linecolor=red!70!white]{\dataMeasurements}
				
		\psset{xunit=\psunit,yunit=\psunit}
		% ===========================================================================================
		
		\rput[lb](-0.16,4.92){\small \textbf{\sf Number of measurements}}
		\rput[ct](4.25,0.14){\small \textbf{\sf Frame index}}
			
		\rput[lt](5.56,2.55){\footnotesize {\sf \color{red!70!black}{measurements $\mathsf{m_k}$}}}
		\psline[linewidth=0.7pt,linecolor=red!70!black](5.6,1.96)(5.6,2.3)
		
		\rput[lb](6.9,1.07){\footnotesize {\sf \color{blue!70!black}{estimate $\mathsf{\phi_k}$}}}
		\psline[linewidth=0.7pt,linecolor=blue!70!black](6.83,1.20)(6.83,1.53)
		
		\rput[lb](0.7,0.97){\footnotesize {\sf \color{green!70!black}{$\mathsf{\ell_1}$-$\mathsf{\ell_1}$ bound (oracle)}}}
		
		\rput[lb](2.45,3.55){\footnotesize {\sf CS bound (oracle)}}
		%\psgrid
	\end{pspicture}
	}
	\hfill
	}
	\subfigure[Relative errors of estimation and reconstruction.]{\label{SubFig:HallErrorsNoisy}
		% ========================================================
	% Volkan
 	\readdata{\estimError}{figures/noisy/Hall_EstimError.dat}
 	\readdata{\recError}{figures/noisy/Hall_RecError.dat}	
	% ========================================================
			
	\psscalebox{0.96}{
	\begin{pspicture}(8.8,5.1)
					
		\def\xMax{300}                                % Maximum value of x
		\def\xMin{0}                                  % Minimum value of x
		\def\xNumTicks{6}                             % Number of x ticks						
		\def\yMax{8}                                  % Maximum value of y
		\def\yMin{0}                                  % Minimum value of y
		\def\yNumTicks{8}                             % Number of y ticks						
		\def\xIncrement{50}                           % =(xMax-xMin)/xNumTicks						
		\def\yIncrement{1}                            % =(yMax-yMin)/yNumTicks			
					
		\def\xOrig{0.60}                              % Origin of the plot: X
		\def\yOrig{0.80}                              % Origin of the plot: Y
		\def\SizeX{8.00}                              % Size of plot: horizontal
		\def\SizeY{3.80}                              % Size of plot: vertical									
		\def\xTickIncr{1.33}                          % = SizeX/xNumTicks
		\def\yTickIncr{0.48}                          % = SizeY/yNumTicks

		%\input{auxFiles/prettyDataPlots}
		
		% ===========================================================================================
		% PrettyDataPlots code: because of log, it has to be tweaked
		
		% ==============================================================================
		% Parameters

    % Colors
    \definecolor{colorXAxis}{gray}{0.55}           % x-Axis color
    \definecolor{colorBackground}{gray}{0.93}      % Background color (0.91)

    % Distance between x-labels and the x-axis
    \def \distXLabels{0.15}

    % Distance between y-labels and the y-axis
    \def \distYLabels{0.12}

    % Width of the x-axis ticks
    \def \xTickWidth{0.08}
    % ==============================================================================

    \fpAdd{\xNumTicks}{1}{\xNumTicksPOne}         % = xNumTicks + 1
    \fpAdd{\yNumTicks}{1}{\yNumTicksPOne}

    \FPadd \xEndPoint \xOrig \SizeX                % xEndPoint = xOrig + SizeX
    \FPadd \yEndPoint \yOrig \SizeY                % yEndPoint = yOrig + SizeY

    \FPset \unit 1

    \FPsub \xRange \xMax \xMin			
    \FPdiv \xScale \SizeX \xRange 
    \FPdiv \xMultByOrigin \unit \xScale
    \FPmul \xDataOrig  \xMultByOrigin \xOrig

    \FPsub \yRange \yMax \yMin
    \FPdiv \yScale \SizeY \yRange 
    \FPdiv \yMultByOrigin \unit \yScale
    \FPmul \yDataOrig  \yMultByOrigin \yOrig

    \fpSub{\yOrig}{\distXLabels}{\xPosLabels}
    \fpSub{\xOrig}{\distYLabels}{\yPosLabels}

    \fpAdd{\yOrig}{\xTickWidth}{\xTickTop}

    % Background
    \psframe*[linecolor=colorBackground](\xOrig,\yOrig)(\xEndPoint,\yEndPoint)
                                        
    % x-labels and x-ticks: dx-position, iB-label			
    \multido{\nx=\xOrig+\xTickIncr, \iB=\xMin+\xIncrement}{\xNumTicksPOne}{	
        \rput[t](\nx,\xPosLabels){\small $\mathsf{\iB}$}
        \psline[linewidth=0.8pt,linecolor=colorXAxis](\nx,\yOrig)(\nx,\xTickTop)	
    }						

    % y-labels and lines			
    \multido{\ny=\yOrig+\yTickIncr, \nB=\yMin+\yIncrement}{\yNumTicksPOne}{	
        %\rput[r](\yPosLabels,\ny){\small $\mathsf{\nB}$}
        \psline[linecolor=white,linewidth=0.8pt](\xOrig,\ny)(\xEndPoint,\ny)
    }						

    % x-axis			
    \psline[linewidth=1.2pt,linecolor=colorXAxis]{-}(\xOrig,\yOrig)(\xEndPoint,\yOrig)
    % ===========================================================================================
		
		% Y-labels
		\rput[r](0.52,0.81){\footnotesize $\mathsf{10^{-9}}$}
		\rput[r](0.52,1.28){\footnotesize $\mathsf{10^{-8}}$}
		\rput[r](0.52,1.76){\footnotesize $\mathsf{10^{-7}}$}
		\rput[r](0.52,2.23){\footnotesize $\mathsf{10^{-6}}$}
		\rput[r](0.52,2.74){\footnotesize $\mathsf{10^{-5}}$}
		\rput[r](0.52,3.20){\footnotesize $\mathsf{10^{-4}}$}
		\rput[r](0.52,3.69){\footnotesize $\mathsf{10^{-3}}$}
		\rput[r](0.52,4.18){\footnotesize $\mathsf{10^{-2}}$}		
		\rput[r](0.52,4.62){\footnotesize $\mathsf{10^{-1}}$}
		            
		% ===========================================================================================
		% Display raw data
		
		\psset{xunit=\xScale\psunit,yunit=\yScale\psunit,linewidth=0.8pt}		
						
		\dataplot[origin={\xDataOrig,\yDataOrig},showpoints=false,linestyle=solid,dotsep=0.8pt,linecolor=red!70!white]{\recError}
		\dataplot[origin={\xDataOrig,\yDataOrig},showpoints=false,linestyle=dashed,dash=4pt 1pt,linecolor=blue!70!white]{\estimError}
				
		\psset{xunit=\psunit,yunit=\psunit}
		% ===========================================================================================
		
		\rput[lb](-0.1,4.92){\small \textbf{\sf Relative error}}
		\rput[ct](4.25,0.14){\small \textbf{\sf Frame index}}
			 		
 		%\rput[lb](5.85,2.75){\footnotesize {\sf \color{red!70!black}{reconstruction $\mathsf{\hat{z}[k]}$}}}
 		\rput[lt](5.85,3.50){\footnotesize {\sf \color{red!70!black}{reconstruction $\mathsf{\hat{z}[k]}$}}} % Volkan
 		\rput[lb](5.85,4.25){\footnotesize {\sf \color{blue!70!black}{estimation $\mathsf{e[k]}$}}}

 		%\psgrid
	\end{pspicture}
	}		
	}
	
	\caption{
		Results for the Hall sequence for the noisy measurement case. The quantities are the same as in \fref{Fig:ResultsHall}. The CS bound curve in \text{(a)} was truncated at $6000$ measurements so that the vertical scale is the same as in \fref{SubFig:HallMeasurements}.
	}
	\label{Fig:ResultsHallNoisy}
	\end{figure*}
% =================================

\mypar{Noisy measurements}
We also applied the version of Algorithm~\ref{Alg:AdaptiveRate} that handles noisy measurements, i.e., $y[k] = A_k x[k] + \eta_k$, with~$\|\eta_k\|_2 \leq \sigma_k$, to the Hall sequence. In this case, $\eta_k$ was a vector of i.i.d.\ Gaussian entries with zero mean and variance~$4/m_k$, and we used~$\sigma_k = 2$ for all frames. The number of measurements was computed as in~\eqref{Eq:L1L1BoundNoisy} with~$\tau = 0.1$. The results are shown in \fref{Fig:ResultsHallNoisy}. It can be seen that all the quantities in \fref{SubFig:HallMeasurementsNoisy} are slightly larger than in \fref{SubFig:HallMeasurements} (we truncated the CS bound in \fref{SubFig:HallMeasurementsNoisy} so that the vertical scales are the same). Yet, all the curves have the same shape. The most noticeable difference between the noisy and the noiseless case is the reconstruction error (\fref{SubFig:HallErrorsNoisy}), which is about~$3$ orders of magnitude larger for the noisy case.

\section{Conclusions}
\label{Sec:Conclusions}

We proposed and analyzed an online algorithm for reconstructing a sparse sequence of signals from a limited number of measurements. The signals vary across time according to a nonlinear dynamical model, and the measurements are linear. Our algorithm is based on $\ell_1$-$\ell_1$ minimization and, assuming Gaussian measurement matrices, it estimates the required number of measurements to perfectly reconstruct each signal, automatically and on-the-fly. We also explored the application of our algorithm to compressive video background subtraction and tested its performance on sequences of real images. It was shown that the proposed algorithm allows reducing the number of required measurements with respect to prior compressive video background subtraction schemes by a large margin.

\appendices

\section{Proof of Lemma~\ref{lem:Delta}}
\label{Sec:AppProofLemmaDelta}

	First, note that condition~\eqref{Eq:LemDelta} is a function only of the parameters of the sequences~$\{x[k]\}$ and~$\{\epsilon[k]\}$ and, therefore, is a deterministic condition. Simple algebraic manipulation shows that it is equivalent to
		\begin{multline*}
			(1+\delta_i)
			\bigg[
				2\overline{h}_{i-1} \log\Big(\frac{n}{u_{i-1}}\Big) + \frac{7}{5}u_{i-1} + 1
			\bigg]
			\geq 2\overline{h}_i \log\Big(\frac{n}{u_i}\Big)  \\ + \frac{7}{5}u_i + 1\,,
		\end{multline*}
		or
		\begin{equation}\label{Eq:ProofLemma1EquivalentCondition}
			(1+\delta_i) \overline{m}_{i-1} \geq \overline{m}_i\,,
		\end{equation} 
		where~$\overline{m}_i$ is the right-hand side of~\eqref{Eq:L1L1Bound} applied to~$x[i]$, that is,
		\begin{equation}\label{Eq:ProofLemma1DefOverlinem}
			\overline{m}_i := 2\overline{h}_i\log\Big(\frac{n}{s_i + \xi_i/2}\Big) + \frac{7}{5}\Big(s_i + \frac{\xi_i}{2}\Big) + 1\,.
		\end{equation}
		Notice that the source of randomness in Algorithm~\ref{Alg:AdaptiveRate} is the set of matrices (random variables)~$A_k$, generated in steps~\ref{SubAlg:step2} and~\ref{SubAlg:step10}. Define the event~$S_i$ as ``perfect reconstruction at time~$i$.'' Since we assume that~$\hat{s}_1$ and~$\hat{s}_2$ are larger than the true sparsity of~$x[1]$ and~$x[2]$, there holds~\cite{Chandrasekaran12-ConvexGeometryLinearInverseProblems}
		\begin{align}
				\mathbb{P}(S_i) 
			\geq 
				1 - \exp\Big[-\frac{1}{2}(m_i - \sqrt{m_i})^2\Big] 
			\notag
			\\			
			\geq 
				1 - \exp\Big[-\frac{1}{2}(\underline{m} - \sqrt{\underline{m}})^2\Big]\,,
			\label{Eq:ProofLemma1ProbS1S2}
		\end{align}
		for~$i=1,2$, where the second inequality is due to $m_i \geq \underline{m}$ and $1-\exp(-(1/2)(x - \sqrt{x})^2)$ being an increasing function.
		
		Next, we compute the probability of the event ``perfect reconstruction at time~$i$'' given that there was ``perfect reconstruction at all previous time instants~$l<i$,'' i.e., $\mathbb{P}(S_i |\bigwedge_{l<i} S_{l})$, for all~$i=3,\ldots,k$. Since we assume~$\alpha = 1$, we have~$\phi_i = \hat{\overline{m}}_{i-1}$ and step~\ref{SubAlg:ChooseMeas} of Algorithm~\ref{Alg:AdaptiveRate} becomes~$m_i = (1+\delta_i)\hat{\overline{m}}_{i-1}$, for all~$i \geq 3$. Under the event~$S_{i-1}$, i.e., $\hat{x}[i-1] = x[i-1]$, we have $\hat{\overline{h}}_{i-1} = \overline{h}_{i-1}$, $\hat{\xi}_{i-1} = \xi_{i-1}$, and~$\hat{\overline{m}}_{i-1} = \overline{m}_{i-1}$, where $\overline{m}_{i-1}$ is defined in~\eqref{Eq:ProofLemma1DefOverlinem}. (The hat variables are random variables.) Consequently, due to our assumption~\eqref{Eq:ProofLemma1EquivalentCondition}, step~\ref{SubAlg:ChooseMeas} can be written as $m_i = (1+\delta_i)\hat{\overline{m}}_{i-1} = (1+\delta_i)\overline{m}_{i-1} \geq \overline{m}_i$. This means~\eqref{Eq:L1L1Bound} is satisfied and hence, for~$i\geq 3$,
		\begin{align}
				\mathbb{P}\Big(S_i\,\big|\, \bigwedge_{l<i} S_{l}\Big) 
			&\geq 
				1 - \exp\Big[-\frac{1}{2}(m_i - \sqrt{m_i})^2\Big]
			\notag
			\\
			&\geq
				1 - \exp\Big[-\frac{1}{2}(\underline{m} - \sqrt{\underline{m}})^2\Big]\,,
			\label{Eq:ProofLowerB}
		\end{align}
		where, again, we used the fact that $m_i \geq \underline{m}$ and that $1-\exp(-(1/2)(x - \sqrt{x})^2)$ is an increasing function.
		
		Finally, we bound the probability that there is perfect reconstruction at all time instants~$1\leq i \leq k$:
		\begin{align}
				\mathbb{P}(S_1 &\wedge S_2 \wedge \cdots \wedge S_k)
			\notag
			\\
			&=
				\mathbb{P}(S_1)\mathbb{P}(S_2|S_1)\prod_{i=3}^k\mathbb{P}(S_i|S_1 \wedge \cdots \wedge S_{i-1})
			\label{Eq:ProofLemma1LastStep1}
			\\
			&=
				\mathbb{P}(S_1)\mathbb{P}(S_2)\prod_{i=3}^k\mathbb{P}\Big(S_i\,\big|\, \bigwedge_{l< i} S_l\Big)
			\label{Eq:ProofLemma1LastStep2}
			\\
			&\geq
				\bigg(1 - \exp \Big[-\frac{1}{2}(\underline{m} - \sqrt{\underline{m}})^2\Big]\bigg)^k\,.
			\label{Eq:ProofLemma1LastStep3}
		\end{align}
		From~\eqref{Eq:ProofLemma1LastStep1} to~\eqref{Eq:ProofLemma1LastStep2} we used the independence between~$S_1$ and~$S_2$. From~\eqref{Eq:ProofLemma1LastStep2} to~\eqref{Eq:ProofLemma1LastStep3}, we used~\eqref{Eq:ProofLemma1ProbS1S2} and~\eqref{Eq:ProofLowerB}.		
		\hfill\qed

\section{Proof of Theorem~\ref{Thm:LaplacianModelingNoise}}
\label{Sec:AppProofLemmaConditions}	
	Recall the definitions of~$\xi$ and~$\overline{h}$ in~\eqref{Eq:QualityParameters}:
		\begin{align*}
				\xi 
			&= 
				\big|\{j\,:\, w_j \neq x_j^\star = 0\}\big| - \big|\{j\,:\, w_j = x_j^\star \neq 0\}\big|
			\\
				\overline{h}
			&= 
				\big|
					\{j\,:\, x_j^\star > 0, \,\,\epsilon_j > 0\} \cup \{j\,:\, x_j^\star < 0, 
						\epsilon_j < 0\}
				\big|\,,
		\end{align*}
		where we rewrote~$\overline{h}$ using~$x^\star = w + \epsilon$. Define the events $\mathcal{A} := \text{``}\,\exists_j\,:\, x_j = w_j = 0\,\text{''}$, $\mathcal{B} := \text{``}\,\overline{h} > 0\,\text{''}$, and
		$$
			\mathcal{C} := \text{``}\,m \geq 2\overline{h}\log\Big(\frac{n}{s + \xi/2}\Big) + \frac{7}{5}\Big(s + \frac{\xi}{2}\Big) + 1\,\text{,''}
		$$
		which are the assumptions of Theorem~\ref{Thm:L1L1}. In~$\mathcal{C}$, $m$ and~$n$  are deterministic, whereas~$s$, $\overline{h}$, and~$\xi$ are random variables. Then,
		\begin{align}
				&\mathbb{P}\big(\hat{x} = x^\star \big)
			%\notag
			%\\
			\geq
				\mathbb{P}\Big(\hat{x} = x^\star\,\Big|\, \mathcal{A} \wedge \mathcal{B} \wedge \mathcal{C}\Big)
				\cdot\mathbb{P}\big(\mathcal{A} \wedge \mathcal{B} \wedge \mathcal{C} \big)
			\notag
			\\
			&\geq
			  \bigg[1 - \exp\Big(-\frac{(m - \sqrt{m})^2}{2}\Big)\bigg]
			  \cdot\mathbb{P}\big(\mathcal{A}\wedge \mathcal{B} \wedge \mathcal{C}\big)\,,
			\label{Eq:ProofLemmaConditionsTheoremL1L1Step1AllTerms}
		\end{align}
		where we used Theorem~\ref{Thm:L1L1}.	The rest of the proof consists of lower bounding~$\mathbb{P}\big(\mathcal{A} \wedge \mathcal{B} \wedge \mathcal{C}\big)$.

		\mypar{Lower bound on \boldmath{$\mathbb{P}\big(\mathcal{A} \wedge \mathcal{B} \wedge \mathcal{C}\big)$}}
		Recall that~$w$ is fixed and that each component~$x_j^\star$ is determined by $x_j^\star = w_j + \epsilon_j$. Due to the continuity of the distribution of~$\epsilon$, with probability~$1$, no component $j \in \Sigma$ (i.e., $\sigma_j^2 \neq 0$) contributes to~$\xi$. When~$j \in \Sigma^c$, we have two cases:
		\begin{itemize}
			\item $j \in \Sigma^c \cap \mathcal{W}$ (i.e., $\sigma_j^2 = 0$ and $w_j \neq 0$): in this case, we have~$x_j^\star = w_j$ with probability~$1$. Hence, these components contribute to the second term of~$\xi$.
			
			\item $j \in \Sigma^c \cap \mathcal{W}^c$ (i.e., $\sigma_j^2 = 0$ and $w_j = 0$): in this case, we also have~$x_j^\star = w_j$ with probability~$1$. However, these components do not contribute to~$\xi$.
		\end{itemize}
		We conclude $\mathbb{P}\big(\mathcal{D}\big) = \mathbb{P}\big(\xi = -\big|\Sigma^c \cap \mathcal{W}\big|\big) = 1$,	where~$\mathcal{D}$ is the event ``$\xi = -\big|\Sigma^c \cap \mathcal{W}\big|$.'' From the second case above we also conclude that our assumption $\Sigma^c \cap \mathcal{W}^c \neq \emptyset$ implies $\mathbb{P}\big(\mathcal{A}\big) = 1$. We can then write
		\begin{align}
				\mathbb{P}\big(\mathcal{A} \wedge \mathcal{B} \wedge \mathcal{C}\big)
			&=
				\mathbb{P}\big(\mathcal{A}\big)\cdot\mathbb{P}\big(\mathcal{B} \wedge \mathcal{C}\,\big|\, \mathcal{A}\big)
			%\notag
			%\\
			=
				\mathbb{P}\big(\mathcal{B} \wedge \mathcal{C}\,\big|\, \mathcal{A}\big)
			\notag
			\\
			&\geq
				\mathbb{P}\big(\mathcal{B} \wedge \mathcal{C}\,\big|\, \mathcal{A},\, \mathcal{D}\big)
				\cdot
				\mathbb{P}\big(\mathcal{D}\,\big|\, \mathcal{A}\big)
			\label{Eq:ProofLemmaConditionsTheoremL1L1StepA-Justification1}
			\\
			&=
				\mathbb{P}\big(\mathcal{B} \wedge \mathcal{C}\,\big|\, \mathcal{A},\, \mathcal{D}\big)
				\cdot
				\mathbb{P}\big(\mathcal{D}\big)
			\label{Eq:ProofLemmaConditionsTheoremL1L1StepA-Justification2}
			\\
			&=
				\mathbb{P}\big(\mathcal{B} \wedge \mathcal{C}\,\big|\, \mathcal{A},\, \mathcal{D}\big)
			\label{Eq:ProofLemmaConditionsTheoremL1L1StepA-Justification3}
			\\
			&=
				\mathbb{P}\big( 0 < \overline{h} \leq \mu + t \,\big|\, \mathcal{A},\, \mathcal{D}\big)\,.
			\label{Eq:ProofLemmaConditionsTheoremL1L1StepA1}
		\end{align}
		From~\eqref{Eq:ProofLemmaConditionsTheoremL1L1StepA-Justification1} to~\eqref{Eq:ProofLemmaConditionsTheoremL1L1StepA-Justification2}, we used the fact that the events~$\mathcal{A} = \text{``}\Sigma^c \cap \mathcal{W}^c \neq \emptyset\text{''}$ and $\mathcal{D} = \text{``}\xi = -|\Sigma^c \cap \mathcal{W}|\text{''}$ are independent. This follows from the independence of the components of~$\epsilon$ and the disjointness of~$\Sigma^c \cap \mathcal{W}^c$ and~$\Sigma^c \cap \mathcal{W}$. From~\eqref{Eq:ProofLemmaConditionsTheoremL1L1StepA-Justification3} to~\eqref{Eq:ProofLemmaConditionsTheoremL1L1StepA1}, we used the fact that event~$\mathcal{C}$ conditioned on~$\mathcal{D}$ is equivalent to $\overline{h} \leq \mu + t$. To see why, note that the sparsity of~$x^\star$ is given by~$s = |\Sigma| + |\Sigma^c \cap \mathcal{W}|$; thus, given~$\mathcal{D}$, $s + \xi/2$ equals~$|\Sigma| + |\Sigma^c \cap \mathcal{W}|/2$; now, subtract the expression in assumption~\eqref{Eq:lemConditionsLaplacianModelNumMeas} from the expression that defines event~$\mathcal{C}$:
		$$
			0 \geq 2(\overline{h} - \mu - t)\log\bigg[\frac{n}{|\Sigma| + \frac{1}{2}|\Sigma^c \cap \mathcal{W}|}\bigg]\,.
		$$
		Using the fact that $n = |\Sigma| + |\Sigma^c| \geq |\Sigma| + |\Sigma^c \cap \mathcal{W}| \geq |\Sigma| + |\Sigma^c \cap \mathcal{W}|/2$, we conclude that~$\mathcal{C}$ is equivalent to the event ``$\overline{h} \leq \mu + t$.''	We now bound~\eqref{Eq:ProofLemmaConditionsTheoremL1L1StepA1} as follows:
		\begin{align}
				&\mathbb{P}\big(0 < \overline{h} \leq \mu + t\,\big|\,\mathcal{A},\,\mathcal{D}\big)
			\notag
			\\
			&\geq
				\mathbb{P}\big(0 < \overline{h} < \mu + t -1\,\big|\,\mathcal{A},\,\mathcal{D}\big)
			\notag
			\\
			&=
				1 - \mathbb{P}\big(\overline{h} \leq 0\,\big|\,\mathcal{A},\,\mathcal{D} \big) - \mathbb{P}\big(\overline{h} \geq \mu + t -1 \,\big|\,\mathcal{A},\,\mathcal{D} \big)
			\notag
			\\
			&=
				1 - \mathbb{P}\big(\overline{h} -\mu \leq -\mu\,\big|\,\mathcal{A},\,\mathcal{D}\big)
				-\mathbb{P}\big(\overline{h} - \mu \geq t -1\,\big|\,\mathcal{A},\,\mathcal{D}\big)
			\label{Eq:ProofLemmaConditionsTheoremL1L1StepB1}
			\\
			&\geq
				1 - \exp\bigg[-\frac{2\mu^2}{\big|\Sigma\big|}\bigg] - \exp\bigg[-\frac{2(t-1)^2}{\big|\Sigma\big|}\bigg]\,,
			\label{Eq:ProofLemmaConditionsTheoremL1L1StepB2}
		\end{align}
		where the last step, explained below, is due to Hoeffding's 
		%inequality~\cite{Hoeffding63-ProbabilityInequalitiesForSumsOfBoundedRandomVariables,Lugosi09-ConcentrationOfMeasureInequalities}.
		inequality~\cite{Lugosi09-ConcentrationOfMeasureInequalities}. 
		Note that once this step is proven, \eqref{Eq:ProofLemmaConditionsTheoremL1L1StepB2} together with~\eqref{Eq:ProofLemmaConditionsTheoremL1L1Step1AllTerms} and~\eqref{Eq:ProofLemmaConditionsTheoremL1L1StepA1} give~\eqref{Eq:lemConditionsLaplacianProb}, proving the theorem.
		
		\mypar{Proof of step~\eqref{Eq:ProofLemmaConditionsTheoremL1L1StepB1}-\eqref{Eq:ProofLemmaConditionsTheoremL1L1StepB2}}
		Hoeffding's inequality states that if~$\{Z_j\}_{j=1}^L$ is a sequence of independent random variables and $\mathbb{P}(0\leq Z_j \leq 1) = 1$ for all~$j$, then 
		%\cite{Hoeffding63-ProbabilityInequalitiesForSumsOfBoundedRandomVariables}, 
		\cite[Th.4]{Lugosi09-ConcentrationOfMeasureInequalities}:
		\begin{align}
				\mathbb{P}\bigg(\sum_{j=1}^L Z_j - \sum_{j=1}^L \mathbb{E}\big[Z_j\big] \geq \tau\bigg) &\leq \exp\Big[-\frac{2\tau^2}{L}\Big]
			\label{Eq:ProofLemmaConditionsTheoremL1L1Hoeffding1}
			\\
				\mathbb{P}\bigg(\sum_{j=1}^L Z_j - \sum_{j=1}^L \mathbb{E}\big[Z_j\big] \leq -\tau\bigg) &\leq \exp\Big[-\frac{2\tau^2}{L}\Big]\,,
			\label{Eq:ProofLemmaConditionsTheoremL1L1Hoeffding2}
		\end{align} 
		for any~$\tau > 0$. We apply~\eqref{Eq:ProofLemmaConditionsTheoremL1L1Hoeffding2} to the second term in~\eqref{Eq:ProofLemmaConditionsTheoremL1L1StepB1} and~\eqref{Eq:ProofLemmaConditionsTheoremL1L1Hoeffding1} to the third term. This is done by showing that $\overline{h}$ is the sum of~$|\Sigma|$ independent random variables, taking values in~$[0,1]$ with probability~$1$, and whose expected values sum to~$\mu$. Note that~$\mu > 0$ by definition, and~$t>1$ by assumption.
			
		We start by noticing that the components of~$\epsilon$ that contribute to~$\overline{h}$ are the ones for 	which~$\sigma_j^2 \neq 0$, i.e., $j \in \Sigma$ (otherwise, $\epsilon_j = 0$ with probability~$1$). Using the relation~$x_j^\star = w_j + \epsilon_j$ [cf.\ \eqref{Eq:IntroStateSpaceModel}], we then have $\overline{h} = \sum_{j \in \Sigma} Z_j$, where $Z_j$ is the indicator of the event
		\begin{equation}\label{Eq:ProofLemmaConditionsTheoremL1L1Indicator}
			\epsilon_j > \max \big\{0, \, -w_j\big\}
			\quad
			\text{or}
			\quad
			\epsilon_j < \min \big\{0, \, -w_j\big\}\,,
		\end{equation} 
		that is, $Z_j = 1$ if~\eqref{Eq:ProofLemmaConditionsTheoremL1L1Indicator} holds, and $Z_j = 0$ otherwise. By construction, $0\leq Z_j\leq 1$ for all~$j$. Furthermore, because the components of~$\epsilon$ are independent, so are the random variables~$Z_j$. All we have left to do is to show that the sum of the expected values of~$Z_j$ conditioned on~$\mathcal{A}$ and~$\mathcal{D}$ equals~$\mu$. This involves just simple integration. Let~$j \in \Sigma$. Then,
		\begin{align}
			&\mathbb{E}\Big[Z_j\,\big|\, \mathcal{A},\,\mathcal{D}\Big]
		%\notag
		%\\
		=			
			\mathbb{P}\big(Z_j = 1\,\big|\, \mathcal{A},\,\mathcal{D}\big)
		\label{Eq:ComputationOfBernoulliParameter0}
		\\			
		&=
			\mathbb{P}\big( \epsilon_j > \max \big\{0, \, -w_j\big\} \big)
		%\notag
		%\\
		%&\qquad
			+
			\mathbb{P}\big( \epsilon_j < \min \big\{0, \, -w_j\big\} \big)
		\label{Eq:ComputationOfBernoulliParameter1}
		\\
		&=
			\frac{1 + \exp\big(-\lambda_j\big|w_j\big|\big)}{2}
		\label{Eq:ComputationOfBernoulliParameter2}
		\\
		&=
			\frac{1 + \exp\big(-\sqrt{2}\big|w_j\big|/\sigma_j\big)}{2}\,.
			\label{Eq:ComputationOfBernoulliParameter3}
	\end{align}
	From~\eqref{Eq:ComputationOfBernoulliParameter0} to~\eqref{Eq:ComputationOfBernoulliParameter1}, we used the fact that the events in~\eqref{Eq:ProofLemmaConditionsTheoremL1L1Indicator} are disjoint for any~$w_j$. From~\eqref{Eq:ComputationOfBernoulliParameter1} to~\eqref{Eq:ComputationOfBernoulliParameter2}, we used the fact that~$\lambda_j$ is finite for~$j \in \Sigma$, and 
	\begin{align*}
		  \mathbb{P}\big( \epsilon_j > \max \big\{0, \, -w_j\big\} \big)
		%\\
		&=	
			\int_{\max\{-w_j, 0\}}^{+\infty} \frac{\lambda_j}{2}\exp(-\lambda_j |u|)\,du
% 		\\
% 		&=
% 			\left\{
% 			\begin{array}{ll}
% 				\int_{0}^{+\infty} \frac{\lambda_j}{2}\exp\big(-\lambda_j u\big)\,du &,\,\, w_j[k] > 0
% 				\vspace{0.3cm}
% 				\\
% 				\int_{-w_j[k]}^{+\infty} \frac{\lambda_j}{2}\exp\big(-\lambda_j u\big)\,du  &,\,\, w_j[k] < 0
% 			\end{array}
% 			\right.
		\\
		&=			
			\left\{
			\begin{array}{ll}
				\frac{1}{2} &,\,\, w_j > 0
				\vspace{0.3cm}
				\\
				\frac{1}{2}\exp\big(\lambda_j w_j\big)  &,\,\, w_j < 0
			\end{array}
			\right.
		\\
			\mathbb{P}\big( \epsilon_j < \min \big\{0, \, -w_j\big\} \big)			
		%\\
		&=
			\int_{-\infty}^{\min\{-w_j, 0\}} \frac{\lambda_j}{2}\exp(\lambda_j |u|)\,du
% 		\\
% 		&=
% 			\left\{
% 				\begin{array}{ll}
% 					\int_{-\infty}^{-w_j[k]} \frac{\lambda_j}{2}\exp\big(\lambda_j u\big)\,du &,\,\, w_j[k] > 0
% 					\vspace{0.3cm}
% 					\\
% 					\int_{-\infty}^{0} \frac{\lambda_j}{2}\exp\big(\lambda_j u\big)\,du &,\,\, w_j[k] < 0
% 				\end{array}
% 			\right.
		\\
		&=
			\left\{
				\begin{array}{ll}
					\frac{1}{2}\exp\big(-\lambda_i w_j\big) &,\,\, w_j > 0
					\vspace{0.3cm}
					\\
					\frac{1}{2}&,\,\, w_j < 0\,.
				\end{array}
			\right.
	\end{align*}
	And from~\eqref{Eq:ComputationOfBernoulliParameter2} to~\eqref{Eq:ComputationOfBernoulliParameter3} we simply replaced~$\lambda_j = \sqrt{2}/\sigma_j$.	The expected value of~$\overline{h}$ conditioned on~$\mathcal{A}$ and~$\mathcal{D}$ is then
	\begin{align*}
			\mathbb{E}\Big[\overline{h}\,\big|\, \mathcal{A},\,\mathcal{D}\Big]
		&=
			\mathbb{E}\Big[\sum_{j \in \Sigma} Z_j\,\big|\, \mathcal{A},\,\mathcal{D}\Big]
		=		
			\sum_{j \in \Sigma} \mathbb{E}\Big[ Z_j\,\big|\, \mathcal{A},\,\mathcal{D}\Big]
		\\
		&=
			\frac{1}{2}\sum_{j \in \Sigma} \Big[1 + \exp\big(-\sqrt{2}\big|w_j\big|/\sigma_j\big)\Big]
		=:
			\mu\,,
	\end{align*}
	where we used~\eqref{Eq:ComputationOfBernoulliParameter3}.

\hfill\qed

\bibliographystyle{IEEEtran}

{ %\isdraft{\singlespace}{}
\bibliography{paper}

\begin{thebibliography}{10}
\providecommand{\url}[1]{#1}
\csname url@samestyle\endcsname
\providecommand{\newblock}{\relax}
\providecommand{\bibinfo}[2]{#2}
\providecommand{\BIBentrySTDinterwordspacing}{\spaceskip=0pt\relax}
\providecommand{\BIBentryALTinterwordstretchfactor}{4}
\providecommand{\BIBentryALTinterwordspacing}{\spaceskip=\fontdimen2\font plus
\BIBentryALTinterwordstretchfactor\fontdimen3\font minus
  \fontdimen4\font\relax}
\providecommand{\BIBforeignlanguage}[2]{{%
\expandafter\ifx\csname l@#1\endcsname\relax
\typeout{** WARNING: IEEEtran.bst: No hyphenation pattern has been}%
\typeout{** loaded for the language `#1'. Using the pattern for}%
\typeout{** the default language instead.}%
\else
\language=\csname l@#1\endcsname
\fi
#2}}
\providecommand{\BIBdecl}{\relax}
\BIBdecl

\bibitem{Mota14-DynamicSparseStateEstimationUsingL1L1}
J.~{Mota}, N.~{Deligiannis}, A.~C. {Sankaranarayanan}, V.~{Cevher}, and
  M.~{Rodrigues}, ``Dynamic sparse state estimation using $\ell_1$-$\ell_1$
  minimization: Adaptive-rate measurement bounds, algorithms, and
  applications,'' 2015, to appear in \ICASSP.

\bibitem{Forsyth02-ComputerVision-AModernApproach}
D.~A. {Forsyth} and J.~{Ponce}, \emph{Computer Vision: A Modern
  Approach}.\hskip 1em plus 0.5em minus 0.4em\relax Prentice Hall, 2002.

\bibitem{Chellappa09-StatisticalMethodsAndModelsForVideoBasedTrackingModelingAndRecognition}
R.~{Chellappa}, A.~C. {Sankaranarayanan}, A.~{Veeraraghavan}, and P.~{Turaga},
  ``Statistical methods and models for video-based tracking, modeling, and
  recognition,'' \emph{Foundations and Trends in Signal Processing}, vol.~3,
  no. 1-2, pp. 1--151, 2009.

\bibitem{Herman09-HighResolutionRadarCS}
M.~{Herman} and T.~{Strohmer}, ``High-resolution radar via compressed
  sensing,'' \emph{\TSP}, vol.~57, no.~6, pp. 2275--2284, 2009.

\bibitem{Eldar14-ApplicationCSLongitudinalMRI}
L.~{Weizman}, Y.~{Eldar}, and D.~{Bashat}, ``The application of compressed
  sensing for longitudinal {MRI},'' 2014,
  \Arxiv{http://arxiv.org/abs/1407.2602}.

\bibitem{Ribeiro10-KalmanFilterInWSN}
A.~{Ribeiro}, I.~D. {Schizas}, S.~I. {Roumeliotis}, and G.~B. {Giannakis},
  ``Kalman filtering in wireless sensor networks,''
  \emph{\ControlSystemsMagazine}, vol.~30, no.~2, pp. 66--86, 2010.

\bibitem{Cevher08-CompressiveSensingForBackgroundSubtraction}
V.~{Cevher}, A.~{Sankaranarayanan}, M.~{Duarte}, D.~{Reddy}, R.~{Baraniuk}, and
  R.~{Chellappa}, ``Compressive sensing for background subtraction,'' in
  \emph{European Conf. on Computer Vision (ECCV)}, 2008.

\bibitem{Maddalena08-SelfOrganizingApproachBackgroundSubtractionVisualSurveillanceApplications}
L.~{Maddalena} and A.~{Petrosino}, ``A self-organizing approach to background
  subtraction for visual surveillance applications,'' \emph{\TIP}, vol.~17,
  no.~7, pp. 1168--1177, 2008.

\bibitem{Brutzer11-EvaluationBackgroundSubtractionTechniquesVideoSurveillance}
S.~{Brutzer}, B.~{H\"oferlin}, and G.~{Heidemann}, ``Evaluation of background
  subtraction techniques for video surveillance,'' in \emph{\CVPR}, 2011, pp.
  1937--1944.

\bibitem{Tseng02-RealTimeVideoSurveillanceForTrafficMonitoringUsingVirtualLineAnalysis}
B.~{Tseng}, C.-Y. {Lin}, and J.~{Smith}, ``Real-time video surveillance for
  traffic monitoring using virtual line analysis,'' in \emph{IEEE Inter. Conf.
  Multimedia and Expo}, 2002, pp. 541--544.

\bibitem{Cheung03-RobustTechniquesForBackgroundSubtractionUrbanTrafficVideo}
S.~{Cheung} and C.~{Kamath}, ``Robust techniques for background subtraction in
  urban traffic video,'' in \emph{Symposium on Electronic Imaging (SPIE)},
  2003, pp. 881--892.

\bibitem{Profio86-DigitalBackgroundSubtractionForFluorescenceImaging}
A.~{Profio}, O.~{Balchum}, and F.~{Carstens}, ``Digital background subtraction
  for fluorescence imaging,'' \emph{Med. Phys.}, vol.~13, no.~5, pp. 717--721,
  1986.

\bibitem{Otazo14-LowRankPlusSparseMatrixDecompositionAcceleratedDynamicMRI}
R.~{Otazo}, E.~{Cand\`es}, and D.~{Sodickson}, ``Low-rank plus sparse matrix
  decomposition for accelerated dynamic {MRI} with separation of background and
  dynamic components,'' 2014, to appear in Magnetic Resonance in Medicine.

\bibitem{Picardi04-BackgroundSubtractionTechniques-AReview}
M.~{Piccardi}, ``Background subtraction techniques: a review,'' in \emph{IEEE
  Inter. Conf. Systems, Man and Cybernetics}, 2004, pp. 3099--3104.

\bibitem{Candes11-RobustPrincipalComponentAnalysis}
E.~{Cand\`es}, X.~{Li}, Y.~{Ma}, and J.~{Wright}, ``Robust principal component
  analysis?'' \emph{\JACM}, vol.~58, no.~3, pp. 11:1--11:37, 2011.

\bibitem{Wakin06-CompressiveImagingForVideoRepresentationAndCoding}
M.~{Wakin}, J.~{Laska}, M.~{Duarte}, D.~{Baron}, S.~{Sarvotham}, D.~{Takhar},
  K.~{Kelly}, and R.~{Baraniuk}, ``Compressive imaging for video representation
  and coding,'' in \emph{Picture Coding Symposium}, 2006.

\bibitem{Sankaranarayanan12-CSMUVI}
A.~{Sankaranarayanan}, C.~{Studer}, and R.~{Baraniuk}, ``{CS-MUVI}: video
  compressive sensing for spatial multiplexing cameras,'' in \emph{Intern.
  Conf. Computation Photography}, 2012, pp. 1--10.

\bibitem{Sankaranarayanan13-CompressiveAcquisitionForLinearDynamicalSystems}
A.~{Sankaranarayanan}, P.~{Turaga}, R.~{Chellapa}, and R.~{Baranuik},
  ``Compressive acquisition of linear dynamical systems,'' \emph{SIAM J.
  Imaging Sciences}, vol.~6, no.~4, pp. 2109--2133, 2013.

\bibitem{Donoho06-CompressedSensing}
D.~{Donoho}, ``Compressed sensing,'' \emph{\TINFO}, vol.~52, no.~4, pp.
  1289--1306, 2006.

\bibitem{Candes06-RobustUncertaintyPrinciplesExactSignalReconstructionHighlyIncomplete}
E.~{Cand\`es}, J.~{Romberg}, and T.~{Tao}, ``Robust uncertainty principles:
  Exact signal reconstruction from highly incomplete frequency information,''
  \emph{\TINFO}, vol.~52, no.~2, pp. 489--509, 2006.

\bibitem{Warnell12-AdaptiveRateCompressiveSensingBackgroundSubtraction}
G.~{Warnell}, D.~{Reddy}, and R.~{Chellappa}, ``Adaptive rate compressive
  sensing for background subtraction,'' in \emph{\ICASSP}, 2012, pp.
  1477--1480.

\bibitem{Warnell14-AdaptiveRateCompressiveSensingUsingSideInformation}
G.~{Warnell}, S.~{Bhattacharya}, R.~{Chellappa}, and T.~{Basar},
  ``Adaptive-rate compressive sensing using side information,'' 2014,
  \Arxiv{http://arxiv.org/abs/1401.0583}.

\bibitem{Park09-MultiscaleFrameworkForCompressiveSensingOfVideo}
J.~{Park} and M.~{Wakin}, ``A multiscale framework for compressive sensing of
  video,'' in \emph{Picture Coding Symposium}, 2009, pp. 1--4.

\bibitem{Reddy11-P2C2}
D.~{Reddy}, A.~{Veeraraghavan}, and R.~{Chellappa}, ``{P2C2}: programmable
  pixel compressive camera for high speed imaging,'' in \emph{\CVPR}, 2011, pp.
  329--336.

\bibitem{Mota14-CSSideInfo}
J.~{Mota}, N.~{Deligiannis}, and M.~{Rodrigues}, ``Compressed sensing with
  prior information: Optimal strategies, geometry, and bounds,'' 2014,
  \Arxiv{http://arxiv.org/abs/1408.5250}.

\bibitem{Mota14-CSwithSideInfo-GlobalSIP}
------, ``Compressed sensing with side information: Geometrical interpretation
  and performance bounds,'' in \emph{\GlobalSIP}, 2014, pp. 675--679.

\bibitem{Donoho98-AtomicDecompositionBasisPursuit}
S.~{Chen}, D.~{Donoho}, and M.~{Saunders}, ``Atomic decomposition by basis
  pursuit,'' \emph{SIAM J. Sci. Comp.}, vol.~20, no.~1, pp. 33--61, 1998.

\bibitem{Berinde08-SparseRecoveryUsingSparseMatrices}
R.~{Berinde} and P.~{Indyk}, ``Sparse recovery using sparse matrices,''
  MIT-CSAIL-TR-2008-001, Tech. Rep., 2008.

\bibitem{Liutkus14-ImagingWithNature-CompressiveImagingUsingMultiplyScatteringMedium}
A.~{Liutkis}, D.~{Martina}, S.~{Popoff}, G.~{Chardon}, O.~{Katz}, G.~{Lerosey},
  S.~{Gigan}, L.~{Daudet}, and I.~{Carron}, ``Imaging with nature: Compressive
  imaging using a multiply scattering medium,'' \emph{Nature}, vol.~4, no.
  5552, pp. 1--7, 2014.

\bibitem{Kalman60-NewApproachLinearFilteringAndPredictionProblems}
R.~{Kalman}, ``A new approach to linear filtering and prediction problems,''
  \emph{Trans. ASME-J. Basic Engineering}, vol.~82, no.~D, pp. 35--45, 1960.

\bibitem{Haykin01-KalmanFilteringAndNeuralNetworks}
S.~{Haykin}, \emph{Kalman Filtering and Neural Networks}.\hskip 1em plus 0.5em
  minus 0.4em\relax John Wiley and Sons, 2001.

\bibitem{Geromel99-OptimalLinearFilteringUnderParameterUncertainty}
J.~{Geromel}, ``Optimal linear filtering under parameter uncertainty,''
  \emph{\TSP}, vol.~47, no.~1, pp. 168--175, 1999.

\bibitem{Ghaoui01-RobustFilteringDiscreteTimeSystemsBoundedNoiseParametricUncertainty}
L.~{El Ghaoui} and G.~{Calafiore}, ``Robust filtering for discrete-time systems
  with bounded noise and parametric uncertainty,'' \emph{\TAC}, vol.~46, no.~7,
  pp. 1084--1089, 2001.

\bibitem{Vaswani08-KalmanFilteredCS}
N.~{Vaswani}, ``Kalman filtered compressed sensing,'' in \emph{\ICIP}, 2008,
  pp. 893--896.

\bibitem{Vaswani09-AnalyzingKalmanFilteredCS}
------, ``Analyzing least squares and {Kalman} filtered compressed sensing,''
  in \emph{\ICASSP}, 2009, pp. 3013--3016.

\bibitem{Ziniel10-TrackingAndSmoothingTimeVaryingSparseSignalsBP}
J.~{Ziniel}, L.~{Potter}, and P.~{Schniter}, ``Tracking and smoothing of
  time-varying sparse signals via approximate belief propagation,'' in
  \emph{\Asilomar}, 2010, pp. 808--812.

\bibitem{Carmi10-MethodsForSparseSignalRecoveryUsingKalmanFiltering}
A.~{Carmi}, P.~{Gurfil}, and D.~{Kanevsky}, ``Methods for sparse signal
  recovery using {Kalman} filtering with embedded pseudo-measurement norms and
  quasi-norms,'' \emph{\TSP}, vol.~58, no.~4, pp. 2405--2409, 2010.

\bibitem{Balzano10-OnlineIdentificationAndTrackingOfSubspacesFromHighlyIncompleteInformation}
L.~{Balzano}, R.~{Nowak}, and B.~{Recht}, ``Online identification and tracking
  of subspaces from highly incomplete information,'' in \emph{\Allerton}, 2010,
  pp. 704--711.

\bibitem{Balzano14-LocalConvergenceOfAnAlgorithmForSubspaceIdentificationFromPartialData-GROUSE}
L.~{Balzano} and S.~J. {Wright}, ``Local convergence of an algorithm for
  subspace identification from partial data,'' \emph{\FoundComptAnalysis}, pp.
  1--36, 2014.

\bibitem{Chi13-PETRELS}
Y.~{Chi}, Y.~C. {Eldar}, and R.~{Calderbank}, ``{PETRELS}: Parallel subspace
  estimation and tracking by recursive least squares from partial
  observations,'' \emph{\TSP}, vol.~61, no.~23, pp. 5947--5959, 2013.

\bibitem{Charles11-SparsityPenaltiesDynamicalSystemEstimation}
A.~{Charles}, M.~{Asif}, J.~{Romberg}, and C.~{Rozell}, ``Sparsity penalties in
  dynamical system estimation,'' in \emph{IEEE Conf. Information Sciences and
  Systems}, 2011, pp. 1--6.

\bibitem{Chandrasekaran12-ConvexGeometryLinearInverseProblems}
V.~{Chandrasekaran}, B.~{Recht}, P.~{Parrilo}, and A.~{Willsky}, ``The convex
  geometry of linear inverse problems,'' \emph{\FoundComptAnalysis}, vol.~12,
  pp. 805--849, 2012.

\bibitem{Duarte08-SinglePixelImagingViaCompressiveSampling}
M.~{Duarte}, M.~{Davenport}, D.~{Takhar}, J.~{Laska}, T.~{Sun}, K.~{Kelly}, and
  R.~{Baraniuk}, ``Single-pixel imaging via compressive sampling,''
  \emph{\SigProcMagazine}, vol.~25, no.~2, pp. 83--91, 2008.

\bibitem{Kanevsky10-KalmanFilteringForCompressedSensing}
D.~{Kanevsky}, A.~{Carmi}, L.~{Horesh}, P.~{Gurfil}, B.~{Ramabhadran}, and
  T.~{Sainath}, ``Kalman filtering for compressed sensing,'' in \emph{Conf.
  Information Fusion (FUSION)}, 2010, pp. 1--8.

\bibitem{Girod05-DistributedVideoCoding}
B.~{Girod}, A.~{Aaron}, S.~{Rane}, and D.~{Rebollo-Monedero}, ``Distributed
  video coding,'' \emph{\ProcIEEE}, vol.~93, no.~1, pp. 71--83, 2005.

\bibitem{Deligiannis12-SideInformationDependentCorrelationChannelEstimationHashBasedDistributedVideoCoding}
N.~{Deligiannis}, J.~{Barbarien}, M.~{Jacobs}, A.~{Munteanu}, A.~{Skodras}, and
  P.~{Schelkens}, ``Side-information-dependent correlation channel estimation
  in hash-based distributed video coding,'' \emph{\TIP}, vol.~21, no.~4, pp.
  1934--1949, 2012.

\bibitem{natario2006extrapolating}
L.~{Nat{\'a}rio}, C.~{Brites}, J.~{Ascenso}, and F.~{Pereira}, ``Extrapolating
  side information for low-delay pixel-domain distributed video coding,'' in
  \emph{Visual Content Processing and Representation}, 2006, pp. 16--21.

\bibitem{wiegand2003overview}
T.~{Wiegand}, G.~J. {Sullivan}, G.~{Bjontegaard}, and A.~{Luthra}, ``Overview
  of the {H. 264/AVC} video coding standard,'' \emph{\CSVT}, vol.~13, no.~7,
  pp. 560--576, 2003.

\bibitem{alparone1996adaptively}
L.~{Alparone}, M.~{Barni}, F.~{Bartolini}, and V.~{Cappellini}, ``Adaptively
  weighted vector-median filters for motion-fields smoothing,'' in
  \emph{\ICASSP}, 1996, pp. 2267--2270.

\bibitem{Deligiannis14-MaximumLikelihoodLaplacianCorrelationChannelEstimationLayeredWynerZivCoding}
N.~{Deligiannis}, A.~{Munteanu}, S.~{Wang}, S.~{Cheng}, and P.~{Schelkens},
  ``Maximum likelihood {Laplacian} correlation channel estimation in layered
  {Wyner-Ziv} coding,'' \emph{\TSP}, vol.~62, no.~4, pp. 892--904, 2014.

\bibitem{fan2010transform}
X.~{Fan}, O.~C. {Au}, and N.~M. {Cheung}, ``Transform-domain adaptive
  correlation estimation {(TRACE)} for {Wyner-Ziv} video coding,''
  \emph{\CSVT}, vol.~20, no.~11, pp. 1423--1436, 2010.

\bibitem{Friedlander08-ProbingParetofrontierBasisPursuit-spgl1}
E.~{van den Berg} and M.~{Friedlander}, ``Probing the {P}areto frontier for
  basis pursuit solutions,'' \emph{SIAM J. Sci. Comput.}, vol.~31, no.~2, pp.
  890--912, 2008.

\bibitem{vanDenBerg11-SparseOptimizationWithLeastSquaresConstraints-spgl1}
------, ``Sparse optimization with least-squares constraints,'' \emph{SIAM J.
  Optim.}, vol.~21, no.~4, pp. 1201--1229, 2011.

\bibitem{TranDinh14-ConstrainedConvexMinimizationViaModelBasedExcessiveGap}
Q.~{Tran-Dinh} and V.~{Cevher}, ``Constrained convex minimization via
  model-based excessive gap,'' in \emph{\NIPS}, 2014.

\bibitem{TranDinh14-APrimalDualAlgorithmicFrameworkForConstrainedConvexMinimization}
------, ``A primal-dual algorithmic framework for constrained convex
  minimization,'' 2014, \Arxiv{http://arxiv.org/abs/1406.5403}.

\bibitem{Lugosi09-ConcentrationOfMeasureInequalities}
G.~{Lugosi}, ``Concentration-of-measure inequalities,'' 2009, lecture notes.

\end{thebibliography}
}

\end{document}